%% file: ncube2fin.tex
\def\version{December 22, 2003}

\documentclass[12pt]{article}
\usepackage[psamsfonts]{amsfonts}
\usepackage{amsmath,amssymb,amsthm}
\usepackage[dvips]{graphicx}
\usepackage{eepic}

\input def99c 
\UseSection  
\renewcommand{\to}      {\rightarrow}

\newcommand{\Pro}{{\mathbb P}_p}

\newcommand{\ben}{\begin{enumerate}}
\newcommand{\een}{\end{enumerate}}

\newcommand{\prob}[2][]   {  {  {\mathbb P}_{#1} ( #2 ) }  }

\newcommand{\probb}[2][]   {  {  {\mathbb P}_{#1} \bigl ( #2 \bigr ) }  }

\newcommand{\Cmax} {{\Ccal}_{\rm max}}

\newcommand{\shift}   {\!\!\!\!}
\newcommand{\SSS}   {\sss}

\newcounter{countC}  
\newcounter{countR}  
\newcommand{\sumtwo}[2]{\sum_{ \mbox{ \scriptsize
    $\begin{array}{c}
                        {#1} \\ {#2}
                        \end{array} $ }
    }
}

\newcommand{\torus}{\mathbb T}
\newcommand{\gr}{\mathbb G}

\newcommand{\cn}{\Omega}
\newcommand{\ver}{{\mathbb V}}

\newcommand{\Exp}{{{\mathbb E}_p}}
\newcommand{\egr}\Exp

\newcommand{\R}{\Rbold}
\newcommand{\Z}{\Zbold}

\newcommand{\conn}{\leftrightarrow}

\newcommand{\dbc}{\Leftrightarrow}
\newcommand{\ct}[1]     { \stackrel{#1}{\conn} }

\newcommand{\ctx}[1]     {\leftarrow\shift\!\xrightarrow{#1}}
\newcommand{\ndbc}      {{ \, \dbc {\hspace{-2.2ex} /} \hspace{+0.8ex} }}
\newcommand{\AND}       {\;\&\;}

\newcommand{\smallsup}[1] {{\scriptscriptstyle{({#1}})}}

\newcommand{\bigo}{O}

\title  {
        Random subgraphs of finite graphs:  \\ II. The lace expansion and the
        triangle condition
        }

\author{Christian Borgs\thanks{Microsoft Research, One Microsoft Way,
Redmond, WA 98052, USA. {\tt borgs@microsoft.com}, {\tt jchayes@microsoft.com}}
\and
Jennifer T.\ Chayes$^*$
\and
Remco van der Hofstad\thanks{Department of Mathematics and Computer Science,
Eindhoven University of Technology, P.O.\ Box  513,
5600 MB Eindhoven, The Netherlands.
{\tt rhofstad@win.tue.nl}}
\and
Gordon Slade\thanks{Department of Mathematics, University of British Columbia,
Vancouver, BC V6T 1Z2, Canada. {\tt slade@math.ubc.ca}}
\and
Joel Spencer\thanks{Department of Computer Science,
Courant Institute of Mathematical Sciences,
New York University, 251 Mercer St., New York, NY 10012, U.S.A.
{\tt spencer@cs.nyu.edu}
}}

\date\version

\begin{document}

\maketitle

\begin{abstract}
In a previous paper, we defined a version of the percolation
triangle condition that is suitable for the analysis of bond
percolation on a finite connected transitive graph, and showed
that this triangle condition implies that the percolation phase
transition has many features in common with the phase transition
on the complete graph. In this paper, we use a new and simplified
approach to the lace expansion to prove quite generally that for
finite graphs that are tori the triangle condition for percolation
is implied by a certain triangle condition for simple random walks
on the graph.

   The latter is readily verified for several graphs with vertex set
    $\{0,1,\ldots, r-1\}^n$, including the Hamming cube on an alphabet
    of $r$ letters (the $n$-cube, for $r=2$),
    the $n$-dimensional torus with nearest-neighbor
    bonds and $n$ sufficiently large, and the $n$-dimensional torus
    with $n>6$ and sufficiently spread-out (long range) bonds.
    The conclusions of our previous paper thus apply to the percolation
    phase transition for each of the above examples.
\end{abstract}


\noindent
{\bf Subject classifications}:  05C80, 60K35, 82B43

\smallskip \noindent
{\bf Keywords}:  random graph, phase transition, lace expansion, triangle
condition, percolation

\section{Introduction and results}
\label{sec-intro}

\subsection{Introduction}

The percolation phase transition on the complete graph
is well understood and forms a central part of modern graph theory
\cite{AS00,Boll01,JLR00}.
In the language of mathematical physics, the phase transition is
{\em mean-field}.  It can be expected that the percolation phase
transition on many other high-dimensional finite graphs will be similar
to that for the complete graph.  In other words, mean-field behaviour
will apply much more generally.

In a previous paper \cite{BCHSS04a}, we introduced the finite-graph
triangle condition, and proved that it is a sufficient condition
for several aspects of
the phase transition on a finite connected transitive graph to be mean-field.
This triangle condition is an adaptation of the well-known triangle
condition of Aizenman and Newman \cite{AN84} for infinite graphs.
In this paper, we verify the finite-graph
triangle condition for a class of graphs
with the structure of high-dimensional tori.  Examples include the $n$-cube,
the Hamming cube
and periodic approximations to $\Z^n$ for large $n$.

Our proof of the triangle condition is based on an adaptation of the
percolation lace expansion of Hara and Slade \cite{HS90a} from $\Z^n$
to finite tori.
We use the same expansion as \cite{HS90a}, but our proof of convergence
of the expansion is new and improved.
This is the first time that the lace expansion has been
applied in a setting where finite-size scaling plays a role.
An advance in our application of the lace expansion is that we prove
a general theorem that the percolation triangle condition on a finite
torus is a consequence of a corresponding condition for {\em random walks}
on the torus.  Thus, we
are able to verify the percolation triangle
condition for our examples by a relatively simple analysis of random
walks on these graphs.

\subsection{The triangle condition on infinite graphs}
\label{sec-tcig}

Let $\mathbb V$ be a finite or infinite set and let
$\mathbb B$ be a subset of the set
of all two-element subsets $\{x,y\}\subset\mathbb V$.
Then $\mathbb G=(\mathbb V,\mathbb B)$ is a finite
or infinite graph with vertex set
$\mathbb V$ and bond (or edge) set $\mathbb B$.
The degree of a vertex $x\in\mathbb V$ is
defined to be the number of edges containing $x$.
A bijective map
$\varphi:\mathbb V\to \mathbb V$ is called a graph-isomorphism if
$\{\varphi(x),\varphi(y)\}\in\mathbb B$ whenever
$\{x,y\}\in\mathbb B$.  We say that $\mathbb G$ is {\it transitive}
if for each pair $x,y \in \mathbb V$ there is a graph-isomorphism
$\varphi$ with $\varphi(x)=y$.
We will always assume that $\mathbb G$ is connected,
and usually assume that $\mathbb G$ is also transitive.
In the latter case, we denote the common degree of each vertex by $\Omega$.

We consider percolation on $\mathbb G$.  That is, we associate
independent Bernoulli random variables to the edges, taking the
value ``occupied'' with probability $p$ and ``vacant'' with probability
$1-p$, where $p \in [0,1]$ is a parameter.
Let $x \conn y$ denote the event that the
vertices $x$ and $y$ are connected by a path in $\mathbb G$ consisting
of occupied bonds, let $C(x) = \{y \in \mathbb V : x \conn y\}$
denote the connected cluster of $x$, and let $|C(x)|$ denote the
cardinality of the random set $C(x)$.
Let
    \eq
    \lbeq{taudef}
    \tau_p(x,y) = \Pbold_p(x \conn y)
    \en
denote the {\em two-point function} and define the {\em susceptibility} by
    \eq\lbeq{chidef}
    \chi(p) ={\mathbb E}_p|C(0)| 
    .
    \en
For many infinite graphs,
such as $\Z^n$ with $n \geq 2$, or
for a regular tree with degree at least three, there is a
$p_c = p_c(\mathbb G) \in (0,1)$ such that
    \eq
    \lbeq{pcG}
    p_c(\mathbb G) = \sup \{ p : \chi(p) < \infty\}
    =
    \inf \{ p : \Pbold_p(|C(0)|=\infty)>0\}.
    \en
Thus $\chi(p)<\infty$
and $\Pbold_p(|C(0)|=\infty)=0$ when $p<p_c$, whereas
$\chi(p)=\infty$ and $\Pbold_p(|C(0)|=\infty)>0$
if $p>p_c$.
The equality
of the infimum and supremum of \refeq{pcG} is a theorem of
\cite{AB87,Mens86}.

Percolation on a tree
is well understood \cite[Chapter~10]{Grim99}, and
infinite graphs whose percolation phase transition
is analogous to the transition on a tree are said
to exhibit \emph{mean-field} behaviour.
In 1984, Aizenman and Newman \cite{AN84} introduced the triangle condition
as a sufficient condition for mean-field behaviour.
The triangle condition is defined in terms of the {\em triangle
diagram}
    \eq
    \lbeq{tri-zd}
    \nabla_p(x,y)
    = \sum_{w,z \in \mathbb V} \tau_p(x,w) \tau_p(w,z) \tau_p(z,y),
    \en
and states that for all $x\in {\mathbb V}$
    \eq
    \lbeq{tc-zd}
    \nabla_{p_c}(x,x) < \infty.
    \en
It is predicted that the triangle condition on $\Z^n$ holds for
all $n>6$.

We write $f(p) = \Theta (g(p))$ if $|f(p)/g(p)|$ is bounded away from zero
and infinity in an appropriate limit.
Aizenman and Newman used a differential inequality for $\chi(p)$ to
show that the triangle condition
implies that
    \eq
    \lbeq{chi-zd}
    \chi(p) = \Theta ((p_c-p)^{-\gamma})
    \quad
    \mbox{uniformly in $p<p_c$,}
    \en
with $\gamma =1$,
and Nguyen \cite{Nguy87} extended this to show that
    \eq
    \frac{\Ebold_p [|C(0)|^{t+1}]}{\Ebold [|C(0)|^t]}
    =
    \Theta ( (p_c-p)^{-\Delta_{t+1}} )
    \quad
    \mbox{uniformly in $p<p_c$,}
    \en
with $\Delta_{t+1}=2$ for $t = 1,2,3,\ldots$.
Subsequently, Barsky and Aizenman \cite{BA91} showed, in
particular, that the triangle
condition also implies that the percolation probability obeys
    \eq
    \lbeq{theta-zd}
    \Pbold_p (|C(0)|=\infty) = \Theta((p-p_c)^{\hat \beta})
    \quad
    \mbox{uniformly in $p \geq p_c$,}
    \en
with $\hat \beta = 1$.

In 1990, Hara and Slade
established the triangle condition for nearest-neighbor bond percolation
on $\Z^n$ for large $n$ (it is now known that $n \geq 19$
is large enough), and for a wide class of long-range models, called
{\em spread-out} models, for $n>6$ \cite{HS90a,HS94}.  Their proof
of the triangle condition was based on the {\em lace expansion},
an adaptation of an expansion introduced
in 1985 by D.C.~Brydges and T.~Spencer
\cite{BS85} to study the self-avoiding walk in high dimensions.
Since the late 1980s, lace expansion methods
have been used to derive detailed estimates on the
critical behaviour of several models in high dimensions;
see \cite{HS94,MS93,Slad99} for reviews.  Recent extensions of
the lace expansion for percolation can be found in \cite{HHS03,HS00b}.

\subsection{The triangle condition on finite graphs}

On a finite graph, $|C(0)| \leq |\mathbb V| < \infty$.
Thus, there cannot be a phase transition characterized by
the divergence to
infinity
of the susceptibility or the existence of an infinite cluster.
Instead, the phase transition takes place
in a small window of $p$ values,
below which clusters are typically small in size and above
which a single giant cluster coexists with many relatively small clusters.
The basic example is the phase transition on the complete graph.

Let $\mathbb G$ be a
connected transitive finite graph, let $V=|\mathbb V| < \infty$
denote its number of vertices, and let $\cn$ denote the common
degree of these vertices.  The susceptibility
$\chi(p) = \Ebold_p|C(0)|$ is an increasing
function of $p$, with $\chi(0)=1$ and $\chi(1)=V$.
In \cite{BCHSS04a}, we defined the {\em critical threshold} $p_c=p_c(\gr)
=p_c(\gr;\lambda)$ to
be the unique solution to the equation
    \eq
    \lbeq{pcdef}
    \chi(p_c(\gr)) =  \lambda V^{1/3},
    \en
where $\lambda$ is a fixed small parameter.
As discussed in more detail in \cite{BCHSS04a}, the power $V^{1/3}$
in \refeq{pcdef} is inspired by the fact that on the complete graph
the critical susceptibility is proportional to $V^{1/3}$, and
we expect \refeq{pcdef}
to be the correct definition only for high-dimensional graphs.
The flexibility in the choice of $\lambda$ in
\refeq{pcdef} is connected with the fact that the phase transition in a finite
system is smeared out over a window rather than occurring at a sharply
defined threshold, and any value in the window could be
chosen as a threshold.

On a finite graph, the triangle diagram
\refeq{tri-zd} is bounded above by $V^2$,
and thus \refeq{tc-zd} is satisfied trivially.
In \cite{BCHSS04a}, we defined the triangle condition for a finite graph
to be the statement that
    \eq
    \lbeq{tcqn}
    \nabla_{p_c(\gr)}(x,y) \leq \delta_{x,y} + a_0,
    \en
where $a_0$ is sufficiently small.  In particular, \refeq{tcqn} implies
that $\nabla_{p_c(\gr)}(x,y)$ is uniformly bounded as $V \to \infty$.
In addition, we defined the {\em stronger} triangle condition to be the
statement that there are constants $K_1$, $K_2$ such that
for $p \leq p_c(\gr)$
\eq
\lbeq{stc}
    \nabla_p(x,y) \leq \delta_{x,y} + K_1 \cn^{-1} + K_2
    \frac{\chi(p)^3}{V}.
\en
Note that \refeq{tcqn} is a consequence of \refeq{stc}, provided
$\cn$ is sufficiently large and $\lambda$ is sufficiently small.
Moreover, since $\sum_y \nabla_{p_c}(x,y) = \chi(p_c)^3 = \lambda^3 V$,
\refeq{tcqn} implies that $\lambda^3 \leq V^{-1}+a_0$ and
hence $\lambda$ must be taken to be small for the triangle condition to hold.

As described in more
detail below, we showed in \cite{BCHSS04a} that the triangle condition
\refeq{tcqn} implies that the percolation phase transition on a finite
graph shares many features
with the transition on the complete
graph.
In this paper, we prove \refeq{stc} and hence \refeq{tcqn} for several
finite graphs, assuming that $\lambda$ is a sufficiently small constant.
These graphs all have vertex set
$\mathbb V = \{0,1,\ldots, r-1\}^n$
for some $r \geq 2$ and $n\geq 1$, with periodic boundary conditions.
We consider various edge
sets, as follows.

\subsection{Periodic tori}
\label{sec-ptori}

There are three levels of generality that we will
use.
First,
we use $\gr$ to denote a finite connected graph, which in general
need not be transitive nor regular.
Our derivation of the lace expansion, and much
of the diagrammatic
estimation of the lace expansion, is valid for general $\gr$.
Second, for our analysis of the lace expansion,
we restrict $\gr$ to have the vertex set of
the torus $\torus = \torus_{r,n} = (\Z_r)^n$, where $\Z_r$ denotes the integers
modulo $r$, for $r = 2,3, \ldots$.  The torus $\torus_{r,n}$ is an additive
group under coordinate-wise addition modulo $r$,
with volume $V=r^n$.
We allow any edge set that respects the symmetries of
translation and $x \mapsto -x$ reflections.
That is, we assume that the edge set is such that
$\{0,x\}$ is an edge if and only if $\{y,y \pm x\}$
is an edge for every vertex $y$.  For the torus (or for any
regular $\gr$) we denote the vertex degree by $\cn$.
Third,
we will verify
the stronger percolation triangle condition \refeq{stc} for the following
specific edge sets:
\begin{enumerate}
\item The nearest-neighbor torus:
an edge joins vertices that differ by $1$ (modulo $r$)
in exactly one component.
For $r=2$, this is
the $n$-cube.
For $n$ fixed and $r$ large, this is a periodic approximation to $\Z^n$.
Here $\Omega =2n$ for $r \geq 3$ and $\Omega = n$ for $r=2$.
We study the limit in which $V=r^n \to \infty$, in any fashion, provided
that $n \geq 7$ and $r \geq 2$.
\item The Hamming torus:
an edge joins vertices that differ in exactly one component
(modulo $r$).
Here $\Omega = (r-1)n$.  For $r=2$, this is again the $n$-cube.
We study the limit in which $V=r^n \to \infty$, in any fashion, provided
that $n \geq 1$ and $r \geq 2$.
\item The spread-out torus:
an edge joins vertices $x=(x_1,\ldots,x_n)$ and $y=(y_1,\ldots,y_n)$
if $0<\max_{i=1,\ldots,n}|x_i-y_i| \leq L$ (with $|\cdot|$ the metric
on $\Zbold_r$).  We study the limit $r \to \infty$,
with $n \geq 7$ fixed and $L$ large (depending on $n$)
and fixed.  This gives a periodic approximation to
range-$L$ percolation on $\Z^n$.
Here $\Omega = [(2L+1)^n-1]$
provided that $r\geq 2L+1$, which we will always assume.
\end{enumerate}

\subsection{Fourier analysis on a torus}
\label{sec-FTgroups}

Our method relies heavily on Fourier analysis.  We denote the
Fourier dual of the torus $\torus_{r,n}$ by
$\torus_{r,n}^* = \frac{2\pi}{r}\torus_{r,n}$.
We will {\em always} identify the dual torus
as $\torus_{r,n}^*
= \frac{2\pi}{r}\{-\lfloor \frac{r-1}{2}\rfloor,\ldots,\lceil
\frac{r-1}{2}\rceil \}^n$, so that each component of $k \in \torus_{r,n}^*$
is between $-\pi$ and $\pi$. The reason for this identification is that the
point $k=0$ plays a
special role, and we do not want to see it mirrored at the point $(2\pi,\ldots,
2\pi)$.
Let $k \cdot x = \sum_{j=1}^n k_jx_j$ denote the dot product
of $k \in \torus_{r,n}^*$ with $x \in \torus_{r,n}$.
The Fourier transform of $f: \torus_{r,n} \to {\mathbb C}$
is defined by
    \eq
    \lbeq{hatft}
    \hat{f}(k) = \sum_{x \in \torus_{r,n}}f(x) e^{ik \cdot x}
    \quad \quad
    (k \in \torus^*_{r,n}),
    \en
with the inverse Fourier transform given by
    \eq
    \lbeq{invftt}
    f(x)
    = \frac{1}{V} \sum_{k \in \torus_{r,n}^*} \hat{f}(k) e^{-ik \cdot x}.
    \en
The {\it convolution} of functions $f,g$ on $\torus_{r,n}$ is defined by
    \eq
    \lbeq{def-conv}
    (f*g)(x) = \sum_{y\in \torus_{r,n}} f(y)g(x-y),
    \en
and the Fourier transform of a
convolution is the product of the Fourier
transforms:
    \eq
    \lbeq{ftconv}
    \widehat{f*g}  = \hat f  \hat g.
    \en

\subsection{The triangle diagram in Fourier form}

It is convenient to use translation invariance to
regard the two-point function or triangle diagram
as a function of a single variable, e.g.,  $\tau_p(x,y)=\tau_p(y-x)$.
With this identification,
\eq
    \hat{\tau}_p(k) = \sum_{x \in \torus_{r,n}} \tau_p(0,x) e^{ik \cdot x},
\en
where $0$ denotes the origin of $\torus_{r,n}$.
It is shown in \cite{AN84} that $\hat\tau_p(k) \geq 0$ for all $k\in
\torus_{r,n}^*$.
The expected cluster size and two-point function are related by
    \eq
    \lbeq{chitau}
    \chi(p) =\Ebold_p|C(0)|
    = \sum_{x\in \torus_{r,n}} \Ebold_p I[x \in C(0)]
    = \sum_{x \in \torus_{r,n}} \tau_p(0,x) = \hat{\tau}_p(0),
    \en
where $I[E]$ denotes the indicator function of the event $E$.
In particular, writing $p_c=p_c(\torus_{r,n})$,
\eq
\lbeq{tauhat0}
    \hat\tau_{p_c}(0) = \chi(p_c) = \lambda V^{1/3}.
\en
Recalling \refeq{def-conv}, the triangle diagram
\refeq{tri-zd} can be written as
    \eq
    \nabla_p(x,y) = (\tau_p *\tau_p * \tau_p)(y-x).
    \en
By \refeq{ftconv} and \refeq{invftt}, this implies that
$\hat{\nabla}_p(k) = \hat{\tau}_p(k)^3$ and
    \eq
    \lbeq{tri-k}
    \nabla_p(x,y) = \frac{1}{V} \sum_{k \in \torus_{r,n}^*}
    \hat{\nabla}_p(k)
    e^{-ik \cdot (y-x)}
    = \frac{1}{V} \sum_{k \in \torus_{r,n}^*}
     \hat{\tau}_p(k)^3 e^{-ik \cdot (y-x)}.
    \en
By \refeq{tauhat0}, when $p=p_c$
the contribution to the right side of \refeq{tri-k}
due to the term $k=0$ is $V^{-1} \lambda^3 V
    = \lambda^3$.
This shows a connection between the
definition $\chi(p_c) = \lambda V^{1/3}$ and the triangle condition,
which in turn is connected to mean-field behavior.



\subsection{Main results}
\label{sec-mainresults}

\subsubsection{The random walk triangle condition}

For $x,y \in \torus_{r,n}$, let
    \eq
    \lbeq{Ddef}
    D(x,y) = D(y-x) = \frac{1}{\Omega} I[ \{x,y\} \in \mathbb B],
    \en
where $\mathbb B$ denotes a particular choice of edge set for the torus.
Thus, $D(x)$ represents the 1-step transition probability for a random walk
to step from $0$ to a neighbor $x$.
As in Section~\ref{sec-ptori}, we assume that
$\mathbb B$ is symmetric in the sense that
$\{0,x\} \in \mathbb B$ if and only if $\{y,y\pm x\} \in \mathbb B$
for every vertex $y$.
We make the following assumptions on $D$,
which can alternatively be regarded as
assumptions on the edge set $\mathbb B$.



\begin{ass}
\label{ass-rw}
There exists $\beta>0$ such that
    \eq
    \lbeq{supbds}
    \max_{x\in \torus_{r,n}} D(x) \leq \beta
    \en
and
    \eq
    \lbeq{rwbd}
    \frac{1}{V}
    \sum_{k \in \torus_{r,n}^*: k\neq 0}
    \frac{\hat{D}(k)^{2}}{[1-\hat{D}(k)]^3} \leq
    \beta .
    \en
\end{ass}

The assumption \refeq{supbds} is
straightforward.  As we will
discuss in more detail in Section~\ref{sec-rw}, the critical two-point function
for random walks is $[1-\hat D(k)]^{-1}$, and comparing with the right side
of \refeq{tri-k}, the assumption \refeq{rwbd} can be
interpreted as a kind of generalized triangle condition for random walks.
Note that the omitted term in \refeq{rwbd}, with $k=0$, is infinite.
For any $D$ defined by \refeq{Ddef},
\refeq{supbds} implies that $\beta\geq \cn^{-1} \geq V^{-1}$.
We will require below that $\beta$ be small.
In particular, the degree of the graph must be large.

Random walks on each of the three tori listed in Section~\ref{sec-ptori}
obey Assumption~\ref{ass-rw} with $\beta$ proportional to $\cn^{-1}$,
as the following proposition shows.
The proof of the proposition is given in Section~\ref{sec-rw}.


\begin{prop}\label{lem-Cbd}
There is an $a >0$ such that
random walks on each of the three tori listed in Section~\ref{sec-ptori}
obey Assumption~\ref{ass-rw} with $\beta = a\cn^{-1}$,
where:
\begin{enumerate}
\item for the nearest-neighbor torus, $a$ is a universal constant,
independent of $r\geq 2$ and $n \geq 7$;
\item for the Hamming torus, $a$ is a universal constant,
independent of $r\geq 2$ and $n \geq 1$;
\item for the spread-out torus, $n\geq 7$ is fixed, $r$ is sufficiently large
depending on $L$ and $n$, and $a$ depends on $n$ but not on $L$ or $r$.
\end{enumerate}
\end{prop}

\subsubsection{The triangle condition and its consequences}

Our main result is that if Assumption~\ref{ass-rw} holds
with appropriately small parameters, then the percolation
triangle condition holds.  By Proposition~\ref{lem-Cbd}, this
establishes the triangle condition for the three tori listed in
Section~\ref{sec-ptori}.

\begin{theorem}[The triangle condition]
\label{thm-tc}
Consider the torus $\torus_{r,n}$ with edge set such that
$\{0,x\}$ is an edge if and only if $\{y,y\pm x\}$ is an edge
for any vertex $y$.
There is an absolute constant $\beta_0>0$,
not depending on $r$, $n$ or the edge set of $\torus_{r,n}$,
such that the stronger triangle condition \refeq{stc} holds in the form
    \eq
    \lbeq{triacon}
    \nabla_{p}(x,y)
    \leq \delta_{x,y} + 13 \beta
    + 10 \frac{\chi(p)^3}{V},
    \en
whenever $\lambda^3 \leq \beta_0$, $p\leq p_c$
and Assumption~\ref{ass-rw} holds with $\beta\leq \beta_0$.
\end{theorem}

This establishes \refeq{stc} for our three tori, since $\beta$ is proportional
to $\cn^{-1}$ in Proposition~\ref{lem-Cbd}.
In particular, the cases covered include:

\begin{itemize}
\item
the $n$-cube $\torus_{2,n}$,
\item the complete graph (Hamming torus with $n=1$ and $r \to \infty$),
\item
nearest-neighbor percolation on $\torus_{r,n}$ with $n \geq 7$
and $r^n \to \infty$ in any fashion, including $n$ fixed and $r \to \infty$,
$r$ fixed and $n \to \infty$, or $r,n \to \infty$ simultaneously,
\item periodic approximations to range-$L$ percolation on $\Z^n$ for
fixed $n \geq 7$ and fixed large $L$.
\end{itemize}

It follows that the various consequences of
the triangle condition established in \cite{BCHSS04a} hold for these
three tori,
provided $\lambda$ and $a\cn^{-1}$ are sufficiently small (as
required by the smallness of the triangle), and $\lambda
V^{1/3}$ is sufficiently large (as required by the additional
condition on $\lambda V^{1/3}$ from
\cite[Theorems~1.2--1.4]{BCHSS04a}).
Note that if $\lambda$ is a fixed positive constant then
the last
condition merely states that $V$ is large. We now summarize these
consequences in this context.  To this end, it will be convenient
to use the standard $O(\cdot)$ notation.  All constants
implicitly in these $O$-symbols are independent of the parameters
of the model, except for an implicit dependence through the constant
$a$ from Proposition~\ref{lem-Cbd}.

The asymptotic behaviour of the critical value $p_c$ is given in
\cite[Theorem~1.5]{BCHSS04a} as follows.

\begin{theorem}[Critical threshold]
\label{thm-pcasy}
For the three tori,
    \eq
    p_c = \frac{1}{\cn}\big[ 1+  \bigo(\cn^{-1})
    + \bigo(\lambda^{-1}V^{-1/3})\big].
   \lbeq{pcasy}
    \en
\end{theorem}



For the subcritical phase, the following results are consequences of
\cite[Theorems~1.2, 1.5]{BCHSS04a}.   A version of
\refeq{chibd} valid for all $p \leq p_c$
is given in \cite[Theorem~1.5]{BCHSS04a}; see
\refeq{chibda} below.
Let $\Cmax$ denote a cluster of
maximal size, and let
    \eq\lbeq{Cmaxdef}
    |\Cmax| = \max\{|C(x)| :  x \in {\mathbb V} \}.
    \en

\begin{theorem}[Subcritical phase]
\label{main-thm-sub}
Let $p=p_c-\cn^{-1}\epsilon$ with $\epsilon\geq 0$.  For the three tori,
the following hold.

\noindent
{\rm i)}  If
$\epsilon \lambda V^{1/3} \to \infty$
as $V \to \infty$, then as
$V \to \infty$,
    \eq\lbeq{chibd}
    \chi(p)
    =\frac 1\epsilon\bigl[1+O(\cn^{-1})+O((\epsilon\lambda V^{1/3})^{-1})\bigr].
    \en

\noindent
{\rm ii)} For all $\epsilon \geq 0$,
    \eq\lbeq{cmaxbd1}
    10^{-4} \chi(p)^2
    \leq
    \Exp\Big(|\Cmax|\Big)
    \leq
    2\chi(p)^2\log(V/\chi(p)^3),
    \en
    \eq\lbeq{cmaxbd2}
    \Pro\Big( 
    |\Cmax|\leq 2\chi(p)^2\log(V/\chi(p)^3)\Big)
    \geq 1-
    \frac {\sqrt{e}}{[2\log(V/\chi(p)^3)]^{3/2}},
\en
and, for $\omega \geq 1$,
\eq
\lbeq{cmaxbdom}
    \Pro\Big(|\Cmax|\geq \frac{\chi(p)^2}{3600\omega}\Big)
    \geq\big(1+\frac {36\chi(p)^3}{\omega V}\Big)^{-1}.
\en
\end{theorem}

Inside a scaling window of width proportional to $V^{-1/3}$,
the following results are consequences of
\cite[Theorem~1.3]{BCHSS04a}.

\begin{theorem}[Scaling Window]
\label{main-thm-critical}
Fix $\lambda>0$ sufficiently small and $\Lambda<\infty$.
For the three tori,
there exist constants $b_1, \ldots, b_8$ such that the following hold
for all
$p=p_c +\cn^{-1}\epsilon$ with $|\epsilon|\leq\Lambda V^{-1/3}$.

\noindent
{\rm i)} If
$k\leq b_1 V^{2/3}$, then
\eq\lbeq{clszdis}
\frac{b_2}{\sqrt k}
\leq
\Pbold_p(|C(0)| \geq k)
\leq
\frac{b_3}{\sqrt k}.
\en

\noindent
{\rm ii)} 
\eq\lbeq{LCEBd1-win}
    {b_4}V^{2/3}
    \leq
    \Exp\big[|\Cmax|\big]
    \leq
    {b_5}V^{2/3}
    \en
and, if $\omega\geq 1$, then
    \eq\lbeq{LCBd1-win}
    \Pro\Big(
    \omega^{-1} V^{2/3}\leq |\Cmax|\leq \omega V^{2/3}
        \Big)
    \geq 1-\frac{b_6}\omega.
    \en

\noindent
{\rm iii)}
    \eq\lbeq{chiasy-win}
    b_7 V^{1/3}
    \leq \chi(p)\leq
    b_8 V^{1/3}.
    \en

\noindent In the above statements, the constants $b_2$ and $b_3$
can be chosen independent of $\lambda$ and $\Lambda$, the constants
$b_5$ and $b_8$ depend on $\Lambda$ and not on $\lambda$, and the constants
$b_1$, $b_4$, $b_6$ and $b_7$ depend on both $\lambda$ and $\Lambda$.

\end{theorem}

For the supercritical phase, the following results are consequences of
\cite[Theorem~1.4]{BCHSS04a}.

\begin{theorem} [Supercritical phase]
\label{main-thm-sup}
Let $p=p_c+\epsilon \cn^{-1}$ with $\epsilon \geq 0$.
For the three tori,
\smallskip

\noindent
{\rm i)}
        \eq\lbeq{bound1-cmax}
        \Exp(|\Cmax|)
        \leq  21\epsilon V+ 7V^{2/3},
        \en
and, for all $\omega >0$,
    \eq
    \lbeq{cmax.2A}
    \Pro\Big(|\Cmax|\leq \omega (V^{2/3}+\epsilon V)\Big)
    \geq 1-\frac{21}\omega.
    \en

\noindent
{\rm ii)}
    \eq\lbeq{chiasysup}
    \chi(p)
    \leq 81(V^{1/3}+\epsilon^2V).
    \en
\end{theorem}

Theorem~\ref{main-thm-sup} provides upper bounds on the size of
clusters in the supercritical phase.  To see that a phase transition
occurs at $p_c$, one wants a {\em lower} bound.  We have not proved
a lower bound at the level of generality of all three tori, but we
have obtained a lower bound for the case of the $n$-cube $\torus_{2,n}
=\{0,1\}^n$.
This is the content of the following theorem, which is proved in
\cite[Theorem~1.5]{BCHSS04c}.  The statement that $E_n$ occurs a.a.s.\
means that $\lim_{n\to \infty} \Pbold (E_n)=0$, assuming that
$\lambda$ is fixed as $n\to\infty$.

\begin{theorem}[Supercritical phase for the $n$-cube]
\label{thm-qnsuper}
There are strictly positive constants $c_0$, $c_1$, $c_2$ such that
the following holds
for
$\gr =\torus_{2,n}$,
all $n$-independent $\lambda$ with $0<\lambda\leq c_0$
and all $p = p_c+\epsilon n^{-1}$ with
$e^{-c_1n^{1/3}} \leq \epsilon\leq 1$:
    \eqalign
    \lbeq{mr1}
    |\Cmax| & \geq c_2 \epsilon 2^n
    \quad\text{ a.a.s. as }n\rightarrow \infty,
    \\
    \lbeq{mr2}
    \chi(p)
    & \geq (c_2 \epsilon)^2 2^n
    \quad\text{ as }n\rightarrow \infty.
    \enalign
\end{theorem}

For the special case of the $n$-cube, Theorem~\ref{thm-pcasy} states
that if $\lambda$ is chosen such that $\lambda^{-1}2^{-n/3} = O(n^{-1})$,
then
$p_c(n) = \frac 1n + O( \frac {1}{n^2})$.  This result has been
extended in \cite{HS04b} to show that there are rational numbers
$a_i$ ($i \geq 1$) such that for all positive integers $M$, all $c,c'>0$,
and all $p$ for which $\chi(p) \in [cn^M,c'n^{-2M}2^n]$,
\eq
    p = \sum_{i=1}^M a_i n^{-i} +O(n^{-M-1}),
\en
where the constant in the error term depends only on $c,c',M$.  It follows from
Theorem~\ref{thm-pcasy} that $a_1=1$, and it is shown in
\cite{HS04a} that $a_2=1$ and $a_3= \frac 72$.

\subsection{Discussion}
\label{sec-H&C}

\subsubsection{Restriction to high-dimensional tori}

Our results show that the phase transition for
percolation on general graphs obeying the triangle condition
shares several features
with the phase transition for the complete graph.
This mean-field behavior 
is expected to apply only
to graphs that are in some sense high-dimensional, and our entire
approach is restricted to high-dimensional graphs.  As discussed
in \cite[Section 3.4.2]{BCHSS04a}, we do not expect the definition
\refeq{pcdef} of the critical threshold to be correct for finite
approximations to low-dimensional graphs, such as $\Z^n$ for $n<6$.
Neither do we expect the triangle condition to be relevant in low
dimensions.

Since every finite abelian group is a direct product of cylic groups,
our restriction to the torus actually covers all abelian groups, apart
from the fact that we consider constant widths in all directions
and make a symmetry assumption.
It would be straightforward to generalize our results to tori with
different widths in different directions.
This leaves open the case of  more general graphs
and non-abelian groups,
which would require a replacement for
both the $x \mapsto -x$ symmetry of the torus
and the commutative law.

\subsubsection{The lace expansion}

The derivation of the lace
expansion in \cite{HS90a} applies immediately to finite
graphs, which need not be transitive nor regular.
Our proof of convergence of the expansion uses
the group structure of the torus for Fourier analysis, as well as the
$x \mapsto -x$ symmetry of the torus.
The proof is an adaptation of the original convergence proof of \cite{HS90a},
but is conceptually simpler and the idea of basing the proof on
Assumption~\ref{ass-rw} is new.  In addition, we
benefit from working on a finite set where Fourier integrals are
simply finite sums.

\subsubsection{Bulk versus periodic boundary conditions}

A natural question for $\Z^n$ is the following.  For $p=p_c(\Z^n)$,
consider the restriction of percolation configurations to a large
box of side $r$, centered at the origin.  How large is the largest
cluster in the box, as $r \to \infty$?  The combined results of
Aizenman \cite{Aize97} and Hara, van der Hofstad and Slade \cite{HHS03}
show that for spread-out models with $n>6$
the largest cluster has size of order $r^4$,
and there are order $r^{d-6}$ clusters of this size.
For the nearest-neighbor model in dimensions $n \gg 6$, the same results
follow from the combined results of \cite{Aize97} and Hara \cite{Hara00}.
These results apply under the {\em bulk}\/ boundary condition, in which
the clusters in the box are defined to be the intersection of the box with
clusters in the infinite lattice (and thus clusters in the box need not be
connected within the box).
In terms of the volume $V=r^n$ of the box, the largest cluster at $p_c(\Z^n)$
therefore
has size $V^{4/n}$,  for $n>6$.
Aizenman \cite{Aize97} raised
the interesting question whether
the $r^4=V^{4/n}$ would change
to $r^{2n/3}=V^{2/3}$ if the periodic boundary condition is used instead
of the the bulk boundary condition.

Theorem~\ref{main-thm-critical} shows that for $p$ within a scaling window of
width proportional to $V^{-1/3}$, centered at
$p_c(\torus_{r,n})$, the largest cluster is of size $V^{2/3}$ both
for the sufficiently spread-out model with
$n>6$ and the nearest-neighbor model with $n$ sufficiently large.
An affirmative answer to
Aizenman's question would then follow if we
could prove that $p_c(\Z^n)$ is
within this scaling window.   It would be interesting to investigate this
further.

\subsection{Organization}
\label{sec-org}

The remainder of this paper is organized as follows.
In Section~\ref{sec-rw}, we analyze random walks on a torus and
verify Assumption~\ref{ass-rw} for the three tori listed in
Section~\ref{sec-ptori}.
In Section~\ref{sec-lexp}, we give a self-contained
derivation of the lace expansion.
In Section~\ref{sec-diagest}, we estimate the Feynman diagrams that
arise in the lace expansion.  The results of
Sections~\ref{sec-lexp} apply on an arbitrary
finite graph $\gr$, which need not be transitive
nor regular.
Parts of Section~\ref{sec-diagest} also apply in this general
context, but in Section~\ref{sec-dbd} we will
specialize to $\torus_{r,n}$.
In Section~\ref{sec-Pibds}, we analyze the lace expansion
on an arbitrary torus that obeys Assumption~\ref{ass-rw},
thereby proving Theorem~\ref{thm-tc}.
Finally, in Section~\ref{sec-pfstPibds}, we establish
a detailed relation between the Fourier transforms
of the two-point functions for percolation and random walks.

\section{Proof of Proposition~\ref{lem-Cbd}}
\label{sec-rw}

\subsection{The random walk two-point function}

Consider a random walk on $\torus_{r,n}$ where the transition
probability for a step from $x$ to $y$ is equal to $D(x,y)$, with $D$ given
by \refeq{Ddef}.  We assume that the edge set of the torus is invariant
under translations  and $x \mapsto -x$ reflections.
The {\em two-point function}\/ for the random walk is defined by
    \eq
    \lbeq{SRW}
    C_{\mu}(0,x) = \sum_{\omega:0\to x} \mu^{|\omega|},
    \en
where $0 \leq \mu < \cn^{-1}$, the sum is over all
random walks $\omega$
from $0$ to $x$ that take any number of steps $|\omega|$, and the
``zero-step'' walk contributes $\delta_{0,x}$.
This is well-defined, because the fact that there are $\cn^m$ nearest-neighbor
random walks of length $m$ starting from the origin implies that
\eq
    C_\mu(0,x) \leq \sum_{x \in \torus_{r,n}}C_\mu(0,x)
    = \sum_{m=0}^\infty \cn^m \mu^m
    = \frac{1}{1-\mu\cn}
    \quad \quad
    (\mu < \cn^{-1}),
\en
i.e., the random walk susceptibility $\sum_{x \in \torus_{r,n}}C_\mu(0,x)$
is finite.
Probabilistically, $C_{1/\cn}(0,x)$ represents the expected number of visits
to $x$ for an infinite random walk starting at $0$.  Since the torus is finite,
the random walk is recurrent, and hence $C_{1/\cn}(0,x)$
is infinite for all $x$.  We therefore must keep $\mu<\cn^{-1}$ when dealing
with $C_\mu(0,x)$.  The value $\mu = \cn^{-1}$ plays the role
of the critical point for random walks.

Using translation invariance, we can write $C_\mu(x,y)=C_\mu(y-x)$.
By conditioning on the first step, we see that
the two-point function obeys the convolution equation
    \eq
    \lbeq{Cconv}
    C_\mu(x) = \delta_{0,x}+ \mu\cn(D * C_\mu) (x).
    \en
Taking the Fourier transform of \refeq{Cconv} gives $\hat{C}_\mu(k) = 1 +
\mu\cn \hat{D}(k)\hat{C}_\mu(k)$ and hence
    \eq
    \lbeq{Cdef}
    \hat C_{\mu}(k)=\frac{1}{1-\mu\cn\hat{D}(k)}.
    \en
Note that $\hat{C}_\mu(0) <\infty$ for $\mu<\cn^{-1}$ but
$\hat{C}_{1/\cn}(0)=\infty$.
Although $C_{1/\cn}(x)$ is infinite,
the formula \refeq{Cdef} does not diverge for $\mu=\cn^{-1}$ for all
$k$ for which $\hat{D}(k) \neq 1$.  Apart from any such singular points
(usually arising only for $k=0$),
the expression
$\hat{C}_{1/\cn}(k)=[1-\hat D(k)]^{-1}$ is finite.
The factor $[1-\hat D(k)]^{-3}$ that appears in \refeq{rwbd}
is thus the same as
$\hat{C}_{1/\cn}(k)^3$.  Comparing with \refeq{tri-k},
we see that \refeq{rwbd} is  closely
related to a triangle condition for random walks.

\subsection{Random walk estimates}
\label{sec-rwe}

In this section we prove Proposition \ref{lem-Cbd}, which for
convenience we restate as Proposition \ref{lem-Cbdbis}.


\begin{prop}\label{lem-Cbdbis}
There is an $a >0$ such that
random walks on each of the three tori listed in Section~\ref{sec-ptori}
obey Assumption~\ref{ass-rw} with $\beta = a\cn^{-1}$,
where:
\begin{enumerate}
\item for the nearest-neighbor torus, $a$ is a universal constant,
independent of $r\geq 2$ and $n \geq 7$;
\item for the Hamming torus, $a$ is a universal constant,
independent of $r\geq 2$ and $n \geq 1$;
\item for the spread-out torus, $n\geq 7$ is fixed, $r$ is sufficiently large
depending on $L$ and $n$, and $a$ depends on $n$ but not on $L$ or $r$.
\end{enumerate}
\end{prop}

The proof is given throughout the remainder of Section~\ref{sec-rwe}.
We first note that \refeq{supbds} is trivial since the maximal value
of $D(x)$ is $\cn^{-1}$ and this is less than $\beta = a\cn^{-1}$ provided
$a \geq 1$.  We verify the substantial assumption \refeq{rwbd} below.
As a first step, we discuss the infrared bound for the random walk models.

\subsubsection{The infrared bound}

For the Hamming torus, $D(x)$
is zero unless exactly one coordinate
of $x$ is different from zero, in which case it is equal to
$\cn^{-1}$.  If we denote the number of non-zero components of $k$ by $m(k)$,
we therefore have
    \eqalign
    \hat{D}(k) &=
    \frac 1\cn\sum_{j=1}^n\sum_{s=1}^{r-1}
    e^{ik_j s}
    =\frac 1\cn\sum_{j=1}^n\left(\sum_{s=0}^{r-1}
    e^{ik_j s} -1\right)
    \\
    &=\frac 1\cn\sum_{j=1}^n (r\delta_{k_j,0}-1)
    =
    1 - \frac {r}{r-1}\frac {m(k)}{n}.
    \lbeq{HTireq1}
    \enalign
This gives the infrared bound
\eqalign
    1-\hat{D}(k)
    =\frac {r}{r-1}\frac {m(k)}{n}\geq \frac {m(k)}{n}.
    \lbeq{HTireq2}
    \enalign

For $k \in \torus_{r,n}^* =  \frac{2\pi}{r}
\{-\lfloor \frac{r-1}{2}\rfloor,\ldots,\lceil
\frac{r-1}{2}\rceil \}^n$, we define
    \eq
    \lbeq{l2}
    |k|^2 = \sum_{j=1}^n k_j^2.
    \en
For
the nearest-neighbor torus,
by the symmetry of $D$ we have
\eq
    \hat{D}(k) = \sum_{x \in \torus_{r,n}} D(x) \cos(k\cdot x)
    =
    \lbeq{Dnn}
     \frac{1}{n} \sum_{j=1}^n \cos k_j.
    \en
Since $1-\cos t \geq 2\pi^{-2} t^2$ for $|t| \leq \pi$, this
implies the infrared bound
 \eq
    \lbeq{nnHTir}
    1-\hat{D}(k)
    \geq \frac{2}{\pi^2}\frac{|k|^2}{n}  .
\en

For the spread-out torus,
we first note that $\hat{D}(k)$
does not depend on $r$ if $r\geq 2L+1$.
Thus, we can
apply bounds on $\hat{D}(k)$ with $D(x)$ regarded as the step distribution
of a random walk on $\Z^n$.
The latter is analyzed in \cite[Appendix~A]{HS02},
where it is shown that there is an $\eta$ depending only on $n$
such that the infrared bound
    \eqalign
    1 - \hat{D}(k) &\geq \eta \big( 1 \wedge  L^2 |k|^2  \big)
    \lbeq{SOTir}
    \enalign
holds for all $k \in \torus_{r,n}^*$.

\subsubsection{The random walk triangle condition \refeq{rwbd}}

\medskip
\noindent
{\em
Proof of\refeq{rwbd} with $\beta = a\cn^{-1}$ for the Hamming torus.
}
Let $m(k)$ denote the number of non-zero components of $k$.
We fix an $\epsilon \in (0,1)$,
and divide the sum
\eq
    \frac{1}{V}
        \sum_{k \in \torus_{r,n}^*: k\neq 0}
        \frac{\hat{D}(k)^2}{[1-\hat{D}(k)]^3}
\en
according to whether $m(k)\leq \epsilon n$
or $m(k)>\epsilon n$.  It follows from
\refeq{HTireq2} that the contribution to the sum due to $m(k)>\epsilon n$
is bounded by
\eq
\lbeq{nHir1}
    \frac{1}{V}
        \sum_{k \in \torus_{r,n}^*: k\neq 0, \, m(k)> \epsilon n}
        \frac{\hat{D}(k)^2}{[1-\hat{D}(k)]^3}
        \leq
    \epsilon^{-3} \frac{1}{V}
    \sum_{k \in \torus_{r,n}^*}\hat{D}(k)^2 = \epsilon^{-3}\cn^{-1},
\en
since $V^{-1}$ times
the summation in the middle term is the probability that a random
walk returns to its starting point after two steps.

Note that the case $m(k)\leq \epsilon n$ does not occur for $n=1$
if we take $\epsilon <1$,
so we may assume that $n \geq 2$.
In this case, since $|\hat{D}(k)|\leq 1$, if follows
from \refeq{HTireq2} that
    \eq
    \lbeq{nH1a}
    \frac{1}{V}
        \sum_{k \in \torus_{r,n}^*: k\neq 0, \, m(k) \leq \epsilon n}
        \frac{\hat{D}(k)^2}{[1-\hat{D}(k)]^3}
    \leq
    \frac{1}{r^n}
    \sum_{m=1}^{\epsilon n} {n \choose m}(r-1)^m
        \frac{n^3}{m^{3}}.
    \en
In \refeq{nH1a}, the binomial coefficient counts the number of ways to choose
$m$ nonzero components from $n$, and the factor $(r-1)^m$ counts the number
of values that each nonzero component can assume.
Discarding two factors of $1/m$ and using
$n[m(n-m)]^{-1}\leq n(n-1)^{-1} \leq 2$,
the right hand side of \refeq{nH1a} is at most
    \eqalign
    \lbeq{nH2}
    n^3
    \sum_{m=1}^{\epsilon n}
    \frac 1m
    {n \choose m}\big(1-\frac{1}{r}\big)^m
    \big(\frac{1}{r}\big)^{n-m}
    & =
    \frac{n^3}{r}
    \sum_{m=1}^{\epsilon n}
    \frac{n}{m(n-m)}{{n-1} \choose m}\big(1-\frac{1}{r}\big)^m
    \big(\frac{1}{r}\big)^{n-1-m}
    \nnb
    & \leq
    2\frac{n^3}{r} \Pbold (X \leq \epsilon n),
    \enalign
where $X$ is a binomial random variable with parameters
$(n-1,1-r^{-1})$.  Let $p=1-r^{-1}$. Since $\epsilon n =
(p-a)(n-1)$ with $a \geq \frac 12 -2\epsilon$, it follows from
\eq
    \Pbold(X \leq (p-a)(n-1)) \leq e^{-(n-1)a^2/2}
\en
(a consequence of the Chernoff bound, see \cite[p.12]{Boll01})
that the right hand side of \refeq{nH2} decays exponentially in $n$,
uniformly in $r \geq 2$,
if we choose an $\epsilon  < \frac 14 $.
This gives the desired result.

\medskip \noindent
{\em Proof of\refeq{rwbd} with $\beta = a\cn^{-1}$ for
the nearest-neighbor torus.}
Since the nearest-neighbor torus is the same as the Hamming torus when
$r=2$ (in which case both are the $n$-cube), we may assume that $r \geq 3$.
Hence, $\cn = 2n$.

We first prove that
\eq
\lbeq{retpr}
    \frac{1}{V} \sum_{k \in \torus_{r,n}^*} \hat{D}(k)^{2i}
    \leq \frac{e 2^i i^{2i}}{\cn^i}.
\en
For this, we
observe that the left side is equal to the probability that a
random walk on $\torus_{r,n}$ that starts at the origin returns to
the origin after $2i$ steps.  This probability is equal to
$\Omega^{-2i}$ times the number of walks that make the transition
from 0 to 0 in $2i$ steps. Each such walk must take an even number
of steps in each coordinate direction, implying that it will live
in a subspace of dimension $\ell\leq \min\{i,n\}$.  If we fix the
subspace, and assume $r \geq 3$, then each step
in the subspace can be chosen from $2\ell$ different directions,
leading to a bound of $(2\ell)^{2i}$ for the number of walks in the
subspace.  Since the number of subspaces of fixed
dimension $\ell$ is given by ${n\choose\ell}\leq n^\ell/\ell!$, we
obtain the bound
\eq
    \sum_{\ell=1}^{i}\frac 1{\ell!}n^\ell (2\ell)^{2i}
    \leq
    \cn^i i^{2i} \sum_{\ell=1}^{i}\frac 1{\ell!}  2^{i}
    \leq
    \cn^i e 2^i i^{2i}
\en
for the number
of walks that make the transition from 0 to 0 in $2i$ steps.
Multiplying by $\cn^{-2i}$ to convert the number of walks into a
probability, this gives \refeq{retpr}.

By \refeq{retpr}
and H\"older's inequality, for any $i' >1$
\eq
\lbeq{rwCS}
    \frac{1}{V}
        \sum_{k \in \torus_{r,n}^*: k\neq 0}
        \frac{\hat{D}(k)^{2}}{[1-\hat{D}(k)]^3}
    \leq
    \left(\frac{e 2^{i'} (i')^{2i'}}{\cn^{i'}}\right)^{1/i'}
    \left(
    \frac 1V
    \sum_{k \in \torus_{r,n}^*: k\neq 0}
        \frac{1}{[1-\hat{D}(k)]^{3i}} \right)^{1/i},
\en
where $i=i'(i'-1)^{-1}$.  We choose $i'$ large enough that
$6i<7$.
By \refeq{rwCS}
and the infrared bound \refeq{nnHTir},
it suffices
to show that
\eq
\lbeq{wt1}
    \frac{1}{V}\sum_{k \in \torus_{r,n}^*: k\neq 0}
        \frac{n^{3i}}{|k|^{6i}}
\en
is bounded uniformly in $n\geq 7$ and $r\geq 3$.
%
Let $B(0,r)= (-\frac{\pi}{r},\frac{\pi}{r}]^n\subset \R^n$.
For each $k \in (-\pi,\pi]^n$,
there is a unique $k_r \in \torus_{r,n}^*$ such that $k\in k_r+B(0,r)$
and we define
    \eq
    \lbeq{Frk}
    F_r(k) =
    \begin{cases}
    \frac{n^{3i}}{|k_r|^{6i}} & (k_r \neq 0) \\
    0& (k_r=0).
    \end{cases}
    \en
Thus, $F_r(k)$ is constant on the cubes $k_r + B(0,r)$
for $k_r\in \torus_{r,n}^*$, the identity $V|B(0,r)|=(2\pi)^n$ holds, and
    \eq
    \lbeq{intF}
    \frac{1}{V}\sum_{k \in \torus_{r,n}^*: k\neq 0}
    \frac{n^{3i}}{|k|^{6i}} = \int_{[-\pi,\pi]^n}
    F_r(k)\frac{d^nk}{(2\pi)^n}.
    \en
We fix $\epsilon \in (0,1)$, and let $Q_1$ and $Q_2$ denote the subsets
of $(-\pi,\pi]^n$ for which $|k| \geq \frac{1}{1-\epsilon}\frac{\pi}{r}\sqrt{n}$
and $|k| \leq \frac{1}{1-\epsilon}\frac{\pi}{r}\sqrt{n}$, respectively.

For $k\in k_r + B(0,r)$,
we have $|k-k_r| \leq \frac \pi r \sqrt{n}$.
For $k \in Q_1$, it follows that $|k_r| \geq \epsilon |k|$, and hence
the contribution to the integral \refeq{intF} due to $k\in Q_1$ is at most
\eq
    \epsilon^{-6i}\int_{Q_1} \frac{n^{3i}}{|k|^{6i}} \frac{d^nk}{(2\pi)^n}
    \leq
    \epsilon^{-6i}\int_{(-\pi,\pi]^n}
    \frac{n^{3i}}{|k|^{6i}} \frac{d^nk}{(2\pi)^n}
    .
\en
The integral on the right hand side
is bounded uniformly in $n \geq 7$, by the following
argument.  For $A >0$ and $m >0$,
\eq
    \frac{1}{A^m} = \frac{1}{\Gamma(m)}\int_0^\infty t^{m-1} e^{-tA} dt.
\en
Applying this with $A=|k|^2/n$ and $m=3i$ gives
\eq
    \int_{[-\pi,\pi]^n} \frac{n^{3i}}{|k|^{6i}}\frac{d^nk}{(2\pi)^n}
    =
    \frac{1}{\Gamma(3i)}\int_0^\infty dt \, t^{3i-1}
    \Big( \int_{-\pi}^\pi \frac{d\theta}{2\pi}
    \left(e^{-t\theta^2 }\right)^{1/n}
    \Big)^n.
\en
The right side is non-increasing in $n$, since $\|f\|_p \leq \|f\|_q$
for $0 < p \leq q \leq \infty$ on a probability space.

For $k \in Q_2$, we use the fact that
$|k_r| \geq \frac{2\pi}{r}$ for all non-zero $k_r$ to obtain
\eq
    \int_{Q_2}F_r(k) \frac{d^nk}{(2\pi)^n}
    \leq
   \frac {n^{3i}r^{6i}}{(2\pi)^{6i}}\int_{Q_2}\frac{d^nk}{(2\pi)^n}
    =
    \frac {n^{3i}r^{6i}}{(2\pi)^{6i}}
    \frac{1}
    {(2\pi)^n} v_n
    \left(\frac{1}{1-\epsilon}\frac{\pi}{r}\sqrt{n}\right)^n,
\en
where $v_n$ denotes the volume of the unit ball in $n$ dimensions.
Since
\eq
    v_n = \frac{\pi^{n/2}}{\Gamma(\frac{n}{2}+1)}
    \leq \left( \frac{2\pi e}{n}\right)^{n/2}
\en
(using $\Gamma(a+1) \geq a^ae^{-a}$), this gives
\eq
\lbeq{Q2int}
    \int_{Q_2}F_r(k) \frac{d^nk}{(2\pi)^n}
    \leq
    \frac {n^{3i}r^{6i}}{(2\pi)^{6i}}
    \left( \frac{e\pi}{(1-\epsilon)^2 2r^2}\right)^{n/2}.
\en
We fix $\epsilon$ so that
\eq
    \frac{e\pi}{(1-\epsilon)^2 2 }<9.
\en
The right hand side of \refeq{Q2int} is then bounded
uniformly in $n \geq 1$ and $r \geq 3$.

\medskip \noindent
{\em Proof of \refeq{rwbd} with $\beta = a\cn^{-1}$ for
the spread-out  torus.}
Now $n\geq 7$ is fixed, $L$ is fixed and large, and $r \to \infty$.
We first use the dominated convergence theorem to show that
\eq
\lbeq{sowt1}
    \lim_{r \to \infty}
    \frac{1}{V}
        \sum_{k \in \torus_{r,n}^*: k\neq 0}
        \frac{\hat{D}(k)^{2}}{[1-\hat{D}(k)]^3}
    =\int_{[-\pi,\pi]^n} \frac{\hat{D}(k)^{2}}{[1-\hat{D}(k)]^3}
    \frac{d^nk}{(2\pi)^n}.
\en
It follows that the expression under the limit on the left hand side is bounded
above by twice the integral on the right hand side,
for $r$ sufficiently large depending on $L,n$.

Recalling the definition of $k_r$ above \refeq{Frk}, we
define
    \eq
    G_r(k) =
    \begin{cases}
    \hat{D}(k_r)^2 [1-\hat{D}(k_r)]^{-3} & (k_r \neq 0) \\
    0 & (k_r=0),
    \end{cases}
    \en
so that
    \eq
    \frac{1}{V}\sum_{k \in \torus_{r,n}^*: k\neq 0}
    \frac{\hat{D}(k)^2}{[1-\hat{D}(k)]^{3}} = \int_{[-\pi,\pi]^n}
    G_r(k) \frac{d^nk}{(2\pi)^n}.
    \en
The function $G_r(k)$ converges pointwise to
$\hat{D}(k)^2[1-\hat{D}(k)]^{-3}$ for $k\neq 0$ and to 0 for
$k=0$.
Also, by the infrared bound \refeq{SOTir},
    \eq
    G_r(k)
    \leq \frac{1}{\eta^{3}L^3 |k_r|^{6}}
    \en
for every nonzero $k_r\in \torus^*_{r,n}$.
For each nonzero $k_r\in \torus^*_{r,n}$ and $k\in k_r + B(0,r)$,
we have $\|k_r\|_\infty \geq \frac{2\pi}{r}$ and
$\|k-k_r\|_\infty \leq \frac{\pi}{r}$, so that
$\|k_r\|_\infty\geq \frac{2}{3}\|k\|_\infty$.   This implies that
$|k_r|^2 \geq \|k_r\|_\infty^2 \geq \frac{4}{9} \|k\|_\infty^2 \geq
\frac{4}{9n}|k|^2$.
Therefore, for every $k\in (-\pi,\pi]^n$,
    \eq
    G_r(k)
    \leq \big(\frac{9}{4}\big)^{3} \frac{n^{3}}{\eta^{3}L^3 |k|^{6}},
    \en
which is integrable when $n\geq 7 $.
Therefore, by dominated convergence, \refeq{sowt1}
holds, and it suffices to bound
the integral on the right hand side of \refeq{sowt1}.

We bound the integral on the right side of
\refeq{sowt1} by considering separately the regions
where $|k|^2 > L^{-2}$ and $|k|^2 \leq L^{-2}$.
For the contribution to the integral on the right side of \refeq{sowt1}
due to $|k|^2 > L^{-2}$, we use \refeq{SOTir} and argue as in
\refeq{nHir1} to obtain
    \eq
    \int_{k\in [-\pi,\pi]^n: |k|^2 > L^{-2}}
    \frac{\hat{D}(k)^{2}}{[1-\hat{D}(k)]^3}
    \frac{d^nk}{(2\pi)^n} \leq \eta^{-3}
    \int_{[-\pi,\pi]^n} \hat{D}(k)^{2} \frac{d^nk}{(2\pi)^n}
    = \eta^{-3}\cn^{-1}.
    \en
For the contribution to the integral in \refeq{sowt1}
due to $|k|^2 \leq L^{-2}$,
we use \refeq{SOTir} and $\hat{D}(k)^{2} \leq 1$ to obtain
    \eq
    \lbeq{227}
    \int_{|k|^2\leq  L^{-2}}
    \frac{\hat{D}(k)^{2}}{[1-\hat{D}(k)]^3}
    \frac{d^nk}{(2\pi)^n} \leq
    \frac{1}{\eta^3 L^{6}} \int_{|k|^2\leq L^{-2}}
    \frac1{|k|^6}\frac{d^nk}{(2\pi)^n}
    = C_{n,\eta}  L^{-n}.
    \en
Summing the two contributions yields
\refeq{rwbd} with $\beta = a\cn^{-1}$.

\subsubsection{A consequence of Assumption~\ref{ass-rw}}

Finally, we note for future reference that \refeq{rwbd}
implies that
\eq
\lbeq{rwbdi0}
    \frac{1}{V}
    \sum_{k \in \torus_{r,n}^*: k\neq 0}
    \frac{1}{[1-\hat{D}(k)]^3} \leq
    1+6\beta .
\en
To see this, we use the identity
\eq
    \frac{1}{[1-\hat{D}]^3} = 1 + 3\hat{D} +
    \frac{3\hat{D}^2}{1-\hat{D}} + \frac{2\hat{D}^2}{[1-\hat{D}]^2}
    + \frac{\hat{D}^2}{[1-\hat{D}]^3}.
\en
The sum of the last three terms on the right side is
at most $6\hat{D}^2[1-\hat{D}]^{-3}$, and their
normalized sum over $k$ is
thus at most $6\beta$, by \refeq{rwbd}.  Since the
normalized sum over {\em all}
$k\in \torus_{r,n}^*$ of $3\hat{D}(k)$ is $3D(0)=0$, its sum over nonzero
$k$ is $-3V^{-1}<0$.  This proves \refeq{rwbdi0}.

\section{The lace expansion}
\label{sec-lexp}

We begin in Section~\ref{sec-le} with a brief overview of the lace
expansion, and then give a self-contained and detailed derivation of
the expansion in Section~\ref{sec-exp}.

The term ``lace'' was used by Brydges and Spencer
\cite{BS85} for a certain
graphical construction that arose in the expansion they invented to
study the self-avoiding walk.  Although the lace expansion for percolation
evolved from the expansion for the self-avoiding walk, this
graphical construction
does not occur for percolation, and so the term ``lace'' expansion is
a misnomer in the percolation context.  However, the name has stuck
for historical reasons.

\subsection{Overview of the lace expansion}
\label{sec-le}

In this section, we give a brief introduction to the lace expansion,
with an indication of how it is used to prove the triangle condition
of Theorem~\ref{thm-tc}.  Since the discussion will
involve the Fourier transform, we restrict attention here
to percolation on the narrow torus
$\torus_{r,n}$, with $r \geq 2$ and $n$ large.  Each vertex has degree
$\cn = 2n$ for $r \geq 3$ and $\cn = n$ for $r=2$.  However, in
Section~\ref{sec-exp} the expansion will be derived on an arbitrary
finite graph $\gr$.

Given a percolation cluster containing $0$ and $x$, we call any bond
whose removal would disconnect $0$ from $x$ a {\em pivotal} bond.
The connected components that remain after removing all pivotal
bonds are called {\em sausages}.
Since they are separated by at least one pivotal bond by definition,
no two sausages can have a common vertex.  Thus, the sausages are constrained
to be mutually avoiding.  However, this is a weak constraint, since
sausage intersections require a cycle, and cycles are unlikely.
In fact, for $p$ asymptotically
proportional to $\cn^{-1}$, and hence for $p=p_c$,
the probability that the origin is in a cycle of length
4 is of order $\cn^2 \cn^{-4} = \cn^{-2}$, and
larger cycles are more unlikely.  The fact that cycles are unlikely
also means that sausages tend to be trees.
This makes it
reasonable to attempt to apply an inclusion-exclusion analysis, where
the connection from 0 to $x$ is treated as a random walk
path, with correction terms taking into account cycles in sausages
and the
avoidance constraint between sausages.
With this in mind, it makes sense to attempt to relate $\tau_p(0,x)$
to the two-point function for random walks.

\begin{figure}
\begin{center}
\setlength{\unitlength}{0.0100in}%
\begin{picture}(185,127)(60,680)
\thicklines
\put(100,760){\line( 1, 0){ 20}}
\put(120,740){\line( 0,-1){ 20}}
\put(120,720){\line( 1, 0){ 20}}
\put(140,720){\line( 1, 0){ 20}}
\put(160,720){\line( 1, 0){ 20}}
\put(180,720){\line( 0, 1){ 20}}
\put(180,740){\line( 1, 0){ 20}}
\thinlines
\put(100,780){\line( 0,-1){ 60}}
\put( 80,760){\line( 1, 0){ 20}}
\put( 80,760){\line( 0,-1){ 20}}
\put( 80,740){\line(-1, 0){ 20}}
\put(120,760){\line( 1, 0){ 20}}
\put(140,760){\line( 0, 1){ 40}}
\put(120,780){\line( 1, 0){ 60}}
\put(120,720){\line( 0,-1){ 20}}
\put(120,740){\line( 1, 0){ 40}}
\put(160,740){\line( 0, 1){ 20}}
\put(160,720){\line( 0,-1){ 40}}
\put(160,700){\line( 1, 0){ 20}}
\put(180,760){\line( 0,-1){ 20}}
\put(180,760){\line( 1, 0){ 40}}
\put(200,780){\line( 0,-1){ 20}}
\put(200,740){\line( 0,-1){ 40}}
\put(200,740){\line( 1, 0){ 40}}
\put(220,740){\line( 0,-1){ 20}}
\put(220,720){\line( 1, 0){ 20}}
\put(240,720){\line( 0,-1){ 20}}
\put(100,760){\makebox(0,0)[lb]{\raisebox{0pt}[0pt][0pt]{\circle*{3}}}}
\put(120,760){\makebox(0,0)[lb]{\raisebox{0pt}[0pt][0pt]{\circle*{3}}}}
\put(120,740){\makebox(0,0)[lb]{\raisebox{0pt}[0pt][0pt]{\circle*{3}}}}
\put(120,720){\makebox(0,0)[lb]{\raisebox{0pt}[0pt][0pt]{\circle*{3}}}}
\put(140,720){\makebox(0,0)[lb]{\raisebox{0pt}[0pt][0pt]{\circle*{3}}}}
\put(160,720){\makebox(0,0)[lb]{\raisebox{0pt}[0pt][0pt]{\circle*{3}}}}
\put(180,720){\makebox(0,0)[lb]{\raisebox{0pt}[0pt][0pt]{\circle*{3}}}}
\put(180,740){\makebox(0,0)[lb]{\raisebox{0pt}[0pt][0pt]{\circle*{3}}}}
\put(200,740){\makebox(0,0)[lb]{\raisebox{0pt}[0pt][0pt]{\circle*{3}}}}
\put( 90,765){\makebox(0,0)[lb]{\raisebox{0pt}[0pt][0pt]{${\sss 0}$}}}
\put(205,745){\makebox(0,0)[lb]{\raisebox{0pt}[0pt][0pt]{${\sss x}$}}}
\put(120,760){\line( 0,-1){ 20}}
\put(140,760){\line( 0,-1){ 20}}
\put(160,760){\line( 0, 1){ 20}}
\put(200,720){\line( 1, 0){ 20}}
\end{picture}
%
\setlength{\unitlength}{0.001in}
\begin{picture}(2602,141)(0,-10)
\put(430,670){\makebox(0,0)[lb]{\raisebox{0pt}[0pt][0pt]{${\sss 0}$}}}
\put(3170,670){\makebox(0,0)[lb]{\raisebox{0pt}[0pt][0pt]{${\sss x}$}}}
\put(2982,700){\ellipse{224}{112}}
\put(2645,700){\ellipse{224}{112}}
\put(2307,700){\ellipse{224}{112}}
\put(1970,700){\ellipse{224}{112}}
\put(1632,700){\ellipse{224}{112}}
\put(1295,700){\ellipse{224}{112}}
\put(957,700){\ellipse{224}{112}}
\put(620,700){\ellipse{224}{112}}
\thicklines
\path(2757,700)(2870,700)
\path(2420,700)(2532,700)
\path(2082,700)(2195,700)
\path(1745,700)(1857,700)
\path(1407,700)(1520,700)
\path(1070,700)(1182,700)
\path(732,700)(845,700)
\end{picture}

\end{center}
\caption{A percolation cluster with a string of 8 sausages joining $0$ to $x$,
and a schematic
representation of the string of sausages.
The 7 pivotal bonds are shown in bold.}
\label{fig-sos}
\end{figure}
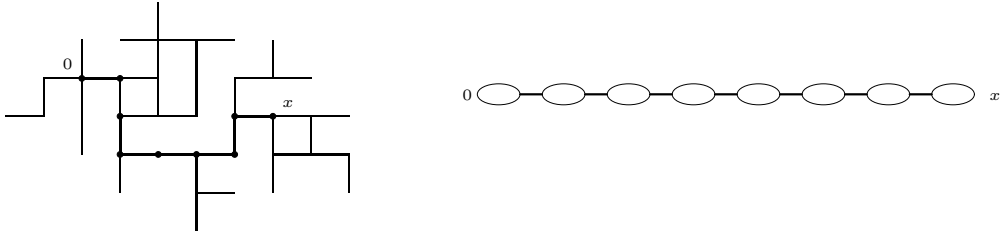

The lace expansion of Hara and Slade \cite{HS90a} makes this procedure
precise.  It produces a convolution equation of the form
    \eq
    \lbeq{taux'}
    \tau_p(0,x) = \delta_{0,x}+ p\cn (D * \tau_p) (0,x) +
        p\cn (\Pi_p * D * \tau_p) (0,x)
        + \Pi_p(0,x)
    \en
for the two-point function, valid for $p \leq p_c(\torus_{r,n})$.
The expansion gives explicit
but complicated formulas for the function
$\Pi_p: \torus_{r,n}\times \torus_{r,n} \to \R$.  It
will turn out that if Assumption~\ref{ass-rw} holds with $\beta = O(\cn^{-1})$,
then $\hat \Pi_p(k)=\bigo(\cn^{-1})$
uniformly
in $p\leq p_c(\gr)$.  Putting $\Pi_p \equiv 0$ in \refeq{taux'} gives
\refeq{Cconv}, and in this sense the percolation two-point function
can be regarded
as a small perturbation of the random walk two-point function.

Using \refeq{ftconv}, \refeq{taux'} can be solved to give
        \eq
        \lbeq{tauk'}
        \hat \tau_p(k) = \frac{1+\hat \Pi_p(k)}
        {1-p\cn \hat{D}(k)[1+\hat \Pi_p(k)]}.
        \en
We will show that under Assumption~\ref{ass-rw},
$\hat \Pi_p(k)$ can be well approximated
by $\hat \Pi_p(0)$.  Since $\hat \Pi_p(k)$ is also small compared to 1,
\refeq{tauk'}
suggests that the approximation
\eq
\lbeq{tauapprox}
    \hat \tau_p(k)
    \approx
    \frac{1}{1-p\cn [1+\hat \Pi_p(0)]\hat{D}(k)}
\en
is reasonable (where $\approx$ denotes an uncontrolled
approximation).  Comparing with \refeq{Cdef}, this suggests that
\eq
\lbeq{tauappC}
    \hat \tau_p(k) \approx \hat{C}_{\mu_p}(k)
    \quad \quad
    \mbox{with} \quad \quad \mu_p\cn = p\cn [1+\hat \Pi_p(0)].
\en
We will make this approximation precise in \refeq{tauMkC6}.
Since $\hat{D}(0)=1$, if we set $k=0$
in \refeq{tauk'} and solve for $p\cn$ then we obtain
        \eq
        \lbeq{pcalmost}
        p\cn  = \frac{1}{1+\hat \Pi_{p}(0)}- \hat{\tau}_p(0)^{-1}.
        \en
For $p=p_c=p_c(\torus_{r,n})$, \refeq{pcalmost} states that
\eq
\lbeq{pcform}
    p_c\cn
    = \frac{1}{1+\hat \Pi_{p_c}(0)} - \lambda^{-1} V^{-1/3},
\en
and hence $\mu_{p_c}\cn \approx 1-\lambda^{-1} V^{-1/3}$.  This should
be compared with the
critical value $\mu\cn =1$ for the random walk.

For the triangle condition, we analyze the Fourier representation of
$\nabla_p(0,x)$ given in \refeq{tri-k}.
Extraction of the $k=0$ term in \refeq{tri-k} gives
\eq
    \nabla_p(0,x) =
    \frac{\chi(p)^3}{V} +
     \frac{1}{V}
     \sum_{k \in \torus_{r,n}^* : \, k \neq 0} \hat{\tau}_{p}(k)^3
    e^{-ik \cdot x}.
\en
The second term can be estimated using \refeq{tauappC} and
Assumption~\ref{ass-rw}, leading to a proof of Theorem~\ref{thm-tc}.
Details are given in Section~\ref{sec-Pibds}.

Equation~\refeq{pcform} is an implicit equation for the critical threshold.
Using $\hat \Pi_p(k)=\bigo(\cn^{-1})$, \refeq{pcform} gives
$p_c = \cn^{-1} + \bigo(\cn^{-2})$ if
$\lambda^{-1}V^{-1/3}\leq \bigo(\cn^{-2})$.
This is the first term in an
asymptotic expansion.  Further terms
will follow from an asymptotic expansion of $\hat{\Pi}_{p_c}(0)$
in powers of $\cn^{-1}$.  Calculation of this sort were
carried out in \cite{HS95} for percolation on
$\Z^n$, and have been extended in \cite{HS04a,HS04b},
as was discussed below Theorem~\ref{thm-qnsuper}.

\subsection{Derivation of the lace expansion}
\label{sec-exp}

In this section, we derive a version of
the lace expansion \refeq{taux'} that
contains a remainder term.
Throughout this section, we let $\gr$
denote an arbitrary finite graph, which need not be transitive
nor regular. We use
the method of \cite{HS90a}, which applies directly in this general
setting,
and we follow the presentation of \cite{MS93}.
We assume for simplicity that $\gr$ is finite, but with minor
modifications the analysis also applies when $\gr$ is infinite
provided there is almost surely no infinite cluster.

Fix $p \in [0,1]$.  We define
\eqalign
    J(x,y) &= pI[\{x,y\} \in \mathbb B]
    \nnb
    & = p\cn D(x,y) \text{   if $\gr$ is regular},
\enalign
with $D$ given by \refeq{Ddef}.
We write $\tau(x,y)=\tau_p(x,y)$ for brevity,
and generally drop subscripts indicating dependence on $p$.
For each $M = 0,1,2,\ldots$, the expansion takes the form
    \eq
    \lbeq{tauxM}
        \tau(x,y) = \delta_{x,y}+  (J * \tau) (x,y) +
          (\Pi_{\SSS M} * J * \tau) (x,y)
        + \Pi_{\SSS M}(x,y) +R_{\SSS M}(x,y),
    \en
where the $*$ product denotes matrix multiplication (this reduces
to convolution when $\gr =\torus_{r,n}$).
The function
$\Pi_{\SSS M} : \ver\times\ver \to \R$ is the key quantity in the expansion, and
$R_{\SSS M}(x,y)$ is a remainder term.
The dependence of $\Pi_{\SSS M}$ on $M$ is given by
    \eq
    \lbeq{2pt.37}
    \Pi_{\SSS M}(x,y) = \sum_{N=0}^{M} (-1)^N
    \Pi^{\smallsup{N}}(x,y),
    \en
with $\Pi^\smallsup{N}(x,y)$ independent of $M$.  The alternating sign
in \refeq{2pt.37} arises via repeated inclusion-exclusion.
In Section~\ref{sec-pfstPibds}, we will prove that for
$\gr = \torus_{r,n}$, $p \leq p_c(\torus_{r,n})$, and assuming
Assumption~\ref{ass-rw} with $\lambda^3 \vee \beta$ sufficiently
small,
    \eq
    \lim_{M \to \infty}
    \sum_y|R_{\SSS M}(x,y)| = 0 .
    \en
This leads to \refeq{taux'} with $\Pi = \Pi_\infty$.  Convergence
properties of \refeq{2pt.37} when $M=\infty$ will also be established
in Section~\ref{sec-pfstPibds}.  The remainder of this section gives
the proof of \refeq{tauxM}.

We need the following definitions.

    \begin{defn}
    \label{def-inon}
    (a)
    Given a bond configuration, and $A \subset \mathbb V$, we
        say $x$ and $y$ are \emph{connected in $A$}, and write
        $\{x \conn y \text{ in } A\}$,
        if $x = y \in A$ or if there is an
        occupied path from $x$ to $y$ having all its endpoints in
        $A$.
        We define a restricted two-point function by
            \eq
            \tau^{A} (x, y) = \prob{\text{$x\conn y$
            in $\mathbb V\backslash A$}}
            .
            \en
        (b)
        Given a bond configuration, and $A\subset \mathbb V$, we
        say $x$ and $y$ are \emph{connected through $A$}, if
    $x \conn y$ and every
        occupied path connecting $x$ to $y$ has at least one bond
        with an endpoint in $A$, or if $x=y\in A$.
        This event is written as
        $x \ct{A} y$.
    \newline
        (c)
        Given a bond configuration, and a bond $b$, we define
        $\tilde{C}^{b}(x)$ to be the set of sites connected to $x$
        in the new configuration obtained by setting $b$ to be vacant.
        \newline
    (d)
Given a bond configuration, we say that $x$ is {\em doubly connected
    to}\/ $y$, and we write $x \dbc y$, if $x=y$ or
    if there are at least two bond-disjoint
    paths from $x$ to $y$
    consisting of occupied bonds.
\newline
(e)
Given a bond configuration, a bond $\{u,v\}$ (occupied or not) is called
\emph{pivotal}
for the connection from $x$ to $y$ if (i) either $x \conn u$ and
$y \conn v$, or $x \conn v$ and $y \conn u$, and (ii)
$y \, \nin \, \tilde{C}^{\{u,v\}}(x)$.
Bonds are not usually regarded
as directed.  However, it will be convenient at times
to regard a bond $\{u,v\}$ as directed from $u$ to $v$, and we will
emphasize this point of view with the notation $(u,v)$.
A directed bond $(u,v)$ is pivotal for the connection from $x$ to
$y$ if $x \conn u$, $v\conn y$ and
$y \, \nin \, \tilde{C}^{\{u,v\}}(x)$.  We denote by $P_{(x,y)}$
the set of directed pivotal bonds for the connection from $x$ to $y$.
    \end{defn}

To begin the expansion, we define
        \eq\lbeq{pi0def}
        \Pi^{\smallsup{0}}(x,y)=\prob{x \dbc  y}
        - \delta_{x, y}
        \en
and distinguish configurations with $x \conn y$ according
to whether or not there is a double connection, to obtain
    \eq
    \lbeq{2ptNm.1}
        \tau(x,y)
        = \delta_{x, y} +
        \Pi^{\smallsup{0}}(x,y)
        +\prob{x \conn y \AND x \ndbc
        y}.
    \en
If $x$ is connected to $y$, but not doubly, then $P_{(x,y)}$
is nonempty.  There is therefore a unique element $(u,v) \in P_{(x,y)}$
(the \emph{first} pivotal bond)
such that $x \dbc u$, and we can  write
    \eq
    \lbeq{2pt.11}
        \prob{x \conn y \AND x \ndbc  y}
        = \sum_{(u, v)} \prob{x \dbc u \text{ and
        $(u, v)$ is occupied and pivotal for } x \conn  y }.
    \en
Now comes the essential part of the expansion.
Ideally, we would like to factor
the probability on the right side of \refeq{2pt.11} as
    \eq
    \lbeq{2pt.14}
        \prob{x \dbc u} \,
        \prob{\text{$(u,v)$ is occupied}}
        \, \prob{v \conn y}
        = \big( \delta_{x, u} + \Pi^{\smallsup{0}}(x,u) \big)
        J(u,v) \tau(v,y).
    \en
The expression \refeq{2pt.14} is the same as
\refeq{tauxM}
with $\Pi_{\SSS M}= \Pi^\smallsup{0}$ and $R_{\SSS M}
=0$.
However, \refeq{2pt.11} does not factor in this way
because the cluster $\tilde{C}^{(u,v)}(u)$ is constrained
not to intersect the cluster $\tilde{C}^{(u,v)}(v)$,
since $(u,v)$ is pivotal.
What we can do is approximate the probability on the right side
of \refeq{2pt.11} by \refeq{2pt.14}, and then attempt to deal with the
error term.

For this, we will use the next lemma, which gives an identity
for the probability on the right hand side of \refeq{2pt.11}.
In fact, we will also need a more general identity, involving
the following generalizations of the event appearing on
the right hand side of \refeq{2pt.11}.
Let $x,u,v,y \in \ver$, and $A \subset \ver$
be nonempty.  Then we define the events
    \eqalign
    \lbeq{317}
    E'(v, y; A) & = \{ v \ct{A} y\} \cap
    \{\not\exists (u', v') \in P_{(v,y)} \text{ such that }
    v \ct{A} u' \}
    \enalign
and
\eqalign
    \lbeq{317a}
    E(x, u, v, y; A) & = E'(x, u; A) \cap
    \{\text{$(u, v)$ is occupied and pivotal for $x \conn
    y$}\}.
    \enalign
Note that $\{x \dbc y\} = E'(x,y; \ver)$,
while $E(x, u, v, y; \ver)$ is the event appearing on the
right hand side of \refeq{2pt.11}.
A version of Lemma~\ref{lem-cut1}, with $E'(x,u;A)$ replaced
by $\{0 \conn u\}$ on both sides of \refeq{percpivineq},
appeared in \cite{AN84}.

\begin{figure}
\begin{center}
\setlength{\unitlength}{0.0080in}%
\begin{picture}(500,100)(100,-50)
\put(120,20){${\sss x}$}
\put(350,25){${\sss u}$}
\put(360,25){${\sss v}$}
\put(520,20){${\sss y}$}

\shade\path(270,-45)(275,-50)(285,-47)(295,-41)
(305,-32)(315,-20)(325,-5)(335,13)(345,34)(346,36)
(335,31)(325,25)(315,20)
(305,12)(295,5)(285,-8)(275,-24)(270,-37)(269,-41)
(270,-45)

\shade\path(155,-65)(160,-70)(170,-67)(180,-61)
(190,-52)(200,-40)(210,-25)(220,-7)(230,14)(230,16)
(220,11)(210,5)(200,0)
(190,-8)(180,-15)(170,-28)(160,-44)(155,-57)(154,-61)
(155,-65)

\put(245,-20){${\scriptstyle A}$}

\qbezier(130,20)(160,-10)(190,20)
\qbezier(130,20)(160,50)(190,20)
\put(190,20){\line(1,0){20}}
\qbezier(210,20)(240,-10)(270,20)
\qbezier(210,20)(240,50)(270,20)
\put(270,20){\line(1,0){20}}
\qbezier(290,20)(320,-10)(350,20)
\qbezier(290,20)(320,50)(350,20)
\put(350,20){\line(1,0){20}}
\qbezier(370,20)(400,-10)(430,20)
\qbezier(370,20)(400,50)(430,20)
\put(430,20){\line(1,0){20}}
\qbezier(450,20)(480,-10)(510,20)
\qbezier(450,20)(480,50)(510,20)
\end{picture}
\end{center}
\caption{\lbfg{perclemevent}
The event
$E(x,u,v,y;A)$
of Lemma~\protect\ref{lem-cut1}.  The
shaded regions
represent the vertices in $A$.  There is no restriction on intersections
between $A$ and $\tilde{C}^{\{u,v\}}(x)$.}
\end{figure}
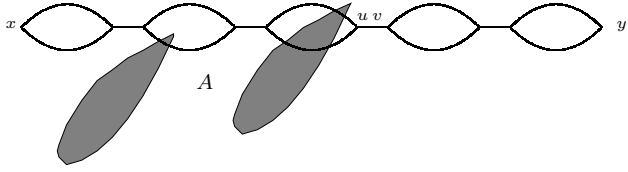

\begin{lemma}
\label{lem-cut1}
Let $\gr$ be a finite graph, $p \in [0,1]$, $u\in \ver$,
and let  $A\subset \ver$ be nonempty.  Then
\eqarray
\lbeq{percpivineq}
    \Ebold \left(I[E(x, u, v, y; A)]\right)
    = p\Ebold\left( I[E'(x,u;A)]\, \tau^{\tilde{C}^{\{u,v\}}(x)}
    (v,y)\right) .
\enarray
\end{lemma}

\proof
The event appearing in the left side of \refeq{percpivineq} is depicted
in Figure~\reffg{perclemevent}.
We first observe that the event
$E'(x,u;A) \cap \{(u,v)\in P_{(x,y)}\}$
is independent of the occupation status of the bond $(u,v)$.
This is true by definition for $\{(u,v) \in P_{(x,y)}\}$, and when
$(u,v)$ is pivotal, the occurrence or not of $E'(x,u;A)$
cannot be affected by $\{u,v\}$ since in this case
$E'(x,u;A)$ is determined by the occupied
paths from $x$ to $u$ and no such path
uses the bond $\{u,v\}$.
Therefore, the left side of the identity in the statement of the
lemma is equal to
\eq
\lbeq{perclem0}
    p \Ebold \left(I[E'(x,u;A)\cap \{(u,v) \in P_{(x,y)}\}]\right) .
\en
By conditioning on $\tilde{C}^{\{u,v\}}(x)$,
\refeq{perclem0} is equal to
\eq
\lbeq{perclem1}
    p \sum_{S:S \ni x}\Ebold\left(
    I[E'(x,u;A)\cap \{ (u,v) \in P_{(x,y)}\}
    \cap\{ \tilde{C}^{\{u,v\}}(x)=S \}]\right),
\en
where the sum is over all finite connected sets of vertices $S$ containing $x$.

In \refeq{perclem1}, we make the replacement
\eq
\lbeq{perclem9}
    \{ (u,v) \in P_{(x,y)}\}
    \cap\{ \tilde{C}^{\{u,v\}}(x)=S \}
    =
    \{ v \conn y \text{ in } \ver\backslash S\}
    \cap\{ \tilde{C}^{\{u,v\}}(x)=S \}.
\en
The event $\{ v \conn y $ in $\ver\backslash S\}$
depends only on the occupation status
of bonds which do not have an endpoint in $S$.
On the other hand, given that $\{ v \conn y$ in $\ver\backslash S\}
    \cap\{ \tilde{C}^{\{u,v\}}(x)=S \}$ occurs,
the event $E'(x,u;A)$
is determined by the occupation status of bonds which have an endpoint in
$S=\tilde{C}^{\{u,v\}}(x)$.
Similarly, the event $\{\tilde{C}^{\{u,v\}}(x)
=S\}$ depends on bonds
which have one or both endpoints in $S$.  Hence, given $S$, the event
$E'(x,u;A) \cap \{\tilde{C}^{\{u,v\}}(x)=S\}$ is independent of the event that
$\{ v \conn y$ in $\ver\backslash S\}$,
and therefore \refeq{perclem1} is equal to
\eq
\lbeq{perclem2}
    p \sum_{S:S \ni x} \Ebold \left(
    I[E'(x,u;a) \cap \{ \tilde{C}^{\{u,v\}}(x) =S\}]\right) \tau_p^{S}(v,y).
\en
Bringing the restricted two-point function inside the expectation,
replacing the superscript $S$ by $\tilde{C}^{\{u,v\}}(x)$, and performing the
sum over $S$, gives the desired result.
\qed

It follows from \refeq{2pt.11} and Lemma~\ref{lem-cut1} that
\eqalign
\lbeq{2pt.21}
    \prob{x \conn y \AND x \ndbc  y}
    & =
    \sum_{(u,v)} J(u,v)
    \Ebold \left(  I[x \dbc u ] \,
    \tau^{\tilde{C}^{(u, v)}(x)}(v, y)
    \right)
    .
\enalign
On the right side,
$\tau^{\tilde{C}^{(u, v)}(x)}(v, y)$ is the
restricted two-point function \emph{given} the cluster $\tilde{C}^{(u,
v)}(x)$ of the outer expectation, so that in the
expectation defining $\tau^{\tilde{C}^{(u,
v)}(x)}(v, y)$, $\tilde{C}^{(u,
v)}(x)$ should be regarded as a \emph{fixed} set. We
stress this delicate point here, as it is crucial also in the
rest of the expansion.  The
inner expectation on the right side effectively introduces a second
percolation model on a second graph, which depends on the original
percolation model via the set $\tilde{C}^{(u, v)}(x)$.

We write
\eqalign
\lbeq{2pt.23}
    & \tau^{\tilde{C}^{(u, v)}(x)}(v, y)
    = \tau(v,y) - \left( \tau(v,y)
    - \tau^{\tilde{C}^{(u, v)}(x)}(v, y)
    \right)
    = \tau(v,y) -
    \probb{v \ctx{\tilde{C}^{(u, v)}(x)} y},
\enalign
insert this into \refeq{2pt.21}, and use \refeq{2ptNm.1}
and \refeq{pi0def} to obtain
\eqalign
\lbeq{2pt.26}
    \tau(x,y)
    & = \delta_{x,y} + \Pi^{\smallsup{0}}(x,y)
    + \sum_{(u,v)}
    \big( \delta_{x, u} + \Pi^{\smallsup{0}}(x,u) \big)
    J(u,v)
    \tau(v,y)
    \nnb & \hspace{5mm}
     - \sum_{(u,v)} J(u,v)
    \Ebold \left(  I [ x \dbc u ] \,
    \prob{v \ctx{\tilde{C}^{(u, v)}(x)} y }
    \right).
\enalign
With $R_{\SSS 0}(x,y)$ equal to the last term on the right side
of \refeq{2pt.26} (including the minus sign), this proves
\refeq{tauxM} for $M=0$.

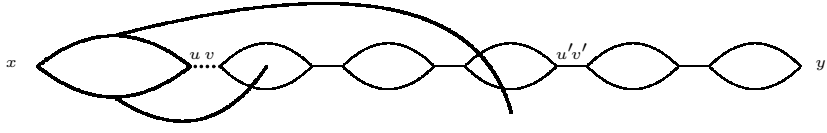
\begin{figure}
\begin{center}
\setlength{\unitlength}{0.0080in}%
\begin{picture}(480,50)
\put(-10,20){${\sss x}$}
\put(110,25){${\sss u}$}
\put(120,25){${\sss v}$}
\put(350,25){${\sss u'}$}
\put(360,25){${\sss v'}$}
\put(520,20){${\sss y}$}
\thicklines
\qbezier(10,20)(60,-20)(110,20)
\qbezier(10,20)(60,60)(110,20)
\qbezier(60,40)(290,100)(320,-10)
\qbezier(60,0)(120,-40)(160,20)
\thinlines
\multiput(110,20)(4,0){6}{\circle*{2}}
\qbezier(130,20)(160,-10)(190,20)
\qbezier(130,20)(160,50)(190,20)
\put(190,20){\line(1,0){20}}
\qbezier(210,20)(240,-10)(270,20)
\qbezier(210,20)(240,50)(270,20)
\put(270,20){\line(1,0){20}}
\qbezier(290,20)(320,-10)(350,20)
\qbezier(290,20)(320,50)(350,20)
\put(350,20){\line(1,0){20}}
\qbezier(370,20)(400,-10)(430,20)
\qbezier(370,20)(400,50)(430,20)
\put(430,20){\line(1,0){20}}
\qbezier(450,20)(480,-10)(510,20)
\qbezier(450,20)(480,50)(510,20)
\end{picture}
\end{center}
\caption{A possible configuration appearing in the second stage of the
expansion.}
\label{fig-2pt2}
\end{figure}

To continue the expansion, we would like to rewrite the final term
of \refeq{2pt.26} in terms of a product with the two-point function.
A configuration contributing to the expectation in the final term of
\refeq{2pt.26} is illustrated schematically in Figure~\ref{fig-2pt2},
in which the bonds drawn with heavy lines should be regarded as living
on a different graph than the bonds drawn with lighter lines, as explained
previously.  Our goal is to extract a factor $\tau(v',y)$, where $v'$
is shown in Figure~\ref{fig-2pt2}.

Given a configuration in which
$v \ctx{A} y $,
the {\em cutting bond}\/ $(u',v')$ is defined
to be the first pivotal bond for $v \conn y$ such that
$v \ctx{A} u'$.
It is possible that no such bond exists, as
for example would be the case in
Figure~\ref{fig-2pt2} if only the leftmost four sausages were included in
the figure (using the terminology of Section~\ref{sec-le}),
with $y$ in the location currently occupied by $u'$.
Recall the definitions of $E'(v, y; A)$
and $E(x, u, v, y; A)$
in \refeq{317}
and \refeq{317a}.
By partitioning $\{ v \ct{A} y\}$
according to the location of the cutting
bond (or the lack of a cutting bond), we obtain the partition
    \eq
    \{ v \ct{A} y\}
    = E'(v, y; A)
    \bigcup^\cdot
    \bigcup_{(u',v')}^{\cdot}
    E(v, u', v', y; A),
    \en
which implies that
\eqalign \lbeq{2pt.31a}
    \prob{v \ct{A} y}
    & =
    \prob{E'(v,y; A)} +
    \sum_{(u',v')}
    \prob{E(v, u', v', y; A)}
    .
    \enalign
Using Lemma~\ref{lem-cut1}, this gives
    \eqalign \lbeq{2pt.31}
    \prob{v \ct{A} y}
    & =
    \prob{E'(v,y; A)} +
     \sum_{(u', v')}J(u',v') \,
    \Ebold \left(  I[E'(v, u'; A)] \,
    \tau^{\tilde{C}^{(u', v')}(v)}(v', y)
    \right) .
    \enalign
Inserting the identity \refeq{2pt.23}
into \refeq{2pt.31}, we obtain
\eqalign
\lbeq{2pt.31b}
    \prob{v \ct{A} y}
    & =
    \prob{E'(v,y; A)}
    + \sum_{(u', v')} J(u',v') \,
    \prob{E'(v, u'; A)} \,
    \tau(v',y)
    \nnb
    &
     \hspace{5mm}
    -  \sum_{(u', v')} J(u',v') \,
    \Ebold_{\sss 1} \big(  I[E'(v, u'; A)] \,
    \prob[\sss 2]{v' \ctx{\tilde{C}_1^{(u', v')}(v)} y}
    \big) .
\enalign
In the last term on the right side, we have introduced subscripts for
$\tilde{C}$ and the expectations,
to indicate to which expectation $\tilde{C}$ belongs.

Let
    \eq
    \lbeq{pi1def}
    \Pi^{\smallsup{1}}(x,y) = \sum_{(u,v)} J(u,v) \,
    \Ebold_{\SSS 0} \left(  I [x \dbc u ]
    \Pbold_{\SSS 1}
    \big(E'(v,y; \tilde{C}_0^{(u, v)}(x))\big)
    \right).
    \en
Inserting \refeq{2pt.31b} into \refeq{2pt.26}, and using \refeq{pi1def},
we have
\eqalign
\lbeq{2pt.33}
    \tau(x,y)
    & = \delta_{x,y} + \Pi^{\smallsup{0}}(x,y)
    - \Pi^{\smallsup{1}}(x,y)
    + \sum_{(u,v)}
    \big( \delta_{x, u} + \Pi^{\smallsup{0}}(x,u)
    - \Pi^{\smallsup{1}}(x,u) \big)
    J(u,v) \,
    \tau(v,y)
    \nnb
    & \hspace{5mm} +  \sum_{(u,v)}J(u,v)
    \sum_{(u',v')} J(u',v')
    \nnb
    & \hspace{10mm}
    \times
    \Ebold_{\sss 0} \Big(  I [ x \dbc u ]
    \Ebold_{\sss 1} \big(
    I[E'(v, u';
    \tilde{C_{0}}^{(u, v)}(x))]
    \prob[\sss 2]{v' \ctx{\tilde{C}_{1}^{(u', v')}(v)} y}
    \big)
    \Big).
\enalign
This proves \refeq{tauxM} for $M=1$, with $R_{\SSS 1}(x,y)$ given by
the last two lines of \refeq{2pt.33}.

We now repeat this procedure recursively, rewriting
$\prob[\sss 2]{v' \ctx{\tilde{C}_{1}^{(u', v')}(v)} y}$
using \refeq{2pt.31b}, and so on. This leads to
\refeq{tauxM}, with $\Pi^\smallsup{0}$ and
$\Pi^\smallsup{1}$ given by \refeq{pi0def} and
\refeq{pi1def}, and, for $N \geq 2$,
    \eqalign
    \lbeq{2pt.38}
    \Pi^{\smallsup{N}}(x,y)
    & = 
    \sum_{(u_0, v_0)}\cdots \sum_{(u_{N-1},v_{N-1})}
    \Big[ \prod_{i=0}^{N-1} J( u_i,v_i) \Big]
    \Ebold_{\sss 0} I[x \dbc u_{0}] \, \,
    \\ \nonumber &
    \quad \quad \times
    \Ebold_{\sss 1} I[E'(v_{0}, u_{1}; \tilde{C}_{0})]
     \cdots
     \Ebold_{\sss N-1} I[E'(v_{N-2}, u_{N-1}; \tilde{C}_{N-2})]
    \Ebold_{\sss N} I[E'(v_{N-1}, y; \tilde{C}_{N-1})],\\
    R_{\SSS M}(x,y)
    & = 
    (-1)^{M+1}
    \sum_{(u_0, v_0)}\cdots \sum_{(u_{M},v_{M})}
    \Big[ \prod_{i=0}^{M} J( u_i,v_i) \Big]
    \Ebold_{\sss 0} I[x \dbc u_{0}] \, \,
    \nonumber \\ \nonumber &
    \quad \quad \times
    \Ebold_{\sss 1} I[E'(v_{0}, u_{1}; \tilde{C}_{0})]
     \cdots
     \Ebold_{\sss M-1} I[E'(v_{M-2}, u_{M-1}; \tilde{C}_{M-2})]
    \nnb &
    \quad \quad \times
    \Ebold_{\sss M} \big[I[E'(v_{M-1}, u_{M}; \tilde{C}_{M-1})]
    {\mathbb P}_{\sss M+1}(v_{M}\ctx{\tilde{C}_{M}} y)\big],
    \lbeq{Rdef}
    \enalign
where we have used the abbreviation
$\tilde{C}_{j} = \tilde{C}_j^{(u_{j}, v_{j})}(v_{j-1})$,
with $v_{-1}=x$.

Since
    \eq
    {\mathbb P}_{\sss M+1}(v_{M}\ctx{\tilde{C}_{M}} y)\leq \tau_p(v_{M},y),
    \en
it follows from \refeq{2pt.38}--\refeq{Rdef}  that
    \eq
    \lbeq{Rxbd}
    |R_{\SSS M}(x,y)|
    \leq \sum_{u_{M},v_{M} \in \ver}
    \Pi^{\smallsup{M}}(x,u_{M}) J(u_M,v_{M}) \tau_p(v_{M},y).
    \en

\section{Diagrammatic estimates for the lace expansion}
\label{sec-diagest}

In this section, we prove bounds on $\Pi^\smallsup{N}$.
These bounds are summarized in Lemma~\ref{prop-Pidiag}.
We refer to the methods of this section
as diagrammatic estimates, as we use
Feynman diagrams to provide a convenient representation for
upper bounds on $\Pi^\smallsup{N}$.

\subsection{The diagrams}

In this section, we show how $\Pi^\smallsup{N}$ of \refeq{2pt.38}
can be bounded in terms of Feynman diagrams.
Our approach here is essentially identical to what is done in
\cite[Section~2.2]{HS90a}, apart from some notational differences,
and we omit some details in the following.
The results of this section apply to any graph
$\mathbb G=(\mathbb V,\mathbb B)$,
finite or infinite,
which need not be transitive nor regular.

Given increasing
events $E,F$, we use the standard
notation $E\circ F$ to denote the event
that $E$ and $F$ occur disjointly.  Roughly speaking, $E \circ F$ is
the set of bond configurations for which there exist two disjoint sets
of occupied bonds such that the first set guarantees the occurrence of $E$
and the second guarantees the occurrence of $F$.
The BK inequality
asserts that $\Pbold (E \circ F) \leq \Pbold (E)\Pbold(F)$,
for increasing events $E$ and $F$.
(See \cite[Section~2.3]{Grim99} for a proof, and for
a precise definition of $E \circ F$.)

\begin{figure}
\vskip6.5cm
\begin{center}
\setlength{\unitlength}{0.0080in}%
\begin{picture}(400,0)
\put(-82,215){$F_0(x,u_0,w_0,z_1)=$}
\qbezier(110,220)(140,250)(170,220)
\qbezier(110,220)(140,190)(170,220)
\put(105,198){$x$}
\put(130,245){$w_{\sss 0}$}
\put(170,198){$u_{\sss 0}$}
\qbezier(140,235)(280,280)(350,265)
\put(350,280){$z_{\sss 1}$}

\put(-167,115){$F'(v_{i-1},t_i,z_i,u_i, w_i, z_{i+1})=$}
\put(100,98){$v_{\sss i-1}$}
\put(285,98){$t_{\sss i}$}
\put(315,88){$z_{\sss i}$}
\put(319,101){$\scriptstyle{\bullet}$}
\put(350,98){$u_{\sss i}$}
\put(300,140){$w_{\sss i}$}
\put(470,172){$z_{\sss i+1}$}
\put(290,120){\line(-1,0){180}}
\qbezier(290,120)(320,90)(350,120)
\qbezier(290,120)(320,150)(350,120)
\qbezier(320,135)(400,165)(470,165)

\put(-167,15){$F''(v_{i-1},t_i,z_i,u_i, w_i, z_{i+1})=$}
\put(100,-2){$v_{i-1}$}
\put(285,-2){$t_{\sss i}$}
\put(315,-12){$z_{\sss i}$}
\put(319,1){$\scriptstyle{\bullet}$}
\put(350,-2){$u_{\sss i}$}
\put(195,35){$w_{\sss i}$}
\put(390,72){$z_{\sss i+1}$}
\put(290,20){\line(-1,0){180}}
\qbezier(290,20)(320,-10)(350,20)
\qbezier(290,20)(320,50)(350,20)
\qbezier(200,20)(280,70)(390,65)

\put(-120,-85){$F_N(v_{N-1},t_N,z_N,y)=$}
\put(110,-102){$v_{\sss N-1}$}
\put(285,-102){$t_{\sss N}$}
\put(315,-112){$z_{\sss N}$}
\put(319,-99){$\scriptstyle{\bullet}$}
\put(350,-102){$y$}
\put(290,-80){\line(-1,0){180}}
\qbezier(290,-80)(320,-50)(350,-80)
\qbezier(290,-80)(320,-110)(350,-80)
\end{picture}
\end{center}
\vskip3cm

\caption{Diagrammatic representations of the events $F_0(x,u_0,w_0,z_1))$,
$F'(v_{i-1},t_i,z_i,u_i, w_i, z_{i+1})$, $F''(v_{i-1},t_i,z_i,u_i,
w_i, z_{i+1})$, $F_N(v_{N-1},t_N,z_N,y)$.
Lines indicate disjoint connections.} \label{fig-F}
\end{figure}
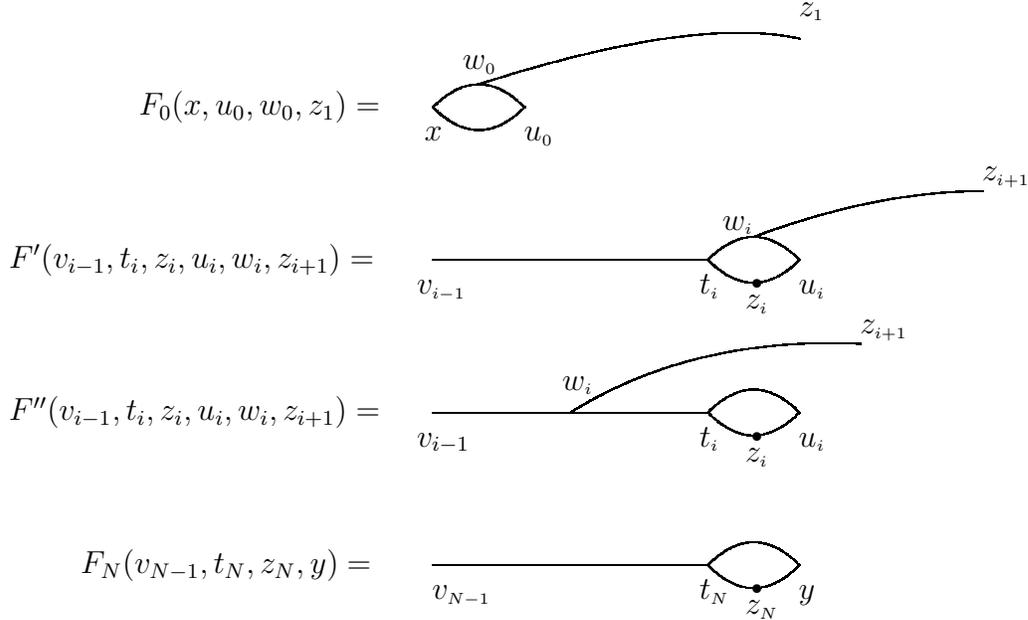

Let ${\mathbb P}^{\smallsup{N}}$ denote the product
measure on $N+1$ copies of percolation on
$\gr$.
By Fubini's Theorem
and \refeq{2pt.38},
    \eqalign
    \Pi^{\smallsup{N}}(x,y)
    &=
    \sum_{(u_0, v_0)}\cdots \sum_{(u_{N-1},v_{N-1})}
    \Big[\prod_{i=0}^{N-1} J( u_i,v_i) \Big]
    \nonumber \\
    &\qquad \times
    {\mathbb P}^{\smallsup{N}}\big(\{x \dbc u_{0}\}_0 \cap
    \big( \bigcap _{i=1}^{N-1}E'(v_{i-1}, u_{i}; \tilde{C}_{i-1})_i \big)
    \cap E'(v_{N-1}, y; \tilde{C}_{N-1})_N\big),
    \enalign
where, for an event $F$, we write $F_i$ to denote
that $F$ occurs on graph $i$.
To estimate $\Pi^\smallsup{N}(x,y)$, it is convenient to define the events
(for $N \geq 1$)
    \eqalign
    F_0(x,u_0,w_0,z_1) &= \{x\conn u_0\}\circ
    \{x\conn w_0\}\circ \{w_0\conn u_0\}
    \circ \{w_0\conn z_1\},
    \lbeq{Fdefa}\\
    \!\!\!F'(v_{i-1},t_i,z_i,u_i, w_i, z_{i+1})
    &= \{v_{i-1}\conn t_i\}\circ \{t_i\conn z_i\}\circ \{t_i\conn w_i\}\nonumber\\
    &\qquad \circ \{z_i\conn u_i\}
    \circ \{w_i\conn u_i\} \circ \{w_i\conn z_{i+1}\},
    \lbeq{F1defa}\\
    \!\!\!F''(v_{i-1},t_i,z_i,u_i, w_i, z_{i+1})
    &= \{v_{i-1}\conn w_i\}\circ \{w_i\conn t_i\}\circ \{t_i\conn z_i\}\nonumber\\
    &\qquad\circ \{t_i\conn u_i\} \circ \{z_i\conn u_i\} \circ \{w_i\conn z_{i+1}\},
    \lbeq{F2defa}\\
    \!\!\!F(v_{i-1},t_i,z_i,u_i, w_i, z_{i+1}) &=F'(v_{i-1},t_i,z_i,u_i, w_i, z_{i+1})
    \cup F''(v_{i-1},t_i,z_i,u_i, w_i, z_{i+1}),
    \lbeq{Fdefsa}
    \\
    F_N(v_{N-1},t_N,z_N,y) &= \{v_{N-1}\conn t_N\}\circ \{t_N\conn z_N\}
    \circ \{t_N\conn x\}
    \circ \{z_N\conn y\}.
    \lbeq{FNdefa}
    \enalign
The events $F_0$, $F'$,
$F''$, $F_N$ are depicted in Figure~\ref{fig-F}.  Note that
\eq
\lbeq{F0FN}
    F_N(v,t,z,y) = F_0(y,z,t,v).
\en

By the definition of $E'$ in \refeq{317},
    \eq
    \lbeq{E'bd}
    E'(v_{N-1}, y; \tilde{C}_{N-1})_N
    \subset \bigcup_{z_N\in \tilde{C}_{N-1}}
    \bigcup_{t_N\in \ver} F_N(v_{N-1},t_N,z_N,y)_N.
    \en
Indeed, viewing the connection from $v_{N-1}$ to $y$ as a string
of sausages beginning at $v_{N-1}$ and ending at $y$, for the event
$E'$ to occur there must be a vertex $z_N \in \tilde{C}_{N-1}$ that lies
on the last sausage, on a path from $v_{N-1}$ to $y$.  (In fact, both
``sides'' of the sausage must contain a vertex in $\tilde{C}_{N-1}$,
but we do not need or use this.)  This leads to \refeq{E'bd}, with
$t_N$ representing the other endpoint of the sausage that terminates at $y$.

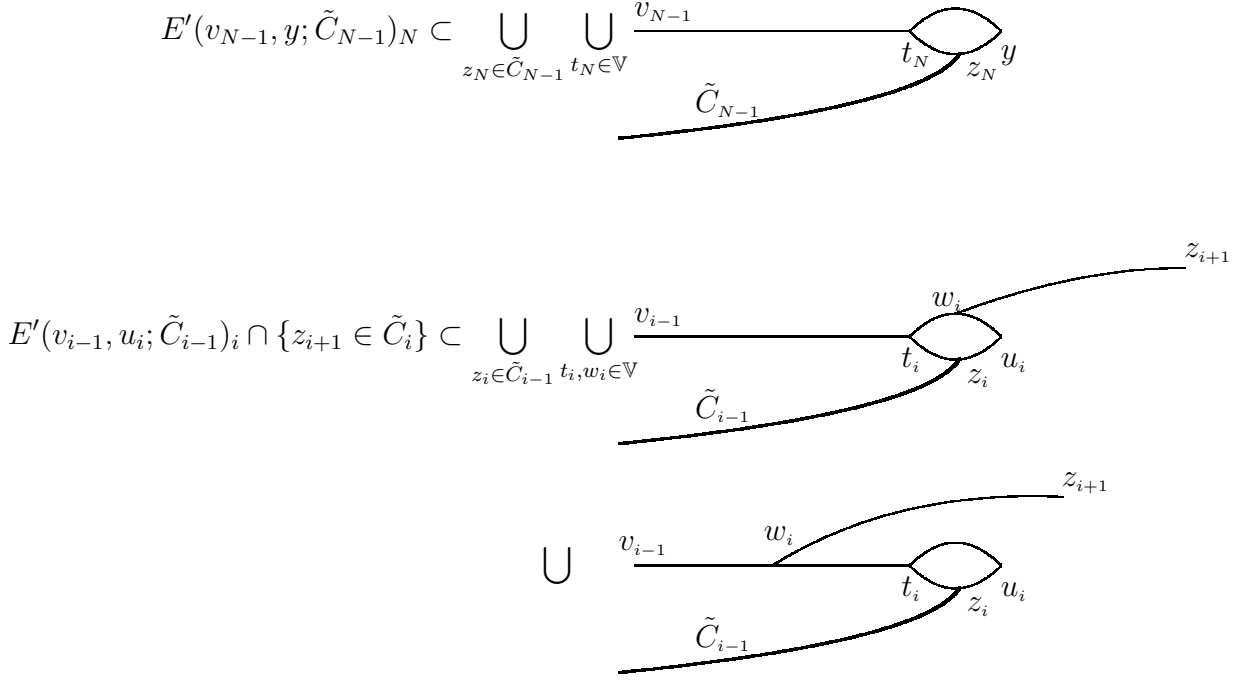
\begin{figure}
\vskip8cm
\begin{center}
\setlength{\unitlength}{0.0080in}%
\begin{picture}(600,0)
\put(0,315){$\displaystyle E'(v_{N-1},y; \tilde C_{N-1})_N
\subset \bigcup_{z_N\in \tilde C_{N-1}}\bigcup_{t_N\in \ver}$}
\put(310,330){$v_{\sss N-1}$}
\put(485,300){$t_{\sss N}$}
\put(527,291){$z_{\sss N}$}
\put(550,302){$y$}
\put(350,267){$\tilde C_{\sss N-1}$}
\put(490,320){\line(-1,0){180}}
\qbezier(490,320)(520,350)(550,320)
\qbezier(490,320)(520,290)(550,320)
\thicklines
\qbezier(300,250)(500,270)(522.5,305.5)
\thinlines

\put(-100,115){$\displaystyle E'(v_{i-1},u_i; \tilde C_{i-1})_i\cap
\{z_{i+1}\in \tilde C_{i}\}
\subset \bigcup_{z_{i}\in \tilde C_{i-1}}\bigcup_{t_i, w_i\in \ver}$}
\put(310,130){$v_{\sss i-1}$}
\put(485,100){$t_{\sss i}$}
\put(528,90){$z_{\sss i}$}
\put(550,100){$u_{\sss i}$}
\put(505,140){$w_{\sss i}$}
\put(670,172){$z_{\sss i+1}$}
\put(350,67){$\tilde C_{\sss i-1}$}
\put(490,120){\line(-1,0){180}}
\qbezier(490,120)(520,90)(550,120)
\qbezier(490,120)(520,150)(550,120)
\qbezier(520,135)(600,165)(670,165)
\thicklines
\qbezier(300,50)(500,70)(522.5,105.5)
\thinlines

\put(250,-35){$\displaystyle\bigcup$}
\put(300,-20){$v_{\sss i-1}$}
\put(485,-50){$t_{\sss i}$}
\put(528,-60){$z_{\sss i}$}
\put(550,-50){$u_{\sss i}$}
\put(395,-13){$w_{\sss i}$}
\put(590,22){$z_{\sss i+1}$}
\put(350,-83){$\tilde C_{\sss i-1}$}
\put(490,-30){\line(-1,0){180}}
\qbezier(490,-30)(520,-60)(550,-30)
\qbezier(490,-30)(520,0)(550,-30)
\qbezier(400,-30)(480,20)(590,15)

\thicklines
\qbezier(300,-100)(500,-80)(522.5,-44.5)
\thinlines
\end{picture}
\end{center}
\vskip2cm
\caption{Diagrammatic representations of the inclusions
in \refeq{E'bd} and \refeq{Eind}.} \label{fig-EsubF}
\end{figure}

Assume, for the moment, that $N \geq 2$.
The condition in \refeq{E'bd} that $z_N\in \tilde{C}_{N-1}$
is a condition on the
graph $N-1$ that must be satisfied in
conjunction with the event $E'(v_{N-2}, u_{N-1}; \tilde{C}_{N-2})_{N-1}$.
It is not difficult to see that for $i \in \{1,\ldots, N-1\}$,
\eq
\lbeq{Eind}
    E'(v_{i-1}, u_{i}; \tilde{C}_{i-1})_{i}
    \cap
    \{ z_{i+1}\in \tilde{C}_{i}\}
    \subset \bigcup_{z_{i}\in \tilde{C}_{i-1}}
    \bigcup_{t_i,w_i\in \ver}
    F(v_{i-1},t_{i},z_{i},u_{i}, w_{i}, z_{i+1})_{i}.
\en
See Figure~\ref{fig-EsubF} for a depiction of the inclusions in
\refeq{E'bd} and \refeq{Eind}.
Further details are given in \cite[Lemma~2.5]{HS90a}
or \cite[Lemma~5.5.8]{MS93}.

With an appropriate treatment for
graph $0$, \refeq{E'bd}
and \refeq{Eind} lead to
    \eqalign
    \lbeq{Fbda}
    &\{x \dbc u_{0}\}_0 \cap
        \Big(\bigcap _{i=1}^{N-1}E'(v_{i-1}, u_{i}; \tilde{C}_{i-1})_i\Big)
        \cap E'(v_{N-1}, y; \tilde{C}_{N-1})_N
    \\
    \nonumber
    &\quad \quad \subset  \bigcup_{\vec{t},\vec{w}, \vec{z}}
    \Big(
    F_0(x, u_0, w_0, z_1)_0
    \cap \big( \bigcap _{i=1}^{N-1} F(v_{i-1},t_i,z_i,u_i, w_i, z_{i+1})_i
    \big) \cap
    F_N(v_{N-1},t_N,z_N,y)_N\Big),
    \enalign
where $\vec{t}= (t_1,\ldots,t_N)$, $\vec{w}= (w_0,\ldots,w_{N-1})$
and $\vec{z}= (z_1,\ldots,z_N)$.
Therefore,
    \eqalign
    \lbeq{PiFs}
    \Pi^{\smallsup{N}}(x,y)
    & \leq
    \sum
    \left[\prod_{i=0}^{N-1}J(u_i,v_i) \right]
    \Pro (F_0(x,u_0,w_0,z_1))
    \nnb
    & \quad \quad \times
    \prod_{i=1}^{N-1} \Pro(F(v_{i-1},t_i,z_i,u_i, w_i, z_{i+1}))
    \Pro(F_N(v_{N-1},t_N,z_N,y)),
    \enalign
where the summation is over $z_1, \ldots, z_N,t_1, \ldots, t_N,
w_0, \ldots, w_{N-1}, u_0, \ldots, u_{N-1},
v_0, \ldots, v_{N-1}$.
The probability in \refeq{PiFs} factors because the
events $F_0,\ldots, F_N$ are events on different
percolation models.
Each probability in \refeq{PiFs} can be estimated using the
BK inequality.  The result is that each of the connections
$\{a \conn b\}$ present in the events $F_0$, $F$ and $F_N$ is replaced by a
two-point function $\tau_p(a,b)$.
This results in a large sum of two-point functions.

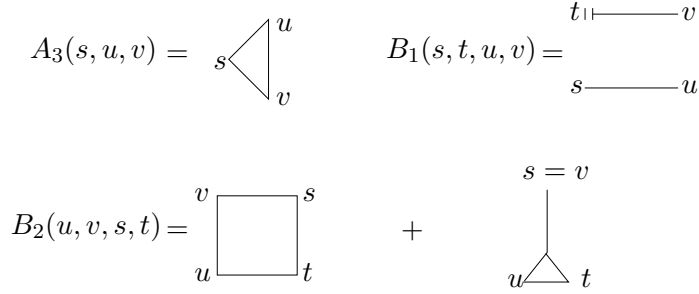
\begin{figure}[t]
\vskip1cm
\begin{center}
\setlength{\unitlength}{0.0004in}
\begingroup\makeatletter\ifx\SetFigFont\undefined%
\gdef\SetFigFont#1#2#3#4#5{%
  \reset@font\fontsize{#1}{#2pt}%
  \fontfamily{#3}\fontseries{#4}\fontshape{#5}%
  \selectfont}%
\fi\endgroup%
{
\begin{picture}(8478,3320)(0,-10)
\path(2325,3000)(2850,3525)(2850,2475)(2325,3000)
\path(7025,2625)(8250,2625)
\path(7125,3600)(8250,3600)
\path(2175,150)(3225,150)(3225,1200)(2175,1200)
\path(2175,1200)(2175,150)
\path(6207,60)(6807,60)(6507,435)
\path(6225,75)(6525,450)(6525,1275)
\path(7125,3675)(7125,3525)
\path(7025,3675)(7025,3525)
{\small
\put(2125,2875){\makebox(0,0)[lb]{$s$}}
\put(2950,3375){\makebox(0,0)[lb]{$u$}}
\put(2950,2400){\makebox(0,0)[lb]{$v$}}
\put(-300,2925){\makebox(0,0)[lb]{$A_3(s,u,v)$}}
\put(4400,2925){\makebox(0,0)[lb]{$B_1(s,t,u,v)$}}
\put(6450,2975){\makebox(0,0)[lb]{$=$}}
\put(-550,600){\makebox(0,0)[lb]{$B_2(u,v,s,t)$}}
\put(1500,650){\makebox(0,0)[lb]{$=$}}
\put(4600,600){\makebox(0,0)[lb]{$+$}}
\put(6825,2550){\makebox(0,0)[lb]{$s$}}
\put(6825,3525){\makebox(0,0)[lb]{$t$}}
\put(8300,2550){\makebox(0,0)[lb]{$u$}}
\put(8300,3525){\makebox(0,0)[lb]{$v$}}
\put(6975,0){\makebox(0,0)[lb]{$t$}}
\put(6000,0){\makebox(0,0)[lb]{$u$}}
\put(6200,1425){\makebox(0,0)[lb]{$s=v$}}
\put(3300,1125){\makebox(0,0)[lb]{$s$}}
\put(3300,75){\makebox(0,0)[lb]{$t$}}
\put(1875,75){\makebox(0,0)[lb]{$u$}}
\put(1875,1125){\makebox(0,0)[lb]{$v$}}
\put(1500,2975){\makebox(0,0)[lb]{$=$}}}
\end{picture}
}
\end{center}
\caption{Diagrammatic representations of $A_3(s,u,v)$, $B_1(s,t,u,v)$
and $B_2(u,v,s,t)$.}
\label{fig-2nc}
\end{figure}

To organize a large sum of this form, we let
    \eqalign
    \tilde{\tau}_p(x,y) &= (J*\tau_p)(x,y)
    \nnb
    & = p\cn (D*\tau_p)(x,y) \text{  if $\gr$ is regular},
    \enalign
and 
define
    \eqalign
        A_3(s,u,v)
        & =
        \tau_{p}(s,v)\tau_{p}(s,u)\tau_{p}(u,v),\lbeq{A3def}\\
        B_1(s,t,u,v) & =  \tilde\tau_{p}(t,v) \tau_p(s,u),
        \lbeq{B1def}
        \\
        B_2(u,v,s,t) & = \tau_{p}(u,v)\tau_{p}(u,t)
        \tau_{p}(v,s)\tau_{p}(s,t)
    \nonumber \\ & \quad
    + \sum_{a\in \ver}
        \tau_{p}(s,a)\tau_{p}(a,u)\tau_{p}(a,t) \delta_{v,s}\tau_{p}(u,t).
        \lbeq{B2def}
    \enalign
The two terms in $B_2$ arise from the two events $F'$ and $F''$
in \refeq{Fdefsa}. We will write them as $B_2^\smallsup{1}$ and
$B_2^\smallsup{2}$, respectively. The above quantities are
represented diagrammatically in Figure~\ref{fig-2nc}. In the
diagrams, a line joining $a$ and $b$ represents $\tau_p(a,b)$. In
addition, small bars are used to distinguish a line that
represents $\tilde{\tau}_p$, as in $B_1$.

Application of the BK
inequality yields
    \eqalign
    \Pro (F_0(x,u_0,w_0,z_1)) &\leq
    A_3(x,u_0,w_0)\tau_p(w_0,z_1),
    \lbeq{F0bd1}\\
    \sum_{v_{N-1}} J(u_{N-1},v_{N-1}) \Pro(F_N(v_{N-1},t_N,z_N,y))&\leq
    \frac{B_1(w_{N-1},u_{N-1},z_N,t_N)}{\tau_p(w_{N-1},z_N)}
    A_3(y,t_N,z_N).\lbeq{F0bd2}
    \enalign
For $F'$ and $F''$,  application of the BK inequality yields
    \eqalign
    &\sum_{v_{i-1}}J(u_{i-1},v_{i-1})\Pro(F'(v_{i-1},t_i,z_i,u_i, w_i,
    z_{i+1})) \nonumber\\
    &\qquad \qquad \leq
    \frac{B_1(w_{i-1}, u_{i-1}, z_i, t_i)}{\tau_p(w_{i-1},z_i)}
    B_2^{\smallsup{1}}(z_i,t_i, w_{i},u_i)
    \tau_p(w_i,z_{i+1}),
    \lbeq{F1bd}\\
    &\sum_{v_{i-1}, t_i} J(u_{i-1},v_{i-1})
    \Pro(F''(v_{i-1}, t_i ,z_i,u_i, w_i,z_{i+1}))\nonumber\\
    &\qquad \qquad \leq
    \frac{B_1(w_{i-1}, u_{i-1}, z_i, w_i)}{\tau_p(w_{i-1},z_i)}
     B_2^{\smallsup{2}}(z_i,w_i,
    w_{i}, u_i) \tau_p(w_i,z_{i+1}).\lbeq{F2bd}
    \enalign
Since the second and the third arguments of
$B_2^{\smallsup{2}}$ are equal by virtue of the Kronecker delta
in \refeq{B2def},
we can combine \refeq{F1bd}--\refeq{F2bd} to obtain
    \eqalign
    &\sum_{v_{i-1}, t_i} J(u_{i-1},v_{i-1})
    \Pro(F(v_{i-1},t_i,z_i,u_i, w_i,
    z_{i+1})) \nonumber\\
    &\qquad \qquad \leq \sum_{t_i}
    \frac{B_1(w_{i-1}, u_{i-1}, z_i, t_i)}{\tau_p(w_{i-1},z_i)}
     B_2(z_i,t_i,
    w_{i},u_i) \tau_p(w_i,z_{i+1}).\lbeq{Fbd}
    \enalign

\begin{figure}
\begin{center}
\setlength{\unitlength}{0.008in}%
\begin{picture}(400,130)(20,635)
\thinlines
\thicklines
\put(105,765){\line( 0,-1){ 10}}
\put(100,765){\line( 0,-1){ 10}}
\thinlines
\put(105,760){\line( 1, 0){ 35}}
\put(100,720){\line( 1, 0){ 40}}
\thicklines
\put(100,685){\line( 0,-1){ 10}}
\put(105,685){\line( 0,-1){ 10}}
\thinlines
\put(105,680){\line( 1, 0){115}}
\put(140,680){\line( 0,-1){ 40}}
\put(180,680){\line( 0,-1){ 40}}
\thicklines
\put(180,645){\line( 0,-1){ 10}}
\put(185,645){\line( 0,-1){ 10}}
\thinlines
\put(185,640){\line( 1, 0){ 35}}
\put(100,640){\line( 1, 0){ 80}}
\put(320,640){\line( 1, 0){ 30}}
\put(360,680){\line( 0,-1){ 30}}
\thicklines
\put(320,685){\line( 0,-1){ 10}}
\put(325,685){\line( 0,-1){ 10}}
\thinlines
\put(325,680){\line( 1, 0){ 75}}
\put(370,645){\line( 0, 1){  0}}
\thicklines
\put(370,645){\line( 0,-1){ 10}}
\put(375,645){\line( 0,-1){ 10}}
\thinlines
\put(375,640){\line( 1, 0){ 25}}
\put(140,720){\line( 1, 1){ 20}}
\put(160,740){\line( 0, 1){  0}}
\put(160,740){\line(-1, 1){ 20}}
\put(140,760){\line( 0,-1){ 40}}
\put(100,760){\line( 0,-1){ 40}}
\put(100,720){\line(-1, 1){ 20}}
\put( 80,740){\line( 1, 1){ 20}}
\put(100,680){\line( 0,-1){ 40}}
\put(100,640){\line(-1, 1){ 20}}
\put( 80,660){\line( 1, 1){ 20}}
\put(220,680){\line( 0,-1){ 40}}
\put(220,640){\line( 1, 1){ 20}}
\put(240,660){\line(-1, 1){ 20}}
\put(360,650){\line(-1,-1){ 10}}
\put(350,640){\line( 1, 0){ 20}}
\put(360,650){\line( 1,-1){ 10}}
\put(400,680){\line( 0,-1){ 40}}
\put(400,640){\line( 1, 1){ 20}}
\put(420,660){\line(-1, 1){ 20}}
\thinlines
\put(320,680){\line( 0,-1){ 40}}
\put(320,640){\line(-1, 1){ 20}}
\put(300,660){\line( 1, 1){ 20}}
\thinlines
\put( 20,750){\makebox(0,0)[lb]{\raisebox{0pt}[0pt][0pt]{(a)}}}
\put( 20,665){\makebox(0,0)[lb]{\raisebox{0pt}[0pt][0pt]{(b)}}}
\put( 60,740){\makebox(0,0)[lb]{\raisebox{0pt}[0pt][0pt]{${x}$}}}
\put( 60,660){\makebox(0,0)[lb]{\raisebox{0pt}[0pt][0pt]{$x$}}}
\put(285,660){\makebox(0,0)[lb]{\raisebox{0pt}[0pt][0pt]{$x$}}}
\put(260,655){\makebox(0,0)[lb]{\raisebox{0pt}[0pt][0pt]{$+$}}}
\put( 170,740){\makebox(0,0)[lb]{\raisebox{0pt}[0pt][0pt]{$y$}}}
\put( 240,670){\makebox(0,0)[lb]{\raisebox{0pt}[0pt][0pt]{$y$}}}
\put(430,660){\makebox(0,0)[lb]{\raisebox{0pt}[0pt][0pt]{$y$}}}
\end{picture}
\end{center}
\caption{\label{fig-pidiag} The diagrams bounding (a) $\Pi^{(1)}(x,y)$
and (b) $\Pi^{(2)}(x,y)$.   }
\end{figure}
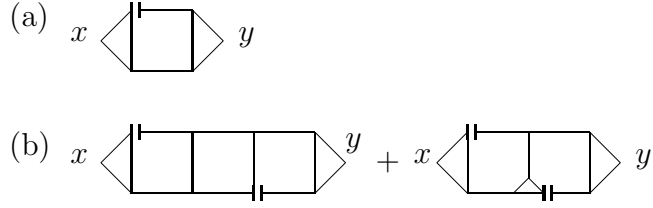

Upon substitution of the bounds on the probabilities in \refeq{F0bd1},
\refeq{F0bd2} and \refeq{Fbd} into \refeq{PiFs}, the ratios of
two-point functions form a telescoping product that disappears.
After relabelling the summation indices,
\refeq{PiFs} becomes
    \eqalign
    \lbeq{PibdAB}
    \Pi^{\smallsup{N}}(x,y)
    & \leq \sum_{\vec{s}, \vec{t},\vec{u},\vec{v}}
    A_3(x,s_1, t_1)
    \prod_{i=1}^{N-1}\big[ B_1(s_i, t_i, u_i, v_i)
    B_2(u_{i}, v_{i}, s_{i+1}, t_{i+1})\big]
    \nonumber \\ & \quad\quad \quad \times
    B_1(s_N,t_N,u_N,v_N)
    A_3(u_N, v_N, y).
    \enalign
The bound \refeq{PibdAB} is valid for $N \geq 1$, and the summation
is over all $s_1,\ldots, s_N$, $t_1, \ldots, t_N$,
$u_1,\ldots, u_N$, $v_1,\ldots ,v_N$.
For $N=1,2$, the right side is represented diagrammatically in
Figure~\ref{fig-pidiag}.  In the diagrams, unlabelled vertices are summed
over $\ver$.

\subsection{The diagrammatic bounds}
\label{sec-dbd}

We now specialize to the case $\gr = \torus_{r,n}$, making use of
the additive structure and the $x \mapsto -x$ symmetry
of the torus.
We will write $\tau_p(y-x)$ in place of $\tau_p(x,y)$,
$p\cn D(y-x)$ in place of $J(x,y)$, and
$\Pi^\smallsup{N}(y-x)$ in place of $\Pi^\smallsup{N}(x,y)$.

The upper bounds we prove are in terms of various
quantities related to the triangle diagram.  Let
        \eq
        \lbeq{T(x)def}
        T_p(x)= \sum_{y,z, u\in \ver}
    \tau_p(y)\tau_p(z-y)p\cn D(u)\tau_p(x+z-u)
        = (\tau_p * \tau_p * \tilde\tau_p )(x),
        \en
        \eq
        \lbeq{Tdef}
        T_p=\max_{x\in \ver} T_p(x),
        \en
        \eq
        \lbeq{T'def}
        T_p'= \max_{x \in \ver}
        \sum_{y,z \in \ver} \tau_p(y)\tau_p(z-y)\tau_p(x-z)
        =\max_{x\in \ver} (\tau_p*\tau_p*\tau_p)(x),
        \en
and, for $k \in \torus_{r,n}^*$, let
       \eq
        \lbeq{Wpydef}
        W_p(y;k) = \sum_{x \in \ver}[1-\cos(k\cdot x)]
        \tilde\tau_p(x) \tau_p(x+y),
        \en
        \eq
        \lbeq{Wpdef}
        W_p(k) = \max_{y \in \ver} W_p(y;k).
        \en
Recall that $B_2^\smallsup{2}$ denotes the second term of
\refeq{B2def}.  For $k \in \torus_{r,n}^*$, we also define
    \eqalign
    H_p(a_1,a_2;k)
    &  =\sum_{u, v, s, t}
    [1-\cos (k \cdot (t-u))]
    B_1(0,a_1, u, s)
    B_2^\smallsup{2}(u, s, s, t)B_1(s, t, v, v+a_2),
    \lbeq{Masdef}
    \enalign
and
    \eq
    \lbeq{Hpdef}
        H_p(k) = \max_{a_1,a_2 \in \ver}H_p(a_1,a_2;k).
    \en
The remainder of this section is devoted to the proof of the following
proposition.

\begin{prop}
\label{prop-Pidiag}
For $N=0$,
    \eq
    \lbeq{Pi0bd}
    \sum_{x \in \ver} \Pi^\smallsup{0}(x)  \leq T_p,
    \en
    \eq
    \lbeq{Pi0x2}
    \sum_{x \in \ver} [1-\cos (k \cdot x)] \Pi^\smallsup{0}(x)  \leq W_p(0;k).
    \en
For $N \geq 1$,
    \eq
    \lbeq{PiNbd}
    \sum_{x \in \ver} \Pi^\smallsup{N}(x)  \leq T_p'(2T_pT_p')^{N},
    \en
    \eqalign
    \lbeq{PiNbdx2}
    \sum_{x \in \ver} [1-\cos (k \cdot x)] \Pi^\smallsup{N}(x)  &\leq
    (4N+3) \Big[
    T_p'W_p(k)\big(2T_p+[1+p\cn] N T_p'\big)(2T_pT_p')^{N-1}
    \nnb & \quad\quad \quad\quad\quad
    + (N-1)\big(T_p^2W_p(k)+H_p(k)\big)(T_p')^2 (2T_pT_p')^{N-2} \Big],
    \enalign
and, for $N=1$, \refeq{PiNbdx2} can also be replaced by
\eqalign
    \sum_{x \in \ver} [1-\cos(k\cdot x)] \Pi^\smallsup{1}(x)
    &\leq W_p(0;k)+ 31 T_pT_p' W_p(k).
    \lbeq{Pi1cos}
    \enalign
\end{prop}

\subsubsection{Proof of \refeq{Pi0bd}--\refeq{Pi0x2}}

By \refeq{pi0def} and the BK inequality,
        \eq
        \Pi^{\smallsup{0}}(x)
    = \Pbold(0 \dbc x) - \delta_{0,x}
    \leq \tau_p(x)^2 - \delta_{0,x}.
        \en
For $x \neq 0$, the event $\{0 \conn x\}$ is the union over neighbors $y$
of the origin of $\{\{0,y\} \text{ occupied}\} \circ \{y \conn x\}$.
Thus, by the
BK inequality,
\eq
\lbeq{tauDtau}
    \tau_p(x) \leq p\cn (D*\tau_p)(x) = \tilde\tau_p(x)
    \quad
    ( x \neq 0).
\en
Therefore,
\eq
    \sum_{x \in \ver} \Pi^\smallsup{0}(x) \leq
    \sum_{x\in \ver }  \tau_p(x)\tilde\tau_p(x)
    \leq
    T_p(0).
\en
Similarly,
\eq
    \sum_{x \in \ver}[1-\cos (k \cdot x)] \Pi^\smallsup{0}(x)
    \leq W_p(0;k).
\en
This proves \refeq{Pi0bd}--\refeq{Pi0x2}.

%

\subsubsection{Proof of \refeq{PiNbd}}

For $N \geq 1$, let
    \eq
    \Psi^{\smallsup{N}}(s_{N+1},t_{N+1})
    =\sum_{\vec{s}, \vec{t},\vec{u},\vec{v}}
    A_3(0,s_1, t_1)
    \prod_{i=1}^N \big[ B_1(s_i, t_i, u_i, v_i)
    B_2(u_{i}, v_{i}, s_{i+1}, t_{i+1})
    \big].
    \lbeq{Psidef}
    \en
For convenience, we define $\Psi^{\smallsup{0}}(x,y)=A_3(0,x,y)$, so that
    \eq
    \lbeq{recPsi}
    \Psi^{\smallsup{N}}(x,y)=\sum_{u_N, v_N, s_N, t_N}
    \Psi^{\smallsup{N-1}}(s_N, t_N)
    B_1(s_N, t_N, u_N, v_N) B_2(u_N, v_N, x, y) \quad (N\geq 1).
    \en
Since
    \eq
    \sum_x A_3(u_N, v_N, x) \leq \sum_{x,y} B_2(u_N, v_N, x,y),
    \en
it follows from \refeq{PibdAB} that
    \eq
    \lbeq{PivsPsi}
    \sum_x \Pi^{\smallsup{N}}(x) \leq \sum_{x,y} \Psi^{\smallsup{N}}(x,y),
    \en
and bounds on $\Pi^{\smallsup{N}}$ can be obtained from bounds on
$\Psi^{\smallsup{N}}$.
We prove bounds on $\Psi^{\smallsup{N}}$, and hence on
$\Pi^\smallsup{N}$, by induction on $N$.

The induction hypothesis is that
    \eq
    \lbeq{indaim}
    \sum_{x,y} \Psi^{\smallsup{N}}(x,y) \leq T_p'(2 T_pT_p')^N.
    \en
For $N=0$, \refeq{indaim} is true since
    \eq
    \sum_{x,y} A_3(0, x, y) \leq T_p'.
    \en
If we assume \refeq{indaim} is valid for $N-1$, then by \refeq{recPsi},
    \eq
    \sum_{x,y} \Psi^{\smallsup{N}}(x,y)
    \leq \Big(\sum_{s_N, t_N}
    \Psi^{\smallsup{N-1}}(s_N, t_N)\Big)
    \Big(\max_{s_N, t_N} \sum_{u_N, v_N, x,y}
    B_1(s_N, t_N, u_N, v_N) B_2(u_N, v_N, x, y)\Big),
    \en
and \refeq{indaim} then follows once we prove that
    \eq
    \max_{s,t} \sum_{u,v, x,y} B_1(s, t, u, v) B_2(u, v, x, y) \leq 2T_pT_p'.
    \lbeq{bdBi1}
    \en

It remains to prove \refeq{bdBi1}.
There are two terms, due to the two terms in \refeq{B2def}, and
    we bound each term separately. The first term is bounded as
    \eqalign
        &\max_{s,t} \sum_{u,v, x,y}
        \tilde\tau_{p}(v-t) \tau_p(u-s)
        \tau_{p}(y-u)\tau_{p}(x-v)\tau_{p}(v-u)\tau_{p}(x-y)
    \nonumber \\
            &\quad= \max_{s,t} \sum_{u,v}
            \tilde\tau_{p}(v-t) \tau_p(u-s) \tau_{p}(v-u)
            \big(\sum_{x,y} \tau_{p}(y-u)\tau_{p}(x-v)\tau_{p}(x-y)\big)
            \nonumber\\
            &\quad \leq T_p' \max_{s,t} \sum_{u,v}
            \tilde\tau_{p}(v-t) \tau_p(u-s) \tau_{p}(v-u)
            \nonumber\\
            &\quad =T_p T_p'.
            \enalign
    The second term is bounded similarly, making use of translation invariance,
    by
    \eqalign
        &\max_{s,t} \sum_{u,v, x,y, a}
        \tilde\tau_{p}(v-t) \tau_p(u-s)\delta_{v,x}\tau_{p}(y-u)
            \tau_{p}(x-a)\tau_{p}(u-a)\tau_{p}(y-a)\nonumber \\
            &\quad= \max_{s,t}
        \sum_{a,y,u} \big(\tilde\tau_{p}*\tau)(a-t) \tau_p(u-s)\big)
            \big(\tau_{p}(y-u)\tau_{p}(u-a)\tau_{p}(y-a)\big)\nonumber\\
            &\quad =\max_{s,t}
        \sum_{y', a'} T_p(a'+s-t) \tau_{p}(y')\tau_{p}(a')
        \tau_{p}(y'-a')
            \nonumber\\
            &\quad \leq \big(\max_{a',s,t} T_p(a'+s-t)\big)
            \big(\sum_{y', a'} \tau_{p}(y')\tau_{p}(a')\tau_{p}(y'-a')\big)
            \nnb &
            \quad
            \leq T_p T_p',
            \enalign
    where $a'=a-u$, $y'=y-u$.
This completes the proof of \refeq{bdBi1} and hence of \refeq{PiNbd}.

\subsubsection{Proof of \refeq{PiNbdx2}}
\label{sec-bdcos}

Next, we estimate $\sum_x [1-\cos (k \cdot x)] \Pi^{\smallsup{N}}(x)$.
In a term in \refeq{PibdAB}, there is a sequence of $2N+1$ two-point
functions along the ``top'' of the diagram, such that the sum of the
displacements of these two-point functions is exactly equal to $x$.
For example, in Figure~\ref{fig-pidiag}(a) there are three displacements
along the top of the diagram, and in
Figure~\ref{fig-pidiag}(b) there are five in the first diagram
and four in the second.  We regard the second diagram as also having
five displacements, with the understanding that the third is constrained
to vanish.  With a similar general convention,
each of the $2^{N-1}$ diagrams bounding $\Pi^\smallsup{N}$ has
$2N+1$ displacements along the top of the diagram.
We denote these displacements by $d_1, \ldots, d_{2N+1}$, so that
$ x=\sum_{j=1}^{2N+1} d_j$.  We will argue as follows to distribute
the factor $1-\cos(k\cdot x)$ among the displacements $d_j$.

Let $t=\sum_{j=1}^Jt_j$.  Taking the real part of the telescoping sum
    \eq
    \lbeq{telexp1}
    1 - e^{i t} = \sum_{j =1}^{J} [1 - e^{i  t_j}]
    e^{i  \sum_{m=1}^{j-1} t_m}
    \en
leads to the bound
    \eq
    1-\cos t \leq
    \sum_{j =1}^{J} [1 - \cos t_j]
    +
    \sum_{j =1}^{J}
    \sin t_j\, \sin \left(\sum_{m=1}^{j-1} t_m \right).
    \en
It follows from the identity $\sin(x+y)=\sin x\cos y+\cos x\sin y$
that $|\sin(x+y)| \leq |\sin x|+|\sin y|$.  Applying this recursively gives
    \eq
    1-\cos t\leq \sum_{j =1}^{J} [1 - \cos t_j] +
    \sum_{j =1}^{J}\sum_{m=1}^{j-1} |\sin t_j||\sin t_m|.
    \en
In the last term we use  $|ab|\leq (a^2+b^2)/2$,
and then $1-\cos^2 a \leq 2[1-\cos a]$, to obtain
    \eqalign
    1-\cos t &\leq \sum_{j =1}^{J} [1 - \cos t_j] +
    \frac 12 \sum_{j =1}^{J}\sum_{m=1}^{j-1}
    [\sin^2 t_j +\sin^2 t_m ] \nnb
    &\leq \sum_{j =1}^{J} [1 - \cos t_j] +J
     \sum_{j =1}^{J}\sin^2 t_j\nonumber
     \\
    &= \sum_{j =1}^{J} [1 - \cos t_j] +J
    \sum_{j =1}^{J}[1-\cos^2 t_j]\nonumber
    \\
    \lbeq{normbd}
    &\leq (2J+1)\sum_{j =1}^{J} [1 - \cos t_j].
    \enalign
We apply \refeq{normbd} with $t= k\cdot x = \sum_{j=1}^{2N+1} k \cdot d_j$
to obtain a sum of $2N+1$ diagrams like the ones
for $\Pi^{\smallsup{N}}(x)$, except now in the $j^{\rm th}$
term, the $j^{\rm th}$ line
in the top of the diagram represents $[1 - \cos(k\cdot d_j)]\tau_p(d_j)$
rather than $\tau_p(d_j)$.

We distinguish three cases:
(a) the displacement $d_j$ is in a line of $A_3$,
(b) the
displacement  $d_j$ is in a line of $B_1$, (c) the
displacement  $d_j$ is in a line of $B_2$.

\smallskip
\noindent {\em Case (a): the displacement is in a line of $A_3$.}
We consider the case where the weight
$[1 - \cos(k\cdot d_j)]$ falls on the last of the
factors $A_3$ in \refeq{PibdAB}.
This contribution is equal to
\eq
\lbeq{x2casea}
    \sum_{u,v} \Psi^{\smallsup{N-1}}(u,v) \sum_{w,x,y}
    B_1(u,v,w,y)\tau_p(y-w)\big[1-\cos\big(k\cdot(x-y)\big)\big]
    \tau_p(x-y) \tau_p(x-w).
\en
Applying \refeq{tauDtau} to $\tau_p(x-y)$, we have
\eq
    \max_{u,v}\sum_{w,x,y}
    B_1(u,v,w,y)\tau_p(y-w)\big[1-\cos\big(k\cdot(x-y)\big)\big]
    \tau_p(x-y) \tau_p(x-w)
    \leq
    T_p W_p(k).
\en
It then follows from \refeq{indaim} that \refeq{x2casea} is bounded above
by $T_p'(2T_pT_p')^{N-1} T_pW_p(k)$.
By symmetry,
the same bound applies when the
weight falls into the first
factor of $A_3$, i.e, when we have a factor
$[1-\cos (k\cdot d_1)]$.  Thus case~(a) leads to an upper bound
\eq
\lbeq{caseabd}
    2T_p'(2T_pT_p')^{N-1} T_pW_p(k).
\en

\smallskip
\noindent
{\em Case (b): the displacement is in a line of  $B_1$.}
Suppose that the factor $[1-\cos(k\cdot d_j)]$
falls on the $i^{\rm th}$ factor $B_1$ in
\refeq{PibdAB}. Depending on $i$, it falls
either on $\tilde \tau_p$ or on $\tau_p$ in \refeq{B1def}.
We write the right side of \refeq{PibdAB} with the extra factor as
    \eq
    \lbeq{x2caseb}
    \sum_x \sum_{s, t, u, v} \Psi^{\smallsup{i-1}}(s, t)
    \tilde B_1(s, t, u, v) \bar\Psi^{\smallsup{N-i}}(u-x, v-x).
    \en
In \refeq{x2caseb}, either
    \eq
    \lbeq{B1tila}
    \tilde B_1(s, t, u, v)= \big[1-\cos\big(k\cdot (u-s)\big)\big]
    \tilde \tau_p(u-s)
    \tau_p(v-t)
    \en
or
    \eq
    \lbeq{B1tilb}
    \tilde B_1(s, t, u, v)= \tilde \tau_p(u-s)
    \big[1-\cos\big(k\cdot (v-t)\big)\big]
    \tau_p(v-t),
    \en
and $\bar\Psi^{\smallsup{N-i}}$ denotes a small variant of
$\Psi^{\smallsup{N-i}}$, defined inductively by
$\bar\Psi^{\smallsup{0}}=\Psi^{\smallsup{0}}$ and
$\bar\Psi^{\smallsup{i}}(x,y)
=\sum_{s,t,u,v}B_2(x,y,s,t)B_1(s,t,u,v) \bar\Psi^{\smallsup{i-1}}(u,v)$.
It can be verified
that $\bar\Psi^{\smallsup{N-i}}$ also obeys \refeq{indaim}.

For \refeq{B1tila}, we let $a_1=t-s$, $a_2=v-u$, and $x'=u-x$.
With this notation, the contribution to \refeq{x2caseb}
due to \refeq{B1tila} is bounded above by
    \eqalign
    &\Big(\sum_{s, a_1} \Psi^{\smallsup{i-1}}(s, s+a_1) \Big)
     \Big(\sum_{x', a_2}\bar\Psi^{\smallsup{N-i}}(x', x'+a_2)\Big)
     \Big(\max_{s, a_1, a_2}\sum_{u}
     \tilde B_1(s, s+a_1, u, u+a_2)\Big)\nonumber\\
     &\quad =\Big(\sum_{s, t} \Psi^{\smallsup{i-1}}(s, t) \Big)
     \Big(\sum_{x,y}\bar\Psi^{\smallsup{N-i}}(x,y)\Big)
     W_p(k)
     \nonumber\\
     &\quad \leq T_p' (2T_p T_p')^i T_p' (2T_p T_p')^{N-i-1} W_p (k)
     = T_p' (2T_p T_p')^{N-1} T_p' W_p(k),
    \enalign
where we used \refeq{indaim}.
For \refeq{B1tilb}, we use \refeq{tauDtau} for $\tau_p(v-t)$,
write $\tilde \tau_p(u-s) = \sum_{y}p\cn D(y) \tau_{p}(u-s-y)$,
estimate the sum over $y$ with a supremum and use $\sum_y p\cn D(y)
=p\cn$.
Since there are $N$ choices of
factors $B_1$, case~(b) leads to an overall upper bound
    \eq
    \lbeq{casebbd}
    N[1+p\cn]T_p' (2T_p T_p')^{N-1} T_p' W_p(k).
    \en

\smallskip
\noindent
{\em Case (c): the displacement is in a line of $B_2$.}
It is sufficient to estimate
    \eqalign
    &
    \max \hspace{-4mm}
    \sumtwo{a, b, u, v}{s, t, w, y, x}
    \shift \shift
    \Psi^{\smallsup{i-1}}(a, b)
    \bar\Psi^{\smallsup{N-i-1}}(w-x, y-x) [1 - \cos(k\cdot d)]
    \nnb
    \lbeq{Hcasec}
    &\hspace{25mm} \times
    B_1(a, b, u, v)
    B_2(u, v, s, t)B_1(s, t, w, y),
    \enalign
where the maximum is over the choices $d=s-v$ or
$d=t-u$.
We consider separately the contributions due to $B_2^\smallsup{1}$
and $B_2^\smallsup{2}$ of \refeq{B2def}, beginning with $B_2^\smallsup{2}$.

Recall the definition of $H(a_1,a_2;k)$ in
\refeq{Masdef}.
The contribution to \refeq{Hcasec} due to $B_2^\smallsup{2}$
can be rewritten, using $x'=w-x$, $a_2=y-w$, $a_1=b-a$, as
    \eqalign
    &\sum_{a, a_1, a_2, x'}
    \Psi^{\smallsup{i-1}}(a, a+a_1)
    \bar\Psi^{\smallsup{N-i-1}}(x', x'+a_2)H(a_1,a_2;k)
    \nonumber \\
    &\qquad \leq H_p (k)
    \Big(\sum_{x,y} \Psi^{\smallsup{i-1}}(x,y)\Big)
    \Big(\sum_{x,y} \bar\Psi^{\smallsup{N-i-1}}(x,y)\Big)
    \nonumber \\
    & \qquad \leq
    H_p(k)(T_p')^2 (2T_pT_p')^{N-2}
    \lbeq{bdinM}.
    \enalign
Since there are $N-1$ factors $B_2$ to choose,
this contribution to case~(c) contributes at most
    \eq
    \lbeq{casecbd}
    (N-1)H_p(k) (T_p')^2 (2T_pT_p')^{N-2}.
    \en
It is not difficult to check that the contribution to case~(c) due
to $B_2^\smallsup{1}$ is at most
        \eq
    \lbeq{casecbdz}
    (N-1)(T_p^2W_p(k))(T_p')^2 (2T_pT_p')^{N-2}.
    \en

The desired estimate \refeq{PiNbdx2} then follows from \refeq{normbd},
\refeq{caseabd}, \refeq{casebbd} and \refeq{casecbd}--\refeq{casecbdz}.

\subsubsection{Proof of \refeq{Pi1cos}}

Recall from \refeq{PibdAB} that
    \eq
    \Pi_p^{\smallsup{1}}(x) \leq \sum_{s,t,u,v\in \ver} A_3(0,s, t)B_1(s,t,u,v)
    A_3(u, v, x).
    \en
We define $A_3'(u, v, x)$ by
    \eq
    A_3'(u, v, x) =  A_3(u, v, x) - \delta_{u,x}\delta_{v,x} .
    \en
Then we have
    \eqalign
    \sum_{x\in \ver}[1-\cos(k\cdot x)]
    \Pi_p^{\smallsup{1}}(x) &\leq \sum_{x\in \ver}[1-\cos(k\cdot x)] B_1(0,0,x,x)\nnb
    &\qquad +\sum_{x,s,t,u,v\in\ver}
    [1-\cos(k\cdot x)] A_3'(0,s, t)B_1(s,t,u,v)
    A_3(u, v, x)\nnb
    &\qquad +\sum_{x,u,v\in \ver}
    [1-\cos(k\cdot x)] B_1(0,0,u,v)A_3'(u, v, x).
    \enalign
The first term equals $W_p(0;k)$.
The second and third terms are bounded above by $7\cdot 3T_pT_p' W_p(k)$
and $5\cdot 2T_p W_p(k) \leq 10T_pT_p' W_p(k)$, respectively,
using \refeq{normbd} (with $J=3$ and $J=2$)
and the methods of Section~\ref{sec-bdcos}.

\medskip
This completes the proof of Proposition~\ref{prop-Pidiag}.
\qed

\section{Analysis of the lace expansion}
\label{sec-Pibds}

In this section, we use the lace expansion to prove the triangle condition
of Theorem~\ref{thm-tc}. The analysis is similar in spirit to the analysis
of \cite{HS90a}, but it has been simplified and reorganized,
and it differs significantly in detail from the presentation of \cite{HS90a}.
Specific improvements include:
(i)~We have reduced the number of functions in the bootstrap argument
from five to three (cf.\ \cite[Proposition~4.3]{HS90a}),
and in the bootstrap we work
directly with the Fourier transform
of the two-point fuction rather than with the triangle and related diagrams.
(ii)~We work with $1-\cos(k\cdot x)$ directly, rather than expanding the
cosine to second order.
(iii)~Our treatment of $H_p(k)$ in Lemma~\ref{lem-Mbd} below is simpler
than the corresponding treatment of \cite[Section~4.4.3(e)]{HS90a}.

We work in this section on an arbitrary torus $\torus_{r,n}$ with
$r \geq 2$, assuming that Assumption~\ref{ass-rw} is satisfied.
As usual, we write the degree of the torus as $\cn$, and we abbreviate
$p_c(\torus_{r,n})$ to $p_c$.

Our analysis actually uses a slightly weaker assumption than the one
stated in Assumption \ref{ass-rw}. Instead of \refeq{rwbd},
we will assume in the proof that
    \eq
    \lbeq{rwbd2}
    \frac{1}{V}
    \sum_{k \in \torus_{r,n}^* : \; k \neq 0}
    \frac{\hat{D}(k)^{2}}{[1-\mu\hat{D}(k)]^3} \leq
    \beta
    \en
holds uniformly in $\mu \in [0, 1- \frac 12 \lambda^{-1}V^{-1/3}]$.
Equation~\refeq{rwbd2} is strictly weaker than \refeq{rwbd},
but not in a significant way.
The analogue of \refeq{rwbdi0} with $\mu$ inserted in the denominator
follows from \refeq{rwbd2} in the same way that \refeq{rwbdi0}
follows from \refeq{rwbd}.

\subsection{The bootstrap argument}

Taking the Fourier transform of \refeq{tauxM} and solving for $\hat{\tau}_p(k)$
gives
        \eq
        \lbeq{tauMk}
        \hat \tau_p(k) = \frac{1+\hat \Pi_{\SSS M}(k)+\hat R_{\SSS M}(k)}
        {1- p\cn \hat{D}(k)[1+\hat \Pi_{\SSS M}(k)]},
        \en
for all $k\in \torus^*_{r,n}$ and all $M=0,1,2,\ldots$.
Recall from \refeq{Cdef} that
$\hat{C}_{\mu}(k)=[1- \mu\cn \hat{D}(k)]^{-1}$.
As explained in Section~\ref{sec-le},
we would like to compare $\hat{\tau}_p(k)$ with
$\hat{C}_\mu(k)$, with
$\mu\cn$ equal to $p\cn [1+\hat \Pi_{\SSS M}(0)]$.
We know that $\hat{\tau}_p(0)=\chi(p) >0$, but we do not yet know
that $1+\hat \Pi_{\SSS M}(0)+\hat R_{\SSS M}(0)$ is positive and thus
we cannot yet be sure that the denominator of \refeq{tauMk} is positive
when $k=0$.
We therefore do not yet know
that our choice of $\mu$ is less than $\cn^{-1}$.
To safeguard against the possibility that
$p\cn [1+\hat \Pi_{\SSS M}(0)] \geq 1$ or
$p\cn [1+\hat \Pi_{\SSS M}(0)] <0$,
we define $\mu_p^{\smallsup{M}}$ by
        \eq
        \lbeq{mudef}
        \mu_p^{\smallsup{M}} \cn
    =\min \{1-\frac 12 \lambda^{-1} V^{-1/3},
    p\cn [1+\hat \Pi_{\SSS M}(0)]^+\},
        \en
where $x^+=\max\{x,0\}$.
Later we will see that in fact
$\mu_p^{\smallsup{M}} \cn =p\cn [1+\hat \Pi_{\SSS M}(0)]$.
We will prove that for all $M$ sufficiently
large (depending on $p$), $\lambda^3 \vee \beta$ sufficiently
small, and for all $p\leq p_c$,
        \eq
        \lbeq{P4P3}
        \max_{k\in \torus^*_{r,n}}
    \frac{\hat \tau_p(k)}{\hat{C}_{\mu_p^{\smallsup{M}}}(k)}
    \leq 3.
        \en
In fact, we will prove that the right hand side
of \refeq{P4P3} can be replaced by $1+c(\lambda^3 \vee \beta)$, where
$c$ is a universal constant.
The inequality \refeq{P4P3} is the key ingredient in the proof of
Theorem~\ref{thm-tc}.

The proof of \refeq{P4P3} is based on the following elementary lemma.
The lemma states that under an appropriate continuity
assumption, if an inequality implies a stronger inequality, then in
fact the stronger inequality must hold.
This kind of bootstrap argument has been applied repeatedly in analyses of the
lace expansion, and goes back to \cite{Slad87} in this context.

\begin{lemma}[The bootstrap]
\label{lem-P4}
Let $f$ be a continuous
function on the interval $[p_1,p_2]$, and assume that $f(p_1) \leq 3$.
Suppose for each $p \in (p_1,p_2)$ that if $f(p) \leq 4$
    then in fact $f(p) \leq 3$.
Then $f(p) \leq 3$ for all $p \in [p_1,p_2]$.
\end{lemma}

\proof
By hypothesis, $f(p)$ cannot be strictly between 3 and 4 for
any $p \in [p_1,p_2)$.
Since $f(p_1) \leq 3$, it follows by continuity that $f(p) \leq 3$
for all $p \in [p_1,p_2]$.
\qed

\smallskip
We will apply Lemma~\ref{lem-P4} with $p_1=0$, $p_2=p_c$, and
    \eq
    \lbeq{fdef}
    f(p) = \max\{f_1(p),f_2(p), f_3(p)\},
    \en
where
    \eq
    f_1(p)=p\cn,
    \quad
    f_2(p) =
    \max_{k\in \torus^*_{r,n}}
    \frac{\hat \tau_p(k)}{\hat{C}_{\mu_p^{\smallsup{M}}}(k)},
    \en
    \eq
    \lbeq{f3def}
    f_3(p) =
    \hspace{-4mm}
    \max_{\mbox{ {\scriptsize
    $\begin{array}{c}
                        k,l\in \torus^*_{r,n} \\ k \neq 0
                        \end{array} $ } }}
    \hspace{-4mm}
    \frac{\hat{C}_{1/\cn}(k)}{8} \frac{|\hat{\tau}_p(l) -
    \frac{1}{2} (\hat{\tau}_p(l-k) + \hat{\tau}_p(l+k))|}
    {\hat{C}_{\mu_p^{\smallsup{M}}}(l-k)\hat{C}_{\mu_p^{\smallsup{M}}}(l)
    +\hat{C}_{\mu_p^{\smallsup{M}}}(l)\hat{C}_{\mu_p^{\smallsup{M}}}(l+k)
    +\hat{C}_{\mu_p^{\smallsup{M}}}(l-k)\hat{C}_{\mu_p^{\smallsup{M}}}(l+k)}.
    \en


As we will see below in Section~\ref{sec-522}, the expression
$\hat{\tau}_p(l) -
    \frac{1}{2} (\hat{\tau}_p(l-k) + \hat{\tau}_p(l+k))$
can be interpreted
as $-\frac 12$ times a discrete Laplacian of $\hat{\tau}_p$.
In addition, this expression
is also the Fourier transform of $[1-\cos (k\cdot x)]\tau_p(x)$,
a quantity which appears implicitly in Proposition~\ref{prop-Pidiag} and in
the following bounds on $\Pi$, which play an essential
role in completing the bootstrap argument.
The proof of Proposition~\ref{cor-Pibds} is deferred to Section~\ref{sec-Pidbds}.

\begin{prop}
\label{cor-Pibds}
Let $M = 0,1,2,\ldots$, and assume that
Assumption~\ref{ass-rw} holds.
If $f(p)$ of \refeq{fdef} obeys $f(p)\leq K$, then there are positive
constants $c_K'$ and $\beta_0=\beta_0(K)$ such that
for $\lambda^3 \vee \beta \leq \beta_0$,
    \eq
    \sum_{x \in \torus_{r,n}}|\Pi_{\SSS M}(x)|
    \leq c_K' (\lambda^3 \vee \beta),
    \lbeq{diffPibda}
    \en
    \eq
    \sum_{x\in \torus_{r,n}} [1-\cos(k\cdot x)]|\Pi_{\SSS M}(x)|
    \leq c_K' (\lambda^3 \vee \beta) [1-\hat{D}(k)],
    \lbeq{diffPibd}
    \en
and for $M$ sufficiently large (depending on $K$ and $V$),
       \eq
        \sum_{x\in \torus_{r,n}} |R_{\SSS M} (x)|
        \leq (\lambda^3 \vee \beta),
        \lbeq{Rbda}
        \en
       \eq
        \sum_{x\in \torus_{r,n}} [1-\cos(k\cdot x)]|R_{\SSS M} (x)|
        \leq (\lambda^3 \vee \beta)[1-\hat{D}(k)].
        \lbeq{Rbd}
        \en
\end{prop}

\subsection{The bootstrap argument completed}
\label{sec-bootcomp}

We now show that $f$ of \refeq{fdef} obeys
the assumptions of Lemma~\ref{lem-P4}, with $p_1=0$ and $p_2=p_c$.

To see that $f(0) \leq 3$, we note that
$\hat{\tau}_0(k)=1$,
$\mu_0^{\smallsup{M}}=0$ and hence
$\hat{C}_{\mu_0^{\smallsup{M}}}(k)=1$, so that $f_2(0) = 1$.
Since $f_1(0)=f_3(0)=0$, we have $f(0) = 1 < 3$.

Next, we verify the continuity of $f$.  Continuity of $f_1$ is clear.
For $f_2$, since $\torus_{r,n}$ is finite it follows
that $\hat{\tau}_p(k)$ is a polynomial in $p$ and hence is
continuous.  Similarly, $\hat \Pi_{\SSS M}(0)$ is a polynomial in $p$.
Therefore
$\mu_p^{\smallsup{M}}$ is continuous in $p$, and hence
$\hat{C}_{\mu_{p}^{\smallsup{M}}}(k)$ also is, since $\hat{C}_\mu(k)$ is
continuous in $\mu$.
The numerator and denominator in the definition of $f_2$ are therefore
both continuous.
There is no division by zero, since the denominator is positive
when $\mu_p^\smallsup{M}<1$, by \refeq{Cdef}.  The maximum over $k$ is
a maximum over a finite set, so $f_2$ is continuous.
Similarly, $f_3$ is continuous, and thus
$f$ is continuous.

The remaining hypothesis of Lemma~\ref{lem-P4} is the substantial one,
and requires the detailed information about $\Pi_{\SSS M}$
and $R_{\SSS M}$ provided by Proposition~\ref{cor-Pibds}.
We fix $p < p_c$ and prove that $f(p) \leq 4$ implies $f(p) \leq 3$.
By the assumption that $f(p) \leq 4$, the hypotheses
of Lemma~\ref{prop-Pibds} are satisfied with $K=4$.
Therefore, assuming that
$M$ is sufficiently large and that $\lambda^3\vee \beta$ is
sufficiently small,
the bounds \refeq{diffPibda}--\refeq{Rbd} hold, with ${c}_K'$
replaced by ${c}_4'$.

Let
    \eq
    \lbeq{lambdadef}
    \lambda_p^{\smallsup{M}} \cn = p\cn[1+\hat \Pi_{\SSS M}(0)].
    \en
We now show that
$\lambda_p^{\smallsup{M}}\cn  \in [0, 1-\frac 12 \lambda^{-1} V^{-1/3}]$,
and hence $\mu_p^{\smallsup{M}}=\lambda_p^{\smallsup{M}}$.
By \refeq{tauMk} with $k=0$,
    \eq
    \chi(p)[1-\lambda_p^{\smallsup{M}}\cn ]
    =1+\hat \Pi_{\SSS M}(0)+\hat R_{\SSS M}(0).
    \en
Therefore,
    \eq
    1-\lambda_p^{\smallsup{M}}\cn
    \geq
    \chi^{-1}(p) \left[ 1 - |\hat\Pi_{\SSS M}(0)| - |\hat R_{\SSS M}(0)| \right]
    \geq  \chi^{-1}(p) \big[1-(c_4'+1)(\lambda^3 \vee \beta) \big].
    \en
Since $\chi(p)\leq \chi(p_c)= \lambda V^{1/3}$, for $\lambda^3 \vee
\beta$ sufficiently small it follows that
        \eq
    \lbeq{lambd}
        \lambda_p^{\smallsup{M}}\cn \leq 1-\frac12 \lambda^{-1} V^{-1/3}.
        \en
In addition, when $\lambda$ and $\beta$ are sufficiently small,
    \eq
    \lambda_p^{\smallsup{M}}\cn = p\cn[1+\hat{\Pi}_{\SSS M}(0)]
    \geq   p\cn\big[ 1- c_4' (\lambda^3 \vee \beta) \big] \geq 0.
    \en
This proves that $\mu_p^{\smallsup{M}}\cn =
\lambda_p^{\smallsup{M}}\cn  = p\cn[1+\hat \Pi_{\SSS M}(0)]$.

\subsubsection{The improved bounds on $f_1(p)$ and $f_2(p)$}

First, we improve the bound on $f_1(p)$.
We have already shown in \refeq{lambd} that
$\mu_p^{\smallsup{M}}\cn  \leq 1$.  Therefore,
by \refeq{diffPibda},
    \eq
    \lbeq{f1p3}
    f_1(p)=  p\cn =\frac{\mu_p^{\smallsup{M}}\cn }{1+\hat \Pi_{\SSS M}(0)}
    \leq \frac{1}{1- c_4' (\lambda^3 \vee \beta)}.
    \en
The right hand side is less than 3,
if $\lambda$ and $\beta$ are small enough.

To improve the bound on $f_2(p)$, we write \refeq{tauMk}
as $\hat{\tau} = \hat N/\hat F$, with
    \eq
    \lbeq{NFdef}
    \hat N(k)=1+\hat{\Pi}_{\sss M}(k)+\hat{R}_{\sss M}(k),
    \qquad \hat F(k)=1-p\cn \hat{D}(k)[1+\hat{\Pi}_{\sss M}(k)].
    \en
This yields
    \eqalign
    \lbeq{tauMkC2}
    \frac{\hat \tau_p(k)}{\hat C_{\mu_p^{\smallsup{M}}}(k)}
    & = \hat N(k)
    + \hat{\tau}_p(k)[1-\mu_p^{\smallsup{M}}\cn \hat{D}(k) -\hat F(k)]
    \nnb
    & =
    [1+\hat \Pi_{\SSS M}(k)+\hat R_{\SSS M}(k)]
    +\hat \tau_p(k)  p\cn\hat{D}(k)
    [\hat{\Pi}_{\SSS M}(k)-\hat{\Pi}_{\SSS M}(0)].
    \enalign
By Proposition~\ref{cor-Pibds}, and by our
assumptions that
$\hat \tau_p(k)\leq 4\hat C_{\mu_p^{\smallsup{M}}}(k)$
and $ p\cn\leq 4$, it follows from \refeq{tauMkC2} that
        \eq
        \lbeq{tauMkC3}
        \frac{\hat \tau_p(k)}{\hat C_{\mu_p^{\smallsup{M}}}(k)}
        \leq 1+ \Big( c_4' +1 +
        4^2 c_4' \hat C_{\mu_p^{\smallsup{M}}}(k)
        [1-\hat{D}(k)]\Big)
        (\lambda^3 \vee \beta).
        \en
Since
    \eq
    \lbeq{C1-D}
    0
    \leq
    \hat C_{\mu_p^{\smallsup{M}}}(k) [1-\hat{D}(k)]
    =
    1+
    \frac{\mu_p^{\smallsup{M}}\cn - 1}{1-\mu_p^{\smallsup{M}}\cn\hat{D}(k)}
    \hat{D}(k)
    \leq 2,
    \en
it follows from \refeq{tauMkC3} that
        \eq
        \lbeq{tauMkC4}
        f_2(p) =
        \max_{k \in \torus_{r,n}^*}
        \frac{\hat \tau_p(k)}{\hat C_{\mu_p^{\smallsup{M}}}(k)}
        \leq 1+ (c_4'+ 1 + 32 c_4')(\lambda^3 \vee \beta).
        \en
This is less than 3, if $\lambda^3 \vee \beta$ is small enough.

\subsubsection{Preliminaries for $f_3(p)$}
\label{sec-522}

Improving the bound on $f_3$ is more involved, and we first
develop some
useful preliminaries.

The expression $\hat{\tau}_p(l) -
    \frac{1}{2} (\hat{\tau}_p(l-k) + \hat{\tau}_p(l+k))$ in \refeq{f3def}
is closely related to a discrete second derivative of $\hat{\tau}_p(l)$.
In fact,
given a function $\hat{f}$ on $\torus_{r,n}^*$ and $k,l\in \torus_{r,n}^*$,
let
\eq
    \partial_k^+ \hat f(l) = \hat f(l+k)-\hat f(l),
\en
\eq
    \partial_k^- \hat f(l) = \hat f(l)-\hat f(l-k),
\en
and $\Delta_k \hat f(l) = \partial_k^- \partial_k^+ \hat f(l)$.
Then
\eq
\lbeq{Delf}
    -\frac 12 \Delta_k \hat f(l)
    =
    \hat f(l) - \frac{1}{2}(\hat{f}(l+k) + \hat{f}(l-k)).
\en
In particular, $-\frac 12 \Delta_k \hat{\tau}_p(l)$ appears in the
numerator of $f_3(p)$.

The following will be useful in computations involving $\Delta_k$.
Let $g$ be a symmetric function on the torus, meaning $g(x)=g(-x)$.
Then the Fourier transform of $g$ is actually the cosine series
$\hat g(l) = \sum_{x}g(x)  \cos(l\cdot x)$.
We define
    \eqalign
    \lbeq{hcos}
    \hat g^{\rm cos}(l,k)&=\sum_{x} g(x) \cos(l\cdot x)\cos(k\cdot x)
    =\frac12 [\hat g(l-k)+\hat g(l+k)],\\
    \lbeq{hsin}
    \hat g^{\rm sin}(l,k)&=\sum_{x} g(x) \sin(l\cdot x)\sin(k\cdot x)
    =\frac{1}{2} [\hat g(l-k)-\hat g(l+k)].
    \enalign
Then
\eq
    -\frac 12 \Delta_k \hat{g}(l) = \hat{g}(l)-\hat{g}^{\rm cos}(l,k).
\en
With this observation, the following lemma can be seen as a kind
of chain rule for the discrete differentiation of $\hat G$.

    \begin{lemma}
    \label{lem-fbd}
Suppose that $g(x)=g(-x)$, and let
$\hat{G}(k) = [1-\hat{g}(k)]^{-1}$.
    For all $k,l\in \torus^*$,
    \eqalign
    \lbeq{aimf}
    -\frac 12 \Delta_k G(l)
    &
    =
    \frac 12 [\hat{G}(l-k)+ \hat{G}(l+k)] \hat{G}(l)
   [\hat{g}(l)-\hat{g}^{\rm cos}(l,k)]
   \nnb & \quad
   - \hat{G}(l-k)\hat{G}(l) \hat{G}(l+k)\hat{g}^{\rm sin}(l,k)^2 .
    \enalign
    \end{lemma}

\proof
Let $\hat g_\pm = \hat{g}(l\pm k)$ and write $\hat g = \hat g(l)$.
Direct computation using \refeq{Delf} gives
    \eqalign
    -\frac 12 \Delta_k G(l)
    & =\frac 12 \hat{G}(l)\hat{G}(l+k)\hat{G}(l-k)\Big[
    [ 2\hat g - \hat g_+ - \hat g_- ] + [ 2 \hat g_+ \hat g_-
    - \hat g \hat g_- -  \hat g\hat g_+ ]\Big]\nonumber
    \\
    &  =
     \hat{G}(l)\hat{G}(l+k)\hat{G}(l-k)\Big[
    [\hat{g}(l)-\hat{g}^{\rm cos}(l,k)] +
    [ \hat g_+ \hat g_- - \hat{g}(l)\hat{g}^{\rm cos}(l,k) ]
    \Big],
    \lbeq{G1st}
    \enalign
using \refeq{hcos} in the last step.
By definition, and using the identity $\cos(u+v) = \cos u \cos v - \sin u \sin
v$,
    \eqalign
    \hat g_- \hat g_+
     & =
    \sum_{x,y} g(x)g(y)
    \cos((l+k)\cdot x)\cos((l-k)\cdot y)
     \nnb & =
     \hat{g}^{\rm cos}(l,k)^2-\hat{g}^{\rm sin}(l,k)^2.
     \lbeq{+-sc}
    \enalign
Substitution in \refeq{G1st} gives \refeq{aimf}.
\qed

Assume that $g(x) =g(-x)$.  Then
\eq
\lbeq{Dfinalpf}
    \frac 12 |\Delta_k \hat{g}(l)| = |\hat{g}(l)-\hat{g}^{\rm cos}(l,k)|
    \leq
    \sum_x [1-\cos (k\cdot x)]|g(x)|.
\en
Also, by the Cauchy-Schwarz inequality and the elementary
estimate $1-\cos^2t \leq 2[1-\cos t]$,
    \eqalign
    \hat{g}^{\rm sin}(k,l)^2
    &\leq \sum_x \sin^2(k\cdot x) |g(x)| \sum_y \sin^2(l\cdot y) |g(y)|
    \nnb
    &=\sum_x [1-\cos^2(k\cdot x)] |g(x)|
    \sum_y [1-\cos^2(l\cdot y)] |g(y)|\nonumber\\
    &
    \leq 4\sum_x [1-\cos(k\cdot x)] |g(x)|
    \sum_y [1-\cos(l\cdot y)] |g(y)|.
    \lbeq{Dsinbdx}
    \enalign
In addition,
\eqalign
    |\partial_k^\pm \hat g(l)|
    & \leq \sum_x  |{\rm Re}\{e^{il\cdot x}[e^{\pm ik\cdot x } -1]\} g(x)|
    \nnb
    & \leq \sum_x
   \big[ [1- \cos (k\cdot x)]
   + |\sin (k\cdot x)||\sin (l\cdot x)| \big]
    |g(x)|
    \nnb
\lbeq{crosspm}
    &\leq
    \sum_x
   [1- \cos (k\cdot x)] |g(x)|
   +
   \left\{
   4\sum_x [1-\cos(k\cdot x)] |g(x)|
    \sum_y [1-\cos(l\cdot y)] |g(y)|
    \right\}^{1/2},
\enalign
using the same technique as in \refeq{Dsinbdx} for the third
inequality.

The definition on $\partial_k^\pm$ leads to
the  quotient and product rules
\eqalign
    \partial_k^+ \frac{b(l)}{d(l)} &= \frac{\partial_k^+ b(l)}{d(l)}
    - \frac{b(l+k)\partial_k^+ d(l)}{d(l)d(l+k)},
\\
    \partial_k^- \frac{b(l)}{d(l)} &= \frac{\partial_k^- b(l)}{d(l)}
    - \frac{b(l-k)\partial_k^- d(l)}{d(l)d(l-k)},
\\
\lbeq{prod+}
    \partial_k^+ [\hat f(l)\hat h(l)] &= \partial_k^+ \hat f(l) \hat h(l+k)
    + \hat f(l) \partial_k^+ \hat h(l),
\\
\lbeq{prod-}
    \partial_k^- [\hat f(l)\hat h(l)] &= \partial_k^- \hat f(l) \hat h(l)
    + \hat f(l-k) \partial_k^- \hat h(l).
\enalign
This gives
\eqalign
    -\frac 12\Delta_k \frac{\hat b(l)}{\hat d(l)}
    & =
    -\frac 12 \partial_k^- \left\{
    \frac{\partial_k^+ \hat b(l)}{\hat d(l)} -
    \frac{\hat b(l+k)\partial_k^+ \hat d(l)}{\hat d(l) \hat d(l+k)}
    \right\}
    \nnb & =
    \frac{-\frac 12\Delta_k \hat b(l)}{\hat d(l)}
    +\frac 12\frac{\partial_k^+ \hat b(l-k) \partial_k^- \hat d(l)}
    {\hat d(l)\hat d(l-k)}
    +\frac 12 \frac{\partial_k^- \hat b(l+k) \partial_k^+ \hat d(l)}
    {\hat d(l)\hat d(l+k)}
    \nnb & \qquad
    +\frac 12 \frac{ \hat b(l) \Delta_k \hat d(l)}{\hat d(l)\hat d(l+k)}
    -\frac 12 \frac{ \hat b(l) \partial_k^+ \hat d(l-k)
        \partial_k^-[\hat d(l)\hat d(l+k)]}
        {\hat d(l-k)\hat d(l)^2\hat d(l+k)}.
\lbeq{Deltagx}
\enalign

\subsubsection{The improved bound on $f_3(p)$}

We now
improve the bound on $f_3(p)$.
We will write
\eq
    A = 1 +\mbox{const}(\lambda^3 \vee \beta),
\en
where the constant is universal and may change from line to line.

First, we recall the definitions of
$\hat{N}$ and $\hat{F}$ in \refeq{NFdef} and write $\hat \tau_p(l)$ as
\eq
    \hat \tau_p(l) = \frac{\hat N(l)}{\hat F(l)}
    =\frac{1}{1-\hat{g}(l)}
\en
with
\eq
\lbeq{gNF}
    \hat g(l)
    =
    1 - \frac{\hat{F}(l)}{\hat{N}(l)}
    =
    1 - \frac{1}{\hat{N}(l)}
    \Big\{
    1-\mu_p^\smallsup{M} \cn \hat D(l)
    + p\cn \hat{D}(l)[\hat\Pi_{\SSS M}(0) - \hat\Pi_{\SSS M}(l)]
    \Big\}.
\en
By Proposition~\ref{cor-Pibds},
\eq
\lbeq{Na}
    |\hat{N}(l)-1| \leq (c_4'+1) (\lambda^3 \vee \beta).
\en
In particular, $\hat{N}(l) >0$.  Since $\hat{\tau}_p(l)\geq 0$ (as proved
in \cite{AN84}), it follows that $\hat{F}(l)> 0$.
Proposition~\ref{cor-Pibds}, \refeq{f1p3} and \refeq{C1-D}
then imply that
\eqalign
\lbeq{Fa}
    0 \leq \hat{F}(l) &\leq
    [1-\mu_p^\smallsup{M}\cn \hat D(l)]
    + A c_4'(\lambda^3\vee \beta)
    [1- \hat D(l)]
    \nnb
    & \leq
    [1+ 2Ac_4'(\lambda^3\vee \beta)]
    [1-\mu_p^\smallsup{M}\cn \hat D(l)].
\enalign

By \refeq{tauMkC4}, \refeq{aimf} implies that
\eqalign
    & \hat{\tau}_p(l) - \frac{1}{2} (\hat{\tau}_p(l+k) + \hat{\tau}_p(l-k))
    =
    -\frac 12 \Delta_k \hat \tau_p(l)
    \nnb
    & \qquad
    \leq A\frac{1}{2}\big(
    \hat{C}_{\mu_p^\smallsup{M}}(l-k)
    + \hat{C}_{\mu_p^\smallsup{M}}(l+k) \big)
    \hat{C}_{\mu_p^\smallsup{M}}(l)|\hat{g}(l)-\hat{g}^{\rm cos}(l,k)|
    \nnb &\qquad
    \qquad
    + A
    \hat{C}_{\mu_p^\smallsup{M}}(l-k)
    \hat{C}_{\mu_p^\smallsup{M}}(l)\hat{C}_{\mu_p^\smallsup{M}}(l+k)
     \hat{g}^{\rm sin}(l,k)^2.
    \lbeq{aimg}
    \enalign
We will prove that
\eqalign
\lbeq{gineq3}
    |\hat{g}^{\rm sin}(l,k)|^2
    &\leq
    8A[1-\hat{D}(k)]\frac{1}{\hat{C}_{\mu_p^\smallsup{M}}(l)},
\\
\lbeq{gineq2}
    |\hat{g}(l)-\hat{g}^{\rm cos}(l,k)|
    &\leq
    A[1-\hat{D}(k)].
\enalign
These inequalities
imply that the right hand side of \refeq{aimg} is bounded above by
\eq
    A [1-\hat{D}(k)]
    \left[
    \frac 12 \big(
    \hat{C}_{\mu_p^\smallsup{M}}(l-k)
    +
    \hat{C}_{\mu_p^\smallsup{M}}(l+k)
    \big)
    \hat{C}_{\mu_p^\smallsup{M}}(l)
    + 8
    \hat{C}_{\mu_p^\smallsup{M}}(l-k)
    \hat{C}_{\mu_p^\smallsup{M}}(l+k)\Big]
    \right].
\en
Recalling that $\hat C_{1/\cn}(k)=[1-\hat D(k)]^{-1}$, this
gives
\eq
\lbeq{f3p3}
    f_3(p)\leq 1+\mbox{const}(\lambda^3 \vee \beta),
\en
so that in particular $f_3(p) \leq 3$.

To prove \refeq{gineq3}, we use \refeq{hsin} and \refeq{gNF} to see that
\eqalign
\lbeq{gsina}
    |\hat{g}^{\rm sin}(l,k)|
    & \leq \left|\frac{\hat{F}^{\rm sin}(l,k)}{\hat{N}(l-k)}\right|
    + \left|\frac{\hat{F}(l+k)\hat{N}^{\rm sin}(l,k)}
    {\hat{N}(l-k)\hat{N}(l+k)}\right|.
\enalign
By \refeq{Na}, the denominators are as close as desired to $1$.
To deal with the first term on the right side of \refeq{gsina}, we use
\refeq{NFdef} and \refeq{hsin} to obtain
\eq
\lbeq{gsinb}
     \hat{F}^{\rm sin}(l,k)
     =
     -p\cn \left[
     \hat{D}^{\rm sin}(l,k)[1+\hat{\Pi}_{\SSS M}(l-k)]
     + \hat{D}(l+k)\hat{\Pi}^{\rm sin}_{\SSS M}(l,k)
     \right].
\en
By \refeq{Dsinbdx},
    \eqalign
    |\hat{D}^{\rm sin}(k,l)|
    \leq \left\{ 4[1-\hat{D}(k)][1-\hat{D}(l)]\right\}^{1/2}.
    \lbeq{Dsinbd}
    \enalign
By \refeq{f1p3} and Proposition~\ref{cor-Pibds}, the first term on
the right hand side of \refeq{gsina} is at most
\eq
    A\left\{ 4[1-\hat{D}(k)][1-\hat{D}(l)]\right\}^{1/2}.
\en
The second term on the right hand
side of \refeq{gsina} can be bounded using the same method,
noting from \refeq{NFdef} that
the factor $\hat{F}(l+k)$ is at most $1+2\cdot1\cdot(1+1)=5$.
In addition, the factor $1-\hat{D}(l)$ can be bounded above
by $2\hat{C}_{\mu_p^\smallsup{M}}(l)^{-1}$, by \refeq{C1-D}.
Therefore, as required,
\eq
    \hat{g}^{\rm sin}(k,l)^2
    \leq
    8A
    [1-\hat{D}(k)]\hat{C}_{\mu_p^\smallsup{M}}(l)^{-1}.
\en

Finally, we estimate $\hat{g}(l)-\hat{g}^{\rm cos}(l,k) =
-\frac 12 \Delta_k \hat{g}(l)$ and prove \refeq{gineq2}.
By \refeq{gNF} and \refeq{Deltagx},
\eqalign
    -\frac 12\Delta_k \hat g(l)
    &=
     \frac{\frac 12\Delta_k \hat F(l)}{\hat N(l)}
    -\frac 12\frac{\partial_k^+ \hat F(l-k) \partial_k^- \hat N(l)}
    {\hat N(l)\hat N(l-k)}
    -\frac 12 \frac{\partial_k^- \hat F(l) \partial_k^+ \hat N(l)}
    {\hat N(l)\hat N(l+k)}
    \nnb & \qquad
    -\frac 12 \frac{ \hat F(l) \Delta_k \hat N(l)}{\hat N(l)\hat N(l+k)}
    +\frac 12 \frac{ \hat F(l) \partial_k^+ \hat N(l-k)
        \partial_k^-[\hat N(l)\hat N(l+k)]}
        {\hat N(l-k)\hat N(l)^2\hat N(l+k)}.
\lbeq{Deltag}
\enalign
The denominators are all as close to $1$ as desired, by \refeq{Na},
and we need to estimate the numerators.
The first term on the right side of \refeq{Deltag} is the main term.
By \refeq{prod+}--\refeq{prod-}, its numerator obeys
\eqalign
    \left| \frac 12\Delta_k \hat F(l) \right|
    & \leq
    p\cn \left|\frac 12 \Delta_k \hat{D}(l)\right|[1+\hat{\Pi}_{\SSS M}(l+k)]
    +\frac 12 p\cn |\partial_k^+
    \hat D(l-k) \partial_k^- \hat{\Pi}_{\SSS M}(l+k)|
    \nnb & \quad
    +\frac 12 p\cn |\partial_k^-
    \hat D(l) \partial_k^+ \hat{\Pi}_{\SSS M}(l)|
    +p\cn \left|\hat{D}(l-k) [\frac 12 \Delta_k  \hat{\Pi}_{\SSS M}(l)]\right|.
\lbeq{DeltaF}
\enalign
We bound the factors $p\cn $ by $A$.  The factor $|\frac 12\Delta_k \hat{D}(l)|$
is bounded above by $1-\hat{D}(k)$, by \refeq{Dfinalpf}.  The last term on
the right side of \refeq{DeltaF} is bounded by a small multiple of
$1-\hat{D}(k)$, by \refeq{Dfinalpf} and Proposition~\ref{cor-Pibds}.
For the cross terms, we use \refeq{crosspm} to obtain
\eqalign
    |\partial_k^\pm D(l)|
    & \leq
    [1-\hat D(k)] + 2 [1-\hat D(k)]^{1/2}[1-\hat D(l)]^{1/2}
    \nnb
    \lbeq{partialD}
    & \leq [1-\hat D(k)] + 2^{3/2} [1-\hat D(k)]^{1/2}.
\enalign
Applying Proposition~\ref{cor-Pibds}, similar
estimates apply to $\partial_k^\pm\hat{\Pi}_{\SSS M}$ and
$\partial_k^\pm\hat R_{\SSS M}$, but with an extra constant multiple of
$\lambda^3 \vee \beta$.  The two cross
terms in \refeq{DeltaF} are therefore bounded by a small multiple of
$1-\hat{D}(k)$.  We have shown that
the first term on the right side of \refeq{Deltag} is bounded above by
$A[1-\hat{D}(k)]$.

It is sufficient to show that the remaining
terms in \refeq{Deltag} are at most $[1-\hat{D}(k)]$ times a multiple
of $\lambda^3\vee \beta$.  The fourth term
on the right side of \refeq{Deltag} obeys this bound, using
\refeq{Fa} to bound $\hat{F}(l)$ by a constant,
and \refeq{Dfinalpf} and Proposition~\ref{cor-Pibds}
to bound $\Delta_k\hat{N}(l) = \Delta_k\hat{\Pi}_{\SSS M}(l)
+ \Delta_k\hat{R}_{\SSS M}(l)$ by $[1-\hat{D}(k)]$ times a
multiple of $\lambda^3\vee \beta$.

The remaining three terms in \refeq{Deltag} each contain a product
of a derivative of $\hat F$ with a derivative of $\hat N$, or a product of
two derivatives of $\hat N$ (using \refeq{prod-} for the last term).
Other factors of $\hat F$ or $\hat N$ are bounded by harmless constants.
The above arguments imply that $\partial_k^\pm \hat{N}(l)$ is bounded
by $\{[1-\hat D(k)] + 2^{3/2} [1-\hat D(k)]^{1/2}\}$ times a multiple
of $\lambda^3\vee \beta$,
as in \refeq{partialD} but with a small factor.  By the definition of $\hat{F}$
in \refeq{NFdef} and by the product rule \refeq{prod+},
we have
\eq
    \partial_k^+ \hat{F}(l)
    = -p\cn \partial_k^{+}\hat D(l)[1+\hat\Pi_{\SSS M}(l+k)]
    - p\cn \hat D(l) \partial_k^+ \hat \Pi_{\SSS M}(l),
\en
which is bounded by a multiple of the right side of \refeq{partialD}
(with no small factor).  The same bound is obeyed by
$\partial_k^- \hat{F}(l)$.  Although the derivative of $\hat F$ does not
produce a small factor, it is accompanied by a derivative of $\hat N$ which
does provide the desired factor $\lambda^3\vee \beta$.  Thus,
each of the remaining three
terms in \refeq{Deltag} is at most $[1-\hat{D}(k)]$ times a multiple
of $\lambda^3\vee \beta$.

This completes the proof that \refeq{Deltag} is bounded above by
$A[1-\hat{D}(k)]$.  Therefore, we have proved \refeq{f3p3}.
In particular, we have obtained the improved bound $f_3(p) \leq 3$.

\medskip
Throughout Section~\ref{sec-bootcomp},
we have relied on Proposition~\ref{cor-Pibds}.
We now prove this proposition.

\subsection{Proof of Proposition~\ref{cor-Pibds}}
\label{sec-Pidbds}

In this section, we prove Proposition~\ref{cor-Pibds}.
The main ingredient is the
following lemma.

\begin{lemma}[Bounds on the lace expansion]
\label{prop-Pibds}
Let $N = 0,1,2,\ldots$, and assume that
Assumption~\ref{ass-rw} holds.
For each $K>0$, there is a constant
$\bar{c}_K$ such that if $f(p)$ of \refeq{fdef} obeys $f(p)\leq K$, then
        \eq
        \sum_{x\in \torus_{r,n}}  {\Pi}^{\smallsup{N}}(x)
        \leq [\bar{c}_K (\lambda^3 \vee \beta)]^{N\vee 1}
        \lbeq{Pibd}
        \en
and
        \eq
        \sum_{x\in \torus_{r,n}}
    [1-\cos(k\cdot x)]
    {\Pi}^{\smallsup{N}}(x)
        \leq [1-\hat{D}(k)] [\bar{c}_K (\lambda^3 \vee \beta)]^{(N-1)\vee 1}.
        \lbeq{Pibd k}
        \en
\end{lemma}

Before proving Lemma~\ref{prop-Pibds}, we show that it implies
Proposition~\ref{cor-Pibds}.

\medskip \noindent
{\em Proof of Proposition~\ref{cor-Pibds}.}
The bounds \refeq{diffPibda}--\refeq{diffPibd} are immediate consequences
of Lemma~\ref{prop-Pibds}.  The constant $c_K'$ can be taken to be
equal to $4\bar c_K$, where the factor 4 comes from summing the geometric
series.

For the remainder term $R_{\SSS M}(x)$, we conclude from \refeq{Rxbd} that
    \eq
    \lbeq{Rxbdz}
    |R_{\SSS M}(x)|\leq 
    K \sum_{u,v}
    \Pi^{\smallsup{M}}(u)D(v-u)  \tau_p(x-v),
    \en
and hence \refeq{Rbda} is bounded above by $K \hat{\Pi}^{\smallsup{M}}(0)
\chi(p) \leq K \lambda V^{1/3}\hat{\Pi}^{\smallsup{M}}(0)$.
This can be made less than $\lambda^3\vee \beta$ by taking $M$
sufficiently large, by Lemma~\ref{prop-Pibds}.
For \refeq{Rbd}, we apply \refeq{normbd} with $J=3$ to obtain
       \eqalign
       &\sum_{x\in \torus_{r,n}} [1-\cos(k\cdot x)]|R_{\SSS M} (x)|\\
       &\quad \leq
       7 K [1-\hat{D}(k)]\hat{\Pi}_p^{\smallsup{M}}(0)\chi(p)
       + 7 K
       \big[\hat{\Pi}_p^{\smallsup{M}}(0)-\hat{\Pi}_p^{\smallsup{M}}(k)]
       \chi(p)
       + 7 K \hat{\Pi}_p^{\smallsup{M}}(0)
       [\hat{\tau}_p(0)-\hat{\tau}_p(k)].\nonumber
       \lbeq{Rbd2}
       \enalign
By Lemma~\ref{prop-Pibds}, we can choose $M$ large
enough that
$7K\hat{\Pi}_p^{\smallsup{M}}(0)\chi(p)\leq \frac 13(\lambda^3 \vee \beta)$.
The second term can be treated
similarly.
For the last term, we apply the bound $f_3(p)\leq K$ for $l=0$
and use \refeq{mudef} to see that
\eqalign
    |\hat{\tau}_p(0)-\hat{\tau}_p(k)|
    & = \frac 12 |\Delta_k \hat \tau_p(0)|
    \leq 24K[1-\hat D(k)][1-\mu_p^\smallsup{M}]^{-2}
    \leq 24K[1-\hat D(k)]4\lambda^2 V^{2/3}
    .
\enalign
Finally, we again take $M$ large and appeal to Lemma~\ref{prop-Pibds}.
\qed

Lemma~\ref{prop-Pibds} will follow from
Proposition~\ref{prop-Pidiag} combined with the following three lemmas.
For these three lemmas,
we recall the quantities defined in \refeq{T(x)def}--\refeq{Hpdef} and also
define
    \eq
    \lbeq{Tp2idef}
    T_p^{\smallsup{2}}
    = \frac 1V \sum_{k\in \torus^*_{r,n}} \hat{D}(k)^{2} \hat{\tau}_p(k)^3.
    \en

\begin{lemma}
\label{lem-TDa}
Fix $p \in (0,p_c)$, assume that $f(p)$ of \refeq{fdef} obeys
$ f(p) \leq K$, and assume that Assumption~\ref{ass-rw} holds.
There is a constant $c_K$, independent of $p$, such that
\eqalign
    &T_p^{\smallsup{2}} \leq c_K (\lambda^3 \vee \beta) ,
    \quad\quad
    T_p \leq  c_K(\lambda^3 \vee \beta),
    \qquad T_p' \leq 1+ c_K(\lambda^3 \vee \beta).
\enalign
The bound on $T_p^\smallsup{2}$ also applies if $\hat{\tau}_p(k)^3$ is
replaced by $\hat{\tau}_p(k)$ or $\hat{\tau}_p(k)^2$ in \refeq{Tp2idef}.
In addition, $\lambda^3$ can be replaced by $V^{-1}\chi(p)^3$ in each
of the above bounds.
\end{lemma}

\proof
We begin with $T_p^\smallsup{2}$.
We extract the term due to $k=0$ in \refeq{Tp2idef} and use $f_2(p) \leq K$
to obtain
        \eq
        \lbeq{T2z}
        T_p^{\smallsup{2}}
        \leq
        V^{-1} \chi(p)^3
        + V^{-1}
        \sum_{k\neq 0} \hat{D}(k)^2 K^3 \hat{C}_{\mu_p^{\smallsup{M}}}(k)^3.
        \en
The first term obeys $V^{-1} \chi(p)^3 \leq V^{-1}\chi(p_c)^3
= \lambda^3$, and the
desired result follows from \refeq{rwbd2} and \refeq{fdef}.
The conclusion concerning replacement of $\hat{\tau}_p(k)^3$
by $\hat{\tau}_p(k)$ or $\hat{\tau}_p(k)^2$ can be obtained
by going to $x$-space and using $\tau_p(x) \leq (\tau_p*\tau_p)(x) \leq
(\tau_p*\tau_p*\tau_p)(x)$.

For $T_p$, we extract the term in \refeq{T(x)def}
due to $y=z=0$ and $u=x$, which is
$p \cn D(x)\leq K\beta$, using $f_1(p) \leq K$
and \refeq{supbds}.  This gives
    \eq
    T_p(x) \leq K \beta +
    \shift\shift\sum_{u,y,z: (y,z-y, x+z-u)\neq (0,0,0)} \shift\shift
    \tau_p(y)\tau_p(z-y)KD(u)\tau_p(x+z-u).
    \en
Therefore, by \refeq{tauDtau},
    \eq
    T_p \leq K\beta  +
    3 K^2 \max_x \sum_{y,z\in \torus_{r,n}}
    \tau_p(y)(D*\tau_p)(z-y)(D*\tau_p)(x+z),
    \en
where the factor 3 comes from the 3 factors $\tau_p$ whose argument
can differ from 0. In terms of the Fourier transform, this gives
    \eq
    \lbeq{TpK}
    T_p \leq K \beta +
    3 K^2 \max_x V^{-1} \sum_{k\in \torus^*_{r,n}}
    \hat{D}(k)^2 \hat{\tau}_p(k)^3
    e^{-ik\cdot x}
    \leq K\beta + 3K^2 T_p^{\smallsup{2}}.
    \en
Our bound on $T_p^\smallsup{2}$ then gives the desired estimate for $T_p$.

The bound on $T_p'$ is a consequence of $T_p' \leq 1+3T_p$.
Here the term 1 is due to the contribution to \refeq{T'def} with
$y=z-y=x-z=0$, so that
$x=y=z=0$.   If at least one of $y, z-y, x-z$
is nonzero, then we can use \refeq{tauDtau}
for the corresponding two-point function.
\qed

\begin{lemma}
\label{lem-TDb}
Fix $p \in (0,p_c)$, assume that $f(p)$ of \refeq{fdef} obeys
$ f(p) \leq K$, and assume that Assumption~\ref{ass-rw} holds.
There is a constant $c_K$, independent of $p$, such that
\eqalign
    W_p(0;k) \leq c_K [1-\hat{D}(k)] (\lambda^3 \vee \beta),
    \quad\quad
    W_p(k) \leq c_K [1-\hat{D}(k)].
\enalign
\end{lemma}

\proof
For the bound on $W_p(0;k)$, we
use \refeq{tauDtau} to obtain
\eqalign
\lbeq{w1}
    \tilde{\tau}_p(x)
    & = p\cn D(x) + \sum_{v: v \neq x} p\cn D(v) \tau(x-v).
    \nnb
    & \leq p\cn D(x)+[p\cn]^2(D*D*\tau_p)(x).
\enalign
We insert \refeq{w1} into the definition \refeq{Wpydef} of $W_p(0;k)$ to get
\eq
\lbeq{w2}
    W_p(0;k) \leq
     p\cn \sum_x [1-\cos(k\cdot x)] D(x)\tau_p(x)
    + [p\cn]^2 \sum_{x} [1-\cos(k\cdot x)]\tau_p(x)(D*D*\tau_p)(x).
\en

We begin with the first term in \refeq{w2}, which receives no contribution
from $x=0$.  Using \refeq{tauDtau} and \refeq{w1} again, we obtain
    \eqalign
    &p\cn  \sum_{x\neq 0} [1-\cos(k\cdot x)] D(x)\tau_p(x)\nnb
    &\qquad\leq [p\cn]^2\sum_x [1-\cos(k\cdot x)] D(x)^2
    +[p\cn]^2 \sum_x [1-\cos(k\cdot x)] D(x)
    \sum_{v\neq x}D(v)\tau_p(x-v)\nnb
    &\qquad\leq [p\cn]^2\sum_x [1-\cos(k\cdot x)] D(x)^2
    +[p\cn]^3\sum_x [1-\cos(k\cdot x)] D(x) (D*D)(x)\nnb
    &\quad \qquad+[p\cn]^3\sum_x [1-\cos(k\cdot x)] D(x) (D*D*\tau_p)(x).
    \lbeq{secsplit}
    \enalign
The first term on the right side is bounded by $K^2\beta [1-\hat{D}(k)]$,
by \refeq{supbds}.
The second term can be bounded similarly,
using $\max_x (D*D)(x) \leq  \beta$.
For the last term in \refeq{secsplit},
we use Parseval's identity,
together with the fact that the Fourier transform
of $[1-\cos(k\cdot x)] D(x)$ is $\hat{D}(l)-\hat{D}^{\rm cos}(k,l)$,
to obtain
    \eq
    \sum_x [1-\cos(k\cdot x)] D(x) (D*D*\tau_p)(x)
    =
    \frac{1}{V} \sum_{l\in \torus^*_{r,n}}
    [\hat{D}(l)-\hat{D}^{\rm cos}(k,l)] \hat{D}(l)^2\hat{\tau}_p(l).
    \en
Applying \refeq{Dfinalpf} and the bound on $T^\smallsup{2}$
(with $\hat\tau_p(k)^3$ replaced by $\hat\tau_p(k)$), this is bounded by
    \eq
    [1-\hat{D}(k)]\frac{1}{V} \sum_{l\in \torus^*_{r,n}}
    \hat{D}(l)^2\hat{\tau}_p(l)\leq c_K (\lambda^3 \vee \beta) [1-\hat{D}(k)].
    \en
This completes the bound on the first term of \refeq{w2}.

For the second term in \refeq{w2}, we again use Parseval's identity to
obtain
    \eq
    \sum_{x} [1-\cos(k\cdot x)]\tau_p(x)(D*D*\tau_p)(x) =
    \frac{1}{V} \sum_{l\in \torus^*_{r,n}}
    \big[\hat{\tau}_p(l)-\frac12(\hat{\tau}_p(l+k)+\hat{\tau}_p(l-k))\big]
    \hat{D}(l)^2\hat{\tau}_p(l).
    \lbeq{mainterm}
    \en
By the assumed bounds on $f_2(p)$ and $f_3(p)$, this is at most
    \eqalign
    &8K^2  [1-\hat{D}(k)] \frac{1}{V} \sum_{l\in \torus^*_{r,n}}
     \hat{D}(l)^2
    \hat{C}_{\mu_p^{\smallsup{M}}}(l)\nnb
    &\qquad \quad \times
    \left[\hat{C}_{\mu_p^{\smallsup{M}}}(l-k)\hat{C}_{\mu_p^{\smallsup{M}}}(l)
    +
    \hat{C}_{\mu_p^{\smallsup{M}}}(l)\hat{C}_{\mu_p^{\smallsup{M}}}(l+k)
    +\hat{C}_{\mu_p^{\smallsup{M}}}(l-k)\hat{C}_{\mu_p^{\smallsup{M}}}(l+k)\right].
    \lbeq{mainterm2}
    \enalign
We set
    \eq
    C_{\mu,k}(x) = \cos(k\cdot x) C_{\mu}(x).
    \lbeq{Ckdef}
    \en
Then
    \eq
    |C_{\mu,k}(x)| \leq C_{\mu}(x),
    \lbeq{Cmukbd}
    \en
and, recalling \refeq{hcos},
    \eq
    \hat{C}_{\mu,k}(l)=\sum_{x\in \ver} \cos(k\cdot x) \cos(l\cdot x)C_{\mu}(x)
    =\hat{C}^{\rm cos}_\mu(l,k).
    \lbeq{Cktransf}
    \en
Also, by \refeq{+-sc},
    \eq
    \hat{C}_{\mu}(l-k)\hat{C}_{\mu}(l+k)
    =\hat{C}^{\rm cos}_{\mu}(l,k)^2 -\hat{C}^{\rm sin}_{\mu}(l,k)^2
    \leq \hat{C}^{\rm cos}_{\mu}(l,k)^2.
    \lbeq{C+-}
    \en
Therefore, using \refeq{Cmukbd} and Parseval's identity,
    \eqalign
    \frac{1}{V} \sum_{l\in \torus^*_{r,n}}
    \hat{D}(l)^2 \hat{C}_{\mu_p^{\smallsup{M}}}(l)
    \hat{C}_{\mu_p^{\smallsup{M}}}(l-k)\hat{C}_{\mu_p^{\smallsup{M}}}(l+k)
    &\leq \frac{1}{V} \sum_{l\in \torus^*_{r,n}}
    \hat{D}(l)^2 \hat{C}_{\mu_p^{\smallsup{M}}}(l)
    \hat{C}^{\rm cos}_{\mu_p^{\smallsup{M}}}(l,k)^2\nnb
    &= (D*D*C_{\mu_p^{\smallsup{M}}}*C_{\mu_p^{\smallsup{M}},k}
    *C_{\mu_p^{\smallsup{M}},k})(0)\nnb
    &\leq (D*D*C_{\mu_p^{\smallsup{M}}}*C_{\mu_p^{\smallsup{M}}}
    *C_{\mu_p^{\smallsup{M}}})(0).
    \enalign
Moreover, by \refeq{rwbd2},
    \eq
    (D*D*C_{\mu_p^{\smallsup{M}}}*C_{\mu_p^{\smallsup{M}}}
    *C_{\mu_p^{\smallsup{M}}})(0)=
    \frac{1}{V} \sum_{l\in \torus^*_{r,n}} \hat{D}(l)^2
    \hat{C}_{\mu_p^{\smallsup{M}}}(l)^3\leq 8\lambda^3 + \beta,
    \en
where the $\lambda^3$ arises from the $l=0$ term together with the fact
that $1 - \mu \geq \frac 12 \lambda^{-1} V^{-1/3}$.
This proves the desired bound on the last term in \refeq{mainterm2}.

To bound the sum of the remaining terms in \refeq{mainterm2}, we consider
    \eq
    \frac{1}{V} \sum_{l\in \torus^*_{r,n}} \hat{D}^2(l)
    \hat{C}_{\mu_p^{\smallsup{M}}}(l)^2
    \left[\hat{C}_{\mu_p^{\smallsup{M}}}(l-k)+
    \hat{C}_{\mu_p^{\smallsup{M}}}(l+k)\right].
    \lbeq{mixed term}
    \en
Applying \refeq{hcos}, \refeq{Cktransf}, \refeq{Cmukbd}, and
\refeq{rwbd2}, \refeq{mixed term} equals
    \eqalign
     \frac{2}{V} \sum_{l\in \torus^*_{r,n}} \hat{D}^2(l)
    \hat{C}_{\mu_p^{\smallsup{M}}}(l)^2
    \hat{C}^{\rm cos}_{\mu_p^{\smallsup{M}}}(l,k)
    &=2
    (D*D*C_{\mu_p^{\smallsup{M}}}*C_{\mu_p^{\smallsup{M}}}*C_{\mu_p^{\smallsup{M}},k})(0)\nnb
    &\leq 2(D*D*C_{\mu_p^{\smallsup{M}}}*C_{\mu_p^{\smallsup{M}}}*C_{\mu_p^{\smallsup{M}}})(0)
    \leq 2(8\lambda^3 + \beta) .
    \enalign
This completes the bound on the second term of \refeq{w2},
and thus the proof
that $W_p(0;k)\leq c_K (\lambda^3 \vee \beta) [1-\hat{D}(k)]$.

Finally, we estimate $W_p(k)$.
Note that no factor $\lambda^3 \vee \beta$ appears in the desired bound.
By \refeq{Wpydef}--\refeq{Wpdef},
    \eqalign
    W_p(k)
    &= p\cn \max_{y\in \ver} \sum_{x,v \in \ver}[1-\cos(k\cdot x)]
    D(v) \tau_p(x-v) \tau_p(x+y).
    \enalign
Let
    \eq
    \lbeq{Dkdef}
    D_k(x) = [1-\cos(k\cdot x)]D(x), \qquad
    \tau_{p,k}(x)= [1-\cos(k\cdot x)]\tau_p(x).
    \en
Applying \refeq{normbd} with
$t=k\cdot v + k\cdot (x-v)$, we obtain
    \eqalign
    W_p(k)
    &\leq 5p\cn \max_{y\in \ver}
    \sum_{x,v \in \ver}[1-\cos(k\cdot v)]D(v)\tau_p(x-v)
    \tau_p(y-x)\nnb
    &\qquad
    +5p\cn \max_{y\in \ver}\sum_{x,v \in \ver}
    D(v)[1-\cos(k\cdot (x-v))]\tau_p(x-v) \tau_p(y-x)
    \nnb
    \lbeq{Wpkbd}
    &\leq5K \max_{y\in \ver} (D_k * \tau_p*\tau_p)(y)
    + 5K\max_{y\in \ver}(D*\tau_{p,k}*\tau_p)(y).
    \enalign
For the first term, we have
    \eq
    (D_k * \tau_p*\tau_p)(y)
    =\frac 1{V}\sum_{l\in \torus^*_{r,n}} e^{-il\cdot y}
    \hat{D}_k(l) \hat{\tau}_p(l)^2
    \leq
    \frac {K^2}{V}\sum_{l\in \torus^*_{r,n}}
    |\hat{D}_k(l)| \hat{C}_{\mu_p^{\smallsup{M}}}(l)^2.
    \en
It follows from \refeq{Dfinalpf} that for all $k,l \in \torus^*_{r,n}$
    \eq
    |\hat{D}_k(l)| = |\hat{D}(l)-\hat{D}^{\rm cos}(k,l)| \leq [1-\hat{D}(k)],
    \en
and hence, by \refeq{rwbdi0},
    \eq
    \max_{y\in \ver} (D_k * \tau_p*\tau_p)(y) \leq [1-\hat{D}(k)] \frac {K^2}{V}
    \sum_{l\in \torus^*_{r,n}} \hat{C}_{\mu_p^{\smallsup{M}}}(l)^2
    \leq c_K(\lambda^3 \vee 1) [1-\hat{D}(k)],
    \en
where the $\lambda^3$ arises from the $l=0$ term.

The remaining term to estimate in \refeq{Wpkbd} is
    \eq
    \lbeq{Wrem}
    \max_{y\in \ver}(D*\tau_{p,k}*\tau_p)(y)= \max_{y\in \ver}
    \frac 1{V}\sum_{l\in \torus^*_{r,n}} e^{-il\cdot y}
    \hat{D}(l) \hat{\tau}_p(l)\hat{\tau}_{p,k}(l).
    \en
Since
    \eq
    \lbeq{taukl+-}
    \hat{\tau}_{p,k}(l)
    = \hat{\tau}_p(l)-
    \frac12(\hat{\tau}_p(l+k)+\hat{\tau}_p(l-k)),
    \en
we can use the bounds on $f_2(p)$
and $f_3(p)$ to see that \refeq{Wrem} is at most
    \eqalign
    & 8 K^2 [1-\hat{D}(k)]
    \frac 1{V}\sum_{l\in \torus^*_{r,n}} |\hat{D}(l)|
    \hat{C}_{\mu_p^{\smallsup{M}}}(l)
    \nnb
    &\quad \times
    \left[\hat{C}_{\mu_p^{\smallsup{M}}}(l-k)\hat{C}_{\mu_p^{\smallsup{M}}}(l)+
    \hat{C}_{\mu_p^{\smallsup{M}}}(l)\hat{C}_{\mu_p^{\smallsup{M}}}(l+k)
    +\hat{C}_{\mu_p^{\smallsup{M}}}(l-k)\hat{C}_{\mu_p^{\smallsup{M}}}(l+k)
    \right].
    \nonumber
    \enalign
The above sums can all be bounded using the methods employed for the
previous term.
For example, the last term can be estimated using $|\hat{D}(l)|\leq 1$,
\refeq{C+-}, \refeq{Cmukbd} and \refeq{rwbdi0}, by
    \eqalign
    \frac 1{V}\sum_{l\in \torus^*_{r,n}} \hat{C}_{\mu_p^{\smallsup{M}}}(l)
    \hat{C}_{\mu_p^{\smallsup{M}}}(l-k)\hat{C}_{\mu_p^{\smallsup{M}}}(l+k)
    &\leq
    \frac 1{V}\sum_{l\in \torus^*_{r,n}}
    \hat{C}_{\mu_p^{\smallsup{M}}}(l)
    \hat{C}_{\mu_p^{\smallsup{M}}}^{\rm \cos}(l,k)^2\nnb
    & =
    (C_{\mu_p^{\smallsup{M}}}*C_{\mu_p^{\smallsup{M}},k}*
    C_{\mu_p^{\smallsup{M}},k})(0)\nnb
    &\leq (C_{\mu_p^{\smallsup{M}}}*C_{\mu_p^{\smallsup{M}}}*
    C_{\mu_p^{\smallsup{M}}})(0).
    \lbeq{CSH}
    \enalign
\qed

    \begin{lemma}
    \label{lem-Mbd}
    Fix $p \in (0,p_c)$, assume that $f(p)$ of \refeq{fdef} obeys
    $f(p) \leq K$, and assume that Assumption~\ref{ass-rw} holds.
    There is a constant $c_K$, independent of $p$, such that
    \eq
    H_p(k)\leq c_K (\lambda^3\vee \beta) [1-\hat{D}(k)].
    \en
    \end{lemma}

\proof
Recall the definition of $H_p(a_1,a_2;k)$ in \refeq{Masdef}.
In terms of the Fourier transform, recalling \refeq{Dkdef},
        \eqalign
    H(a_1,a_2;k)
        & =\frac{1}{V^{3}} \sum_{l_1, l_2, l_3\in \torus^*_{r,n}}
        e^{-il_1\cdot a_1}
        e^{-il_2\cdot a_2}
        \hat{D}(l_1) \hat{\tau}_p(l_1)^2 \hat{D}(l_2) \hat{\tau}_p(l_2)^2
        \hat{\tau}_{p,k}(l_3)
        \nonumber \\
        \lbeq{Hbd1}
        &\hspace{30mm} \times\hat{\tau}_p(l_1-l_2)
        \hat{\tau}_p(l_2-l_3)\hat{\tau}_p(l_1-l_3).
        \enalign
We use $f(p) \leq K$ to replace $\hat{\tau}_p(k)$ by
$K\hat{C}_{\mu_p^{\smallsup{M}}}(k)$
and (recalling \refeq{taukl+-})
$\hat{\tau}_{p,k}(l_3)$ by
    \eqalign
    & 8K [1-\hat{D}(k)] \left[
    \hat{C}_{\mu_p^{\smallsup{M}}}(l_3-k)\hat{C}_{\mu_p^{\smallsup{M}}}(l_3)+
    \hat{C}_{\mu_p^{\smallsup{M}}}(l_3)\hat{C}_{\mu_p^{\smallsup{M}}}(l_3+k)
    +\hat{C}_{\mu_p^{\smallsup{M}}}(l_3-k)\hat{C}_{\mu_p^{\smallsup{M}}}(l_3+k)
    \right].
    \enalign
This gives an upper bound for \refeq{Hbd1} consisting
of a sum of 3 terms.

The last of these
terms can be bounded by
    \eqalign
    &8K^8 [1-\hat{D}(k)]\frac{1}{V^{3}} \sum_{l_1, l_2, l_3\in \torus^*_{r,n}}
        |\hat{D}(l_1)| \hat{C}_{\mu_p^{\smallsup{M}}}(l_1)^2 |\hat{D}(l_2)|
        \hat{C}_{\mu_p^{\smallsup{M}}}(l_2)^2
        \nnb
        &\qquad \times
        \hat{C}_{\mu_p^{\smallsup{M}}}(l_3-k)
        \hat{C}_{\mu_p^{\smallsup{M}}}(l_3+k)
        \hat{C}_{\mu_p^{\smallsup{M}}}(l_1-l_2)
        \hat{C}_{\mu_p^{\smallsup{M}}}(l_2-l_3)
        \hat{C}_{\mu_p^{\smallsup{M}}}(l_1-l_3).
    \lbeq{Hlast}
    \enalign
Using H\"older's inequality with $p=3$ and $q=3/2$,
\refeq{Hlast} is bounded
above by $8K^8$ times
        \eqalign
        &
        [1-\hat{D}(k)]
        \Big( \frac{1}{V^{3}} \sum_{l_1, l_2, l_3}
        |\hat{D}(l_1)|^{3/2} \hat{C}_{\mu_p^{\smallsup{M}}}(l_1)^3
        |\hat{D}(l_2)|^{3/2}
        \hat{C}_{\mu_p^{\smallsup{M}}}(l_2)^3
        \hat{C}_{\mu_p^{\smallsup{M}}}(l_3+k)^{3/2}
        \hat{C}_{\mu_p^{\smallsup{M}}}(l_1-l_3)^{3/2}
        \Big)^{2/3}
        \nonumber\\
        \lbeq{Hbd2}
        &\qquad \times\Big(\frac{1}{V^{3}} \sum_{l_1, l_2, l_3}
        \hat{C}_{\mu_p^{\smallsup{M}}}(l_1-l_2)^3
        \hat{C}_{\mu_p^{\smallsup{M}}}(l_2-l_3)^3
        \hat{C}_{\mu_p^{\smallsup{M}}}(l_3-k)^3
        \Big)^{1/3}.
        \enalign
Let
    \eq
    \lbeq{Salphadef}
    S^\smallsup{\alpha}_p = V^{-1}\sum_{l \in \torus^*_{r,n}}
    |\hat{D}(l)|^\alpha \hat{C}_{\mu_p^{\smallsup{M}}}(l)^3.
    \en
The Cauchy--Schwarz inequality implies that for all $k$ and $l_1$,
    \eq
    \frac 1V
    \sum_{l_3}
    \hat{C}_{\mu_p^{\smallsup{M}}}(l_3+k)^{3/2}
    \hat{C}_{\mu_p^{\smallsup{M}}}(l_1-l_3)^{3/2}
    \leq S_p^{\smallsup{0}}.
    \en
Therefore, \refeq{Hbd2} is bounded above by
        \eq
        \lbeq{Hbd3}
        [1-\hat{D}(k)]
        \big(S_p^{{\smallsup{0}}}\big)^{5/3}
        \big(S_p^{{\smallsup{3/2}}}\big)^{4/3}.
        \en
To complete the proof, we note that by H\"older's inequality,
        \eq
        S_p^{{\smallsup{3/2}}}\leq
        \big(S_p^{{\smallsup{2}}}\big)^{3/4}
        \big(S_p^{{\smallsup{0}}}\big)^{1/4}.
        \en
        Thus \refeq{Hbd3} is bounded above by
        $[1-\hat{D}(k)]S_p^\smallsup{2} (S_p^\smallsup{0})^2$.
        The latter factor can be bounded using \refeq{rwbdi0},
        and the former with \refeq{rwbd2}. This gives a bound of
    the desired form,
        with the $\lambda^3$ arising as usual from the $l=0$ term
    of $S_p^\smallsup{2}$.

Routine bounds can be used to deal with the other two terms
in a similar fashion.
    \qed

\medskip \noindent
{\em Proof of Lemma~\ref{prop-Pibds}.}
This is an immediate consequence of Proposition~\ref{prop-Pidiag} and
Lemmas~\ref{lem-TDa}--\ref{lem-Mbd}.
The bound \refeq{Pi1cos} is used for \refeq{Pibd k}
when $N=1$ (as \refeq{PiNbdx2} is not sufficient).
\qed

\subsection{The triangle condition}
\label{sec-tricon}

The hypotheses of Lemma~\ref{lem-P4} have all been verified,
and we conclude from the lemma that $f(p) \leq 3$ for all
$p \leq p_c$.
Moreover, we have seen
in \refeq{f1p3}, \refeq{tauMkC4} and \refeq{f3p3} that
it follows from $f(p) \leq 4$ that in fact
\eq
\lbeq{fK0}
    f(p) \leq K_0 = 1 +\mbox{const}(\lambda^3 \vee \beta),
\en
where the constant is universal.
Therefore, \refeq{fK0} indeed holds for all $p \leq p_c$.
In particular, the bounds of
Proposition~\ref{cor-Pibds} and
Lemmas~\ref{prop-Pibds}--\ref{lem-Mbd} all hold,
with $K$ equal to the $K_0$ of \refeq{fK0}.

\medskip \noindent
{\em Proof of Theorem~\ref{thm-tc}.}
By definition, $1 \leq \nabla_p(x,x) \leq T_p'$, and it was
noted in the proof of Lemma~\ref{lem-TDa} that $T_p' \leq 1+3T_p$.
For $x \neq y$ it follows from \refeq{tauDtau} that
    \eq
    \nabla_p(x,y) = (\tau_p * \tau_p *\tau_p)(y-x)
    \leq 3T_p(y-x)
    \quad\quad (x \neq y),
    \en
where the factor $3$ arises since there are three factors $\tau_p$
that could have a nonzero argument and hence permit application of
\refeq{tauDtau}.  Thus, it suffices to show that $3T_p \leq 10V^{-1}\chi(p)^3
+ 13\beta$.
By \refeq{T2z} and \refeq{TpK} with $K=K_0$ of \refeq{fK0},
\eq
    T_p \leq K_0\beta + 3K_0^2 \frac{\chi(p)^3}{V} + 3K_0^5\beta,
\en
and the desired result follows from the fact that $K_0$ can
be taken to be as close as desired to $1$ by taking $\lambda^3\vee \beta$
sufficiently small.
\qed

\section{Asymptotics for $\hat{\tau}_p(k)$}
\label{sec-pfstPibds}

In this section, we restrict attention to the torus $\torus_{r,n}$,
and assume that Assumption~\ref{ass-rw} holds with $\lambda^3 \vee \beta$
small.  We will show that it is possible to extend \refeq{chibd}
to an asymptotic formula for $\hat{\tau}_p(k)$ when $p \leq p_c$,
for all $k \in \torus^*_{r,n}$.  This result is not used
elsewhere in the paper.

The observation below \refeq{fK0} can be used in conjunction with
\refeq{Pibd} and \refeq{Rxbdz} to see
that $\lim_{M\to \infty} \sum_x |R_{\sss M}(0,x)|=0$, so
that by \refeq{tauMk} we have
        \eq
        \lbeq{tauinfk}
        \hat \tau_p(k) = \frac{1+\hat \Pi_{p}(k)}
        {1- p\cn \hat{D}(k)[1+\hat \Pi_{p}(k)]},
        \en
where $\Pi_p$ denotes $\Pi_{M=\infty}$.  Similarly, the limit $M \to \infty$
can be taken in \refeq{tauMkC2} to conclude that
        \eq
        \lbeq{tauMkC5}
        \frac{\hat \tau_p(k)}{\hat C_{\mu_p}(k)}-1 =
        \hat{\Pi}_p(k) + \hat \tau_p(k)
        p\cn \hat{D}(k)\big[\hat{\Pi}_p(k)-\hat{\Pi}_p(0)\big],
        \en
where
    $\mu_p = \mu_p^\smallsup{\infty}
    = p [1+\hat{\Pi}_p(0)]$.

\begin{theorem} [Asymptotics for the two-point function]
Suppose that Assumption~\ref{ass-rw} holds for percolation
on $\torus_{r,n}$, with $\lambda^3 \vee \beta$ sufficiently small.
For $p\leq p_c$,
\label{thm-main2}
    \eq\lbeq{taubd}
    \hat{\tau}_p(k)
    = (1+\bigo(\lambda^3\vee \beta))\hat{C}_{m_p}(k)
    = \frac{1+\bigo(\lambda^3 \vee \beta)}{1-m_p\cn \hat{D}(k)}
    ,
    \en
where
$ m_p\cn  = 1 - \cn (p_c-p) - \lambda^{-1}V^{-1/3}$ and
the error term is uniform in $k\in \torus^*_{r,n}$ and $p \leq p_c$.
\end{theorem}

\proof
Let $\epsilon = \cn(p_c-p)\geq 0$.
We first consider the case $k=0$.
The combination of \cite[Theorem~1.2~i)]{BCHSS04a} and Theorem~\ref{thm-tc}
implies that for all
$p \leq p_c$,
\eq
\lbeq{chibda}
    \frac{1}{\lambda^{-1}V^{-1/3}+\epsilon} \leq \chi(p)
    \leq
    \frac{1}{\lambda^{-1}V^{-1/3}+[1-\bigo(\lambda^3\vee \beta)]\epsilon}.
\en
This implies \refeq{taubd} for $k=0$, and
we therefore assume $k \neq 0$ henceforth.

Using \refeq{P4P3}, \refeq{f1p3}, Proposition~\ref{cor-Pibds} and
\refeq{C1-D}, \refeq{tauMkC5} leads to
        \eq
        \lbeq{tauMkC6}
        \big| \frac{\hat \tau_p(k)}{\hat C_{\mu_p}(k)}-1\big|
    = \bigo(\lambda^3 \vee \beta).
        \en
Since
\eq
    \frac{\hat \tau_p(k)}{\hat C_{m_p}(k)}-1 =
    \left( \frac{\hat \tau_p(k)}{\hat C_{\mu_p}(k)}-1 \right)
    \frac{\hat{C}_{\mu_p}(k)}{\hat{C}_{m_p}(k)}
    + \left( \frac{\hat{C}_{\mu_p}(k)}{\hat{C}_{m_p}(k)} - 1 \right)
\en
and since
\eq
    \big| \frac{\hat{C}_{\mu_p}(k)}{\hat{C}_{m_p}(k)} - 1 \big|
    =
    \frac{ |(\mu_p - m_p)\cn\hat{D}(k)|}{1-\mu_p \cn\hat{D}(k)}
    \leq
    \frac{|\mu_p - m_p|\cn}{1-\mu_p \cn},
\en
it suffices to show that
\eq
    \frac{|\mu_p - m_p|\cn}{1-\mu_p \cn} = \bigo(\lambda^3 \vee \beta).
\en

But by definition and \refeq{tauinfk},
\eq
    \mu_p\cn = 1 - [1+\hat{\Pi}_p(0)]\chi(p)^{-1}.
\en
Also, by \refeq{chibda},
\eq
    m_p \cn= 1 - [1+\bigo(\lambda^3 \vee \beta)]\chi(p)^{-1}.
\en
Therefore, as required,
\eq
    \frac{| \mu_p - m_p|\cn}{1-\mu_p \cn}
    = \frac{\bigo(\lambda^3 \vee \beta) \chi(p)^{-1}}
    {[1+\hat{\Pi}_p(0)]\chi(p)^{-1}}
    = \bigo(\lambda^3 \vee \beta).
\en
\qed

\section*{Acknowledgements}
This work began during a conversation at afternoon tea,
while RvdH, GS and JS were visiting Microsoft Research.
The work of GS was supported in part by NSERC of Canada.
The work of RvdH was carried out in part at the University
of British Columbia and in part at Delft University of Technology.
We thank Akira Sakai for helpful comments on
a previous version of the manuscript.


\end{document}

%% file: def99c.tex
\def\UseSection{
        \numberwithin{equation}{section}
    \theoremstyle{plain}
        \newtheorem{theorem}    {Theorem}[section]
        \DefineTheorems 
}

\def\DefineTheorems{
    
    \newtheorem{lemma}      [theorem] {Lemma}
    
    \newtheorem{prop}       [theorem] {Proposition}
    
    \newtheorem{cor}        [theorem] {Corollary}
    \newtheorem{ass}        [theorem] {Assumption}

    \theoremstyle{definition}
    \newtheorem{defn}       [theorem] {Definition}

    \theoremstyle{definition}

}

\newcommand{\bt}   {\begin{theorem}}
\newcommand{\et}   {\end  {theorem}}
\newcommand{\bl}   {\begin{lemma}}
\newcommand{\el}   {\end  {lemma}}
\newcommand{\bp}   {\begin{prop}}
\newcommand{\ep}   {\end  {prop}}
\newcommand{\bc}   {\begin{cor}}
\newcommand{\ec}   {\end  {cor}}
\newcommand{\bd}   {\begin{defn}}
\newcommand{\ed}   {\end  {defn}}

\newcommand{\ba}   {\begin{array}}
\newcommand{\ea}   {\end  {array}}
\newcommand{\be}   {\begin{enumerate}}
\newcommand{\ee}   {\end  {enumerate}}
\newcommand{\bi}   {\begin{itemize}}
\newcommand{\ei}   {\end  {itemize}}

\def\eq#1\en{\begin{equation}#1\end{equation}}  
\def\eqsplit#1\ensplit{
    \begin{equation}\begin{split}#1\end{split}\end{equation}
    }
\def\eqalign#1\enalign{
    \begin{align}#1\end{align}
    }
\def\eqmul#1\enmul{
    \begin{multline}#1\end{multline}
    }
\newcommand{\eqarrstar} {\begin{eqnarray*}} 
\newcommand{\enarrstar} {\end{eqnarray*}} 
\newcommand{\eqarray}   {\begin{eqnarray}} 
\newcommand{\enarray}   {\end{eqnarray}} 
\newcommand{\nnb}   {\nonumber \\} 

\newcommand{\lbeq}[1]  {\label{e:#1}}
\newcommand{\refeq}[1] {\eqref{e:#1}}    
\newcommand{\lbfg}[1]  {\label{fg: #1}}
\newcommand{\reffg}[1] {\ref{fg: #1}}

\makeatletter
\newcommand{\labelcounter}[2]{{%
    \stepcounter{#1}
    \protected@write\@auxout{}%
    {\string\newlabel{#2}{{\csname the#1\endcsname}{\thepage}}}%
    {\ref{#2}}
    }}
\makeatother
\newcommand{\sss}   { \scriptscriptstyle } 

\newcommand{\Ebold} {{\mathbb E}}

\newcommand{\Pbold} {{\mathbb P}}

\newcommand{\Rbold} {{\mathbb R}}

\newcommand{\Zbold} {{\mathbb Z}}

\newcommand{\Ccal}   {\mathcal{C}}

\newcommand{\spose}[1] {{\hbox to 0pt{#1\hss}} }
\newcommand{\ltapprox} {\mathrel{\spose{\lower 3pt\hbox{$\mathchar"218$}}
 \raise 2.0pt\hbox{$\mathchar"13C$}}}
\newcommand{\gtapprox} {\mathrel{\spose{\lower 3pt\hbox{$\mathchar"218$}}
 \raise 2.0pt\hbox{$\mathchar"13E$}}}

\newcommand{\nin}  {{ \not\in }}





%% file: ncube2fin.bbl
\begin{thebibliography}{10}

\bibitem{Aize97}
M.~Aizenman.
\newblock On the number of incipient spanning clusters.
\newblock {\em Nucl. Phys. B [FS]}, {\bf 485}:551--582, (1997).

\bibitem{AB87}
M.~Aizenman and D.J. Barsky.
\newblock Sharpness of the phase transition in percolation models.
\newblock {\em Commun. Math. Phys.}, {\bf 108}:489--526, (1987).

\bibitem{AN84}
M.~Aizenman and C.M. Newman.
\newblock Tree graph inequalities and critical behavior in percolation models.
\newblock {\em J. Stat. Phys.}, {\bf 36}:107--143, (1984).

\bibitem{AS00}
N.~Alon and J.H. Spencer.
\newblock {\em The Probabilistic Method}.
\newblock Wiley, New York, 2nd edition, (2000).

\bibitem{BA91}
D.J. Barsky and M.~Aizenman.
\newblock Percolation critical exponents under the triangle condition.
\newblock {\em Ann. Probab.}, {\bf 19}:1520--1536, (1991).

\bibitem{Boll01}
B.~Bollob\'as.
\newblock {\em Random Graphs}.
\newblock Cambridge University Press, Cambridge, 2nd edition, (2001).

\bibitem{BCHSS04a}
C.~Borgs, J.T. Chayes, R.~van~der Hofstad, G.~Slade, and J.~Spencer.
\newblock Random subgraphs of finite graphs: {I}. {The} scaling window under
  the triangle condition.
\newblock Preprint, (2003).

\bibitem{BCHSS04c}
C.~Borgs, J.T. Chayes, R.~van~der Hofstad, G.~Slade, and J.~Spencer.
\newblock Random subgraphs of finite graphs: {III}. {The} phase transition for
  the $n$-cube.
\newblock Preprint, (2003).

\bibitem{BS85}
D.C. Brydges and T.~Spencer.
\newblock Self-avoiding walk in 5 or more dimensions.
\newblock {\em Commun. Math. Phys.}, {\bf 97}:125--148, (1985).

\bibitem{Grim99}
G.~Grimmett.
\newblock {\em Percolation}.
\newblock Springer, Berlin, 2nd edition, (1999).

\bibitem{Hara00}
T.~Hara.
\newblock Critical two-point functions for nearest-neighbour high-dimensional
  self-avoiding walk and percolation.
\newblock In preparation.

\bibitem{HHS03}
T.~Hara, R.~van~der Hofstad, and G.~Slade.
\newblock Critical two-point functions and the lace expansion for spread-out
  high-dimensional percolation and related models.
\newblock {\em Ann. Probab.}, {\bf 31}:349--408, (2003).

\bibitem{HS90a}
T.~Hara and G.~Slade.
\newblock Mean-field critical behaviour for percolation in high dimensions.
\newblock {\em Commun. Math. Phys.}, {\bf 128}:333--391, (1990).

\bibitem{HS94}
T.~Hara and G.~Slade.
\newblock Mean-field behaviour and the lace expansion.
\newblock In G.\ Grimmett, editor, {\em Probability and Phase Transition},
  Dordrecht, (1994). Kluwer.

\bibitem{HS95}
T.~Hara and G.~Slade.
\newblock The self-avoiding-walk and percolation critical points in high
  dimensions.
\newblock {\em Combin. Probab. Comput.}, {\bf 4}:197--215, (1995).

\bibitem{HS00b}
T.~Hara and G.~Slade.
\newblock The scaling limit of the incipient infinite cluster in
  high-dimensional percolation. {II}. {Integrated} super-{Brownian} excursion.
\newblock {\em J.\ Math.\ Phys.}, {\bf 41}:1244--1293, (2000).

\bibitem{HS02}
R.~van~der Hofstad and G.~Slade.
\newblock A generalised inductive approach to the lace expansion.
\newblock {\em Probab. Th. Rel. Fields}, {\bf 122}:389--430, (2002).


\bibitem{HS03a}
R.~van~der Hofstad and G.~Slade.
\newblock The lace expansion on a tree with application to networks of
  self-avoiding walks.
\newblock {\em Adv.\ Appl.\ Math.}, {\bf 30}:471--528, (2003).

\bibitem{HS04a}
R.~van~der Hofstad and G.~Slade.
\newblock Expansion in $n^{-1}$ for percolation critical values on the $n$-cube
  and ${\mathbb Z}^n$: the first three terms.
\newblock Preprint, (2003).

\bibitem{HS04b}
R.~van~der Hofstad and G.~Slade.
\newblock Asymptotic expansions in $n^{-1}$ for percolation critical values on
  the $n$-cube and ${\mathbb Z}^n$.
\newblock Preprint, (2003).

\bibitem{JLR00}
S.~Janson, T.~{\L}uczak, and A.~Ruci\'{n}ski.
\newblock {\em Random Graphs}.
\newblock John Wiley and Sons, New York, (2000).

\bibitem{MS93}
N.~Madras and G.~Slade.
\newblock {\em The Self-Avoiding Walk}.
\newblock Birkh{\"a}user, Boston, (1993).

\bibitem{Mens86}
M.V. Menshikov.
\newblock Coincidence of critical points in percolation problems.
\newblock {\em Soviet Mathematics, Doklady}, {\bf 33}:856--859, (1986).

\bibitem{Nguy87}
B.G. Nguyen.
\newblock Gap exponents for percolation processes with triangle condition.
\newblock {\em J. Stat. Phys.}, {\bf 49}:235--243, (1987).

\bibitem{Slad87}
G.~Slade.
\newblock The diffusion of self-avoiding random walk in high dimensions.
\newblock {\em Commun. Math. Phys.}, {\bf 110}:661--683, (1987).

\bibitem{Slad99}
G.~Slade.
\newblock Lattice trees, percolation and super-{Brownian} motion.
\newblock In M.~Bramson and R.~Durrett, editors, {\em Perplexing Problems in
  Probability: Festschrift in Honor of Harry Kesten}, Basel, (1999).
  Birkh\"auser.

\end{thebibliography}
